%% file: main.tex
\documentclass[onefignum,onetabnum]{siamart220329}

\input{include}

\newcommand{\mvp}{matrix-vector product }

\begin{document}

\maketitle



    
    
\begin{abstract}
    We present two new \textbf{algebraic} multilevel hierarchical matrix algorithms to perform fast matrix-vector product (MVP) for $N$-body problems in $d$ dimensions, namely efficient \rk{$\mathcal{H}^2_{*}$} (fully nested algorithm, i.e., $\mathcal{H}^2$ matrix-like algorithm) and \rk{$\bkt{\mathcal{H}^2 + \mathcal{H}}_{*}$} (semi-nested algorithm, i.e., cross of $\mathcal{H}^2$ and $\mathcal{H}$ matrix-like algorithms). \rk{The efficient $\mathcal{H}^2_{*}$ and $\bkt{\mathcal{H}^2 + \mathcal{H}}_{*}$ hierarchical representations are based on our recently introduced \emph{weak admissibility} condition in higher dimensions \cite{khan2022numerical}, where the admissible clusters are the far-field and the \textbf{vertex-sharing} clusters}. Due to the use of nested form of the bases, the proposed hierarchical matrix algorithms are more efficient than the non-nested algorithms ($\mathcal{H}$ matrix algorithms). We rely on purely algebraic low-rank approximation techniques (e.g., ACA \cite{aca} and NCA \cite{bebendorf2012constructing, zhao2019fast, gujjula2022new}) and develop both algorithms in a black-box (kernel-independent) fashion. The initialization time of the proposed algorithms scales quasi-linearly \footnote{$\mathcal{O} \bkt{N\log^{\alpha}(N)}$, $\alpha \geq 0$ and small.}. Using the proposed hierarchical representations, one can perform the MVP that scales at most quasi-linearly. Another noteworthy contribution of this article is that we perform a comparative study of the proposed algorithms with different \textbf{algebraic} (NCA or ACA-based compression) fast MVP algorithms (e.g., $\mathcal{H}^2$, $\mathcal{H}$, etc) in $2$D and $3$D $(d=2,3)$. The fast algorithms are tested on various kernel matrices and applied to get fast iterative solutions of a dense linear system arising from the discretized integral equations and radial basis function interpolation. The article also discusses the scalability of the algorithms and provides various benchmarks. Notably, all the algorithms are developed in a similar fashion in $\texttt{C++}$ and tested within the same environment, allowing for \textbf{\emph{meaningful comparisons}}. The numerical results demonstrate that the proposed algorithms are competitive to the NCA-based standard $\mathcal{H}^2$ matrix algorithm \cite{bebendorf2012constructing, zhao2019fast, gujjula2022new} (where the admissible clusters are the far-field clusters) with respect to the memory and time. The $\texttt{C++}$ implementation of the proposed algorithms is available at \texttt{\url{https://github.com/riteshkhan/H2weak/}}.
\end{abstract}


\section{Introduction}
Kernel matrices are frequently encountered in many fields such as PDEs \cite{greengard1987fast, ho2013hierarchical, massei2022hierarchical}, Gaussian processes ~\cite{gp_book}, machine learning~\cite{gray2000, ker_den}, inverse problems \cite{saibaba2012application}, etc. These kernel matrices are usually large and dense. A direct evaluation of the product of a $N \times N$ kernel matrix with vector is prohibitive as its time and space complexity scale as $\mathcal{O} \bkt{N^2}$. However, these matrices possess block low-rank structures, which can be leveraged to store and perform matrix operations. The literature on the block low-rank matrices is vast, and we do not intend to review it here. Instead, we refer to a few selected articles \cite{chandrasekaran2007fast, borm2003introduction, halko2011finding, fmm_ref} and the books by Hackbusch \cite{hmatrix_book} and Bebendorf \cite{bebendorf08}.

In the past decades, various algorithms have been developed to perform the matrix-vector product efficiently. One of the first works in this area was the Barnes-Hut algorithm \cite{barnes1986hierarchical} or Tree code, which reduces the \mvp complexity from $\mathcal{O} \bkt{N^2}$ to $\mathcal{O} \bkt{N\log(N)}$. Greengard and Rokhlin propose the Fast Multipole Method (from now on abbreviated as FMM) \cite{greengard1987fast,greengard1997new}, which further reduces the \mvp cost to $\mathcal{O} \bkt{N}$. After that, various FMM-like kernel-independent algorithms \cite{ying2004kernel} were proposed, primarily based on analytic expansions. Hackbusch and collaborators \cite{borm2003introduction,borm2003hierarchical} are the pioneers in interpreting certain sub-blocks (sub-matrices) of the matrices arising from $N$-body problems as low-rank and provide an important theoretical framework. In \cite{borm2003introduction,borm2003hierarchical,grasedyck2003construction}, they discuss the standard (or strong) admissibility condition, i.e., where the separation distance between two clusters exceeds the diameter of either cluster. It is to be noted that the FMM (equivalent to $\mathcal{H}^2$ matrix with standard admissibility) and the Tree code (equivalent to $\mathcal{H}$ matrix with standard admissibility) are based on the standard or strong admissibility condition. In their subsequent work \cite{hackbusch2004hierarchical}, they introduce the notion of weak admissibility condition for one-dimensional problems. This article shows that the rank of interaction between the neighboring intervals in $1$D does not scale as a positive power of $N$, and consequently, the neighboring or non-overlapping intervals are admissible. HODLR~\cite{ambikasaran2013mathcal,ambikasaran2019hodlrlib}, HSS~\cite{xia2010fast,chandrasekaran2007fast} and HBS~\cite{gillman2012direct} matrices belong to the category of $\mathcal{H}$ matrix based on this weak admissibility condition. However, the work \cite{hackbusch2004hierarchical} does not discuss the notion of weak admissibility condition in higher dimensions. A straightforward extension of this idea, i.e., compressing all the non-self interactions, does not result in a quasi-linear \mvp algorithm in $d$ dimensions $(d>1)$ since the rank of the nearby clusters grows as $\mathcal{O} \bkt{N^{(d-1)/d} \log(N)}$ \cite{ho2013hierarchical,kandappan2022hodlr2d,khan2022numerical}. 

In our recent work \cite{khan2022numerical}, we have shown that the rank of the far-field and the vertex-sharing interactions do not scale with any positive power of $N$. Hence, admissibility of the far-field and the vertex-sharing clusters \emph{could be} a way to extend the notion of \emph{weak admissibility} condition in higher dimensions. Based on this \emph{weak admissibility} condition, we developed a \mvp algorithm in $d$ dimensions with quasi-linear complexity, namely the hierarchically off-diagonal low-rank matrix in $d$ dimensions or HODLR$d$D. \rk{HODLR$d$D is a special type of $\mathcal{H}$ matrix, and in this article, we denote it as $\mathcal{H}_{*}$ from now on. 

\emph{Our main focus in this article is to develop nested bases algorithms in a \textbf{purely} algebraic way, which are based on our weak admissibility condition (admissible clusters are the far-field and vertex-sharing clusters)}}.

\rk{In this article, we present two new hierarchical matrix algorithms with \textbf{nested bases} for fast MVP in $d$ dimensions. The first algorithm we propose is efficient $\mathcal{H}^2_{*}$, which is a \emph{special} subclass of the $\mathcal{H}^2$ matrices}. Due to the use of nested form of the bases, the time and space complexities are reduced compared to the $\mathcal{H}_{*}$ algorithm. The numerical experiments in $2$D and $3$D show that the efficient $\mathcal{H}^2_{*}$ algorithm is competitive to the NCA-based standard $\mathcal{H}^2$ matrix algorithm \cite{zhao2019fast, gujjula2022new} with respect to the MVP time and space complexity. \rk{The second algorithm we develop in this article is a \emph{semi-nested} \footnote{We term it semi-nested because nested bases are employed for the \textbf{partial} interaction list of a cluster, while non-nested bases are used for the remaining interaction list} algorithm, denoted as $\bkt{\mathcal{H}^2 + \mathcal{H}}_{*}$}. The $\bkt{\mathcal{H}^2 + \mathcal{H}}_{*}$ algorithm can be thought of as a cross of $\mathcal{H}^2$ and $\mathcal{H}$ like matrices. We develop all the algorithms in this article using purely algebraic techniques (ACA \cite{aca} or NCA \cite{zhao2019fast, gujjula2022new}), making them kernel-independent. 

\section{Related works and novelty of our work} The kernel matrices arising out of $N$-body problems usually possess block low-rank structure, i.e., certain sub-blocks of these matrices can be approximated by a low-rank representation. The construction of the low-rank approximation has been studied extensively. One can categorize the low-rank approximation techniques into two classes: (i) Analytic techniques (Analytic series expansions \cite{greengard1987fast,greengard1997new}, Interpolation \cite{fong2009black}, etc.) and (ii) Algebraic techniques (SVD, RRQR \cite{chandrasekaran1994rank}, RRLU \cite{pan2000existence}, ACA \cite{aca}, etc.). Both the analytic and algebraic techniques have their advantages and disadvantages. 

The popular FMM algorithms \cite{greengard1987fast, greengard1997new, ying2004kernel, fong2009black} and the widely used FMM libraries (FMM$3$D, bbFMM \cite{fong2009black}, PVFMM \cite{malhotra2015pvfmm}) are mainly based on analytic techniques. All these analytic FMM algorithms require very minimal pre-computation/initialization/assembly time \RK{for translation invariant and homogeneous kernels}. The majority of their time is spent on performing the MVP. Also, for the translation invariant and homogeneous kernels the storage requirement is less. However, in higher dimensions with a general kernel function, the storage requirement becomes high since the analytic compression rank is high. In this case, the translation matrices/operators need special attention (further compression using algebraic techniques like SVD) to reduce the storage, resulting in an increase in initialization time.

On the other hand, algebraic techniques are popular due to their wide range of applicability (in terms of domain and problem specification) and better compression than their \RK{purely} analytic counterparts. If the algebraic techniques are employed to construct the standard $\mathcal{H}^2$ matrix, the $\mathcal{H}^2$ matrix-vector product can be viewed as an algebraic variant of the FMM. Since the algebraic compression rank is usually lower than the \RK{purely} analytic compression rank, algebraic $\mathcal{H}^2$ MVP is faster than a \RK{purely} analytic FMM. However, the downside of algebraic techniques is that they usually have a longer $\mathcal{H}^2$ matrix initialization time than the analytic techniques for \RK{translation invariant and homogeneous kernels}. Nevertheless, the algebraic $\mathcal{H}^2$ matrix is more effective when performing multiple MVPs \RK{for a wide range of problems}.

We want to summarize the reasons for choosing algebraic techniques.
\begin{enumerate}
    \item Algebraic techniques only require access to the matrix entries and can be used in a kernel-independent (black-box) fashion; they neither need any series expansions of the underlying kernel function nor the knowledge of the kernel function. 
    \item The rank of the compressed blocks obtained using algebraic techniques is usually lower than that obtained from \RK{purely} analytic techniques, as algebraic techniques are domain and problem-specific.
    \item Algebraic techniques are more efficient than analytical techniques when dealing with data sets of higher dimensions.
    \item This article focuses on developing efficient nested hierarchical matrix representations based on a \emph{weak admissibility} condition \cite{khan2022numerical}, \RK{where the distance between two admissible clusters can be zero}. Algebraic techniques are useful in this case.
\end{enumerate}

 \textbf{The main focus of this article is to develop \emph{new} nested hierarchical MVP algorithms in a \emph{purely} algebraic way. Hence, we believe that it would be fair to compare the proposed algorithms with the related algebraic $\mathcal{H}^2$ MVP algorithms \cite{bebendorf2012constructing, zhao2019fast, gujjula2022new}.}
 
This section discusses some existing works that study different algebraic techniques and those that are in line with this article. 

The $\mathcal{H}$-matrix \cite{hackbusch2004hierarchical,hmatrix_book} is one of the most commonly used frameworks to handle the large dense matrices arising from the $N$-body problems. Though the $\mathcal{H}$ matrix algorithm is fast, it can be further accelerated using the nested form of the bases, which leads to $\mathcal{H}^2$ matrix algorithm \cite{borm2009construction,borm2003hierarchical}. It has been demonstrated that the $\mathcal{H}^2$ matrix is a more efficient framework than the $\mathcal{H}$ matrix regarding time and storage. 
Several kernel-independent algorithms have been developed to construct the $\mathcal{H}^2$ matrix efficiently. B\"{o}rm \cite{borm2009construction} proposes an $\mathcal{O} \bkt{N \log(N)}$ algorithm to construct $\mathcal{H}^2$ matrix in a purely algebraic way without storing the entire matrix. The projection-based algebraic Nested Cross Approximation (NCA) was introduced in \cite{bebendorf2012constructing}. The NCA is a variant of the Adaptive Cross Approximation (ACA) \cite{aca} that gives the nested basis for the $\mathcal{H}^2$ matrix. Bebendorf et al. \cite{bebendorf2012constructing} choose indices of points close to the Chebyshev grids and perform ACA upon them to obtain the pivots in a top-bottom fashion, leading to complexity $\mathcal{O} \bkt{N \log(N)}$. However, \cite{bebendorf2012constructing} required the geometrical information of clusters to find the pivots.
Zhao et al. \cite{zhao2019fast} develop two \textbf{purely} algebraic NCA algorithms. They first propose an $\mathcal{O} \bkt{N \log(N)}$ algorithm to find the pivots, which involves a similar approach (top-bottom tree traversal) like \cite{bebendorf2012constructing} but without the geometrical information of clusters. They also develop another algorithm that involves a bottom-top tree traversal to find the local pivots followed by a top-bottom tree traversal to find the global pivots. The overall complexity of this algorithm reduces to $\mathcal{O} \bkt{N}$. Gujjula et al. \cite{gujjula2022new} propose a NCA with complexity $\mathcal{O} \bkt{N}$, where they show that one can eliminate the global pivots selection step in \cite{zhao2019fast}, i.e., only the bottom-top traversal is sufficient (local pivots are enough), the subsequent top-bottom tree traversal could be redundant (which finds global pivots) without substantial change in the relative error. This makes their NCA slightly computationally faster than the \cite{zhao2019fast}. They also compare their NCA algorithm with \cite{bebendorf2012constructing,zhao2019fast} and present various benchmarks. We use similar NCAs as in \cite{zhao2019fast,gujjula2022new} to construct our algorithms, \rk{but directly applying them on our \emph{weak admissibility} condition in higher dimensions \cite{khan2022numerical} (far-field and vertex-sharing clusters are the admissible clusters) will not produce an \textbf{efficient} nested bases algorithm}.

The most important difference between the existing works and ours is that in $d$ dimensions ($d>1$), the existing algorithms convert the $\mathcal{H}$ matrix based on the strong/standard admissibility condition to the $\mathcal{H}^2$ matrix with the same strong/standard admissibility condition or directly construct the $\mathcal{H}^2$ matrix. \rk{In comparison, this article presents an efficient way to construct a hierarchical matrix algorithm with nested bases ($\mathcal{H}^2_{*}$), which is based on a \emph{weak admissibility} condition in higher dimensions \cite{khan2022numerical} (the admissible clusters include the vertex-sharing clusters and the distance between two clusters sharing a vertex is zero)}. It is by no means a trivial task to construct such a nested hierarchical representation efficiently. We have shown in \cite{khan2022numerical} that a \emph{weak admissibility}-based $\mathcal{H}_{*}$ representation is competitive with the strong admissibility-based standard $\mathcal{H}$ representation (ACA-based compression). This article will also investigate whether the proposed efficient $\mathcal{H}^2_{*}$ representation is competitive with the NCA-based standard $\mathcal{H}^2$ representation.
\rk{Furthermore, we present a semi-nested hierarchical representation, $\bkt{\mathcal{H}^2 + \mathcal{H}}_{*}$, which is also based on our \emph{weak admissibility} condition in higher dimensions \cite{khan2022numerical}}.

\vspace{1cm}
\boxed{\textbf{\emph{Main highlights of the article:}}}
The following are the main highlights of this article.
\begin{enumerate}
    \item \RK{We propose a fully nested \textbf{efficient} hierarchical matrix algorithm based on our \emph{weak admissibility} condition in higher dimensions \cite{khan2022numerical}. To construct it, we use a combination of NCA with bottom-top and top-bottom pivot selection with an appropriate partitioning of the interaction list. We refer to this algorithm as efficient $\mathcal{H}^2_*$ / $\mathcal{H}^2_*$(b$+$t)}. The numerical results show that the $\mathcal{H}^2_{*}$(b$+$t) algorithm performs much better than the $\mathcal{H}_{*}$ algorithm \cite{khan2022numerical} and is competitive to the NCA-based standard $\mathcal{H}^2$ matrix \cite{zhao2019fast} in $d$ dimensions with respect to the memory and MVP time.
    \item \rk{The second hierarchical matrix algorithm we propose is the $\bkt{\mathcal{H}^2 + \mathcal{H}}_{*}$, which is a \textbf{semi-nested} algorithm based on our \emph{weak admissibility} condition in higher dimensions \cite{khan2022numerical}}. The $\bkt{\mathcal{H}^2 + \mathcal{H}}_{*}$ could be considered as a cross of $\mathcal{H}^2$ and $\mathcal{H}$ matrix-like algorithms. We rely on the NCA and ACA to build the $\bkt{\mathcal{H}^2 + \mathcal{H}}_{*}$ representation. The numerical results show that the $\bkt{\mathcal{H}^2 + \mathcal{H}}_{*}$ algorithm performs better than the $\mathcal{H}_{*}$ algorithm, and in $3$D, it is competitive to the NCA-based standard $\mathcal{H}^2$ matrix \cite{zhao2019fast} with respect to the memory and MVP time.
    \item We report the performance of the proposed algorithms for various kernel matrix-vector products. In $2$D, we choose the single-layer Laplacian. In $3$D, we choose the single-layer Laplacian, Matérn covariance kernel and Helmholtz kernel (wave number $k=1$).
    \item \RK{We apply the proposed fast MVP algorithms to solve integral equations and radial basis function (RBF) interpolation in $2$D and $3$D using fast iterative solver. In each experiment, we compare the performance of the proposed algorithms with standard $\mathcal{H}^2$ matrix algorithm \cite{zhao2019fast, gujjula2022new} and discuss their scalability. Additionally, in $3$D, we report the performance of the proposed algorithms for various non-translation invariant kernels}. All the algorithms are implemented in the same fashion using $\texttt{C++}$ and tested within the same environment, allowing for meaningful comparisons.
    \item Finally, as a part of this article, we would also like to release the $\texttt{C++}$ implementation of the proposed algorithms made available at \texttt{\url{https://github.com/riteshkhan/H2weak/}}.
\end{enumerate}

\section*{\textit{Outline of the article}} The rest of the article is organized as follows. In \Cref{prelim}, we discuss different hierarchical representations used in this article with a summary of NCA. In \Cref{b2t_vs_t2b}, we discuss the applicability of NCA with various pivot selection strategies on our \emph{weak admissibility} condition. The \Cref{proposed_algo} discusses the proposed hierarchical matrix algorithms to perform fast MVP and their complexity analysis. In \Cref{num_results}, we perform various numerical experiments in $2$D and $3$D and compare the performance of the proposed algorithms with related purely algebraic hierarchical matrix algorithms. Finally, we conclude in \Cref{conclusion}.

\section{Preliminaries} \label{prelim} 
In this section, we describe the tree data structure used to build the different hierarchical matrix algorithms and briefly discuss various hierarchical representations used in this article. A brief discussion on Nested Cross Approximation (NCA) to construct $\mathcal{H}^2$ matrices is also included.
\subsection{Construction of the \texorpdfstring{$2^d$}{two power d} tree} \label{tree_construction} 
Let $C \subset \mathbb{R}^d$ be a compact hyper-cube, which contains the $N$ particles. To exploit the hierarchical representations, we need to subdivide the hyper-cube $C$ (computational domain). Depending upon the particle distribution or requirement, one can use different trees like $2^d$ uniform tree (or balanced tree), adaptive $2^d$ tree, level restricted $2^d$ tree, K-d tree, etc. \textbf{But for simplicity, we consider the $2^d$ uniform tree (for $d=1,2$ and $3$, it is binary, quad and oct tree, respectively) in this article}. At level $0$ of the tree is the hyper-cube $C$ itself (root level). A hyper-cube at a level $l$ is subdivided into $2^d$ finer hyper-cubes belonging to level $(l+1)$ of the tree. The former is called the parent of the latter, and the latter (finer hyper-cubes) are the children of the former (coarser hyper-cube). We define cluster $\mathcal{C}^{(l)}$ as containing all the particles inside a hyper-cube at level $l$. In \Cref{fig:tree_div}, we illustrate the uniform quad tree at different levels in $2$D. We stop the sub-division at a level $\kappa$ of the tree when each finest hyper-cube contains at most $n_{max}$ particles, where $n_{max}$ is a user-specified threshold that defines the maximum number of particles at each finest hyper-cube.  The total $2^{d \kappa}$ finest hyper-cubes at level $\kappa$ are called the leaves. Note that, $N \leq n_{max} 2^{d \kappa} \implies \kappa =  \ceil{\log_{2^d} \bkt{N/n_{max}}}$.

\begin{tikzfadingfrompicture}[name=rndpts]
    \fill[transparent!0] foreach ~ in {1,...,800}{(rand,rand) circle (0.01)};
\end{tikzfadingfrompicture}
\begin{tikzfadingfrompicture}[name=pivots]
    \fill[transparent!0] foreach ~ in {1,...,250}{(rand,rand) circle (0.02)};
\end{tikzfadingfrompicture}
\begin{figure}[H]
    \centering
    \subfloat[At level 0]{
    \begin{tikzpicture}[scale=2]
        \node[text=blue] at (0.5,0.5) {$0$};
        \draw[line width=0.2mm,  black] (0, 0) grid (1, 1);
        \draw (0,0) rectangle (1,1);
        \fill[black!20,scope fading=rndpts, fit fading=true] (0,0) rectangle (1,1);
        \node[text=blue] at (0.5,0.5) {$0$};
    \end{tikzpicture}
    } \qquad \qquad 
    \subfloat[At level 1]{
    \begin{tikzpicture}[scale=1]
        \node[text=blue] at (0.5,0.5) {$0$};
        \node[text=blue] at (1.5,0.5) {$1$};
        \node[text=blue] at (1.5,1.5) {$2$};
        \node[text=blue] at (0.5,1.5) {$3$};
        \draw[line width=0.2mm,  black] (0, 0) grid (2, 2);
        \draw (0,0) rectangle (2,2);
        \fill[black!20,scope fading=rndpts, fit fading=true] (0,0) rectangle (2,2);
        \draw[line width=0.2mm,  black] (0, 0) grid (2, 2);
        
        \node[text=blue] at (0.5,0.5) {$0$};
        \node[text=blue] at (1.5,0.5) {$1$};
        \node[text=blue] at (1.5,1.5) {$2$};
        \node[text=blue] at (0.5,1.5) {$3$};
    \end{tikzpicture}
    } \qquad \qquad 
    \subfloat[At level 2]{
    \begin{tikzpicture}[scale=0.5]
        \node[text=blue] at (0.5,0.5) {$0$};
        \node[text=blue] at (1.5,0.5) {$1$};
        \node[text=blue] at (1.5,1.5) {$2$};
        \node[text=blue] at (0.5,1.5) {$3$};

        \node[text=blue] at (2.5,0.5) {$4$};
        \node[text=blue] at (3.5,0.5) {$5$};
        \node[text=blue] at (3.5,1.5) {$6$};
        \node[text=blue] at (2.5,1.5) {$7$};
        
        \node[text=blue] at (2.5,2.5) {$8$};
        \node[text=blue] at (3.5,2.5) {$9$};
        \node[text=blue] at (3.5,3.5) {$10$};
        \node[text=blue] at (2.5,3.5) {$11$};

        \node[text=blue] at (0.5,3.5) {$15$};
        \node[text=blue] at (1.5,3.5) {$14$};
        \node[text=blue] at (1.5,2.5) {$13$};
        \node[text=blue] at (0.5,2.5) {$12$};
        
        \draw[line width=0.2mm,  black] (0, 0) grid (4, 4);
        \draw (0,0) rectangle (4,4);
        \fill[black!20,scope fading=rndpts, fit fading=true] (0,0) rectangle (4,4);
        \draw[line width=0.2mm,  black] (0, 0) grid (4, 4);
        
        \node[text=blue] at (0.5,0.5) {$0$};
        \node[text=blue] at (1.5,0.5) {$1$};
        \node[text=blue] at (1.5,1.5) {$2$};
        \node[text=blue] at (0.5,1.5) {$3$};

        \node[text=blue] at (2.5,0.5) {$4$};
        \node[text=blue] at (3.5,0.5) {$5$};
        \node[text=blue] at (3.5,1.5) {$6$};
        \node[text=blue] at (2.5,1.5) {$7$};
        
        \node[text=blue] at (2.5,2.5) {$8$};
        \node[text=blue] at (3.5,2.5) {$9$};
        \node[text=blue] at (3.5,3.5) {$10$};
        \node[text=blue] at (2.5,3.5) {$11$};

        \node[text=blue] at (0.5,3.5) {$15$};
        \node[text=blue] at (1.5,3.5) {$14$};
        \node[text=blue] at (1.5,2.5) {$13$};
        \node[text=blue] at (0.5,2.5) {$12$};
    \end{tikzpicture}
    } \qquad \qquad 
    \subfloat[At level 3]{
    \begin{tikzpicture}[scale=0.25]
        \draw[line width=0.2mm,  black] (0, 0) grid (8, 8);
        \draw (0,0) rectangle (8,8);
        \fill[black!20,scope fading=rndpts, fit fading=true] (0,0) rectangle (8,8);
        \draw[line width=0.2mm,  black] (0, 0) grid (8, 8);
    \end{tikzpicture}
    } \qquad 
    \caption{Hierarchical subdivision of the computational domain $C$ using $2^d$ uniform tree and the numbering convention followed in this article at different levels $(d=2)$.}
    \label{fig:tree_div}
\end{figure}

\subsection{Definition} 
\rk{The diameter of a cluster and the distance between two clusters are defined as follows:   
\begin{itemize}
    \item The diameter of a cluster $X$ is given by
        \begin{align}
            \text{diam}(X)= \sup\{\|u-v\|_{2}:u,v\in X\}
        \end{align}
    \item The distance between two clusters $X$ and $Y$ is given by
        \begin{align} \label{usual_distance}
            \text{dist}(X,Y) = \inf\{\|x-y\|_{2}:x\in X, y\in Y\}
        \end{align}
\end{itemize}}

\subsection{\texorpdfstring{$\mathcal{H}$}{H}-matrix based on the standard/strong admissibility condition} \label{strong_admiss} \rk{The} $\mathcal{H}$-matrix with strong admissibility condition is the hierarchical low-rank representation whose admissible clusters are the far-field or well-separated clusters (at least one box/cluster away). 

\boxed{\textit{\textbf{Standard/strong admissibility condition.}}}  \emph{Two clusters $X$ and $Y$ at the same level of the $2^d$ tree are admissible iff}
\begin{align} \label{H_admiss}
    \min(\text{diam}(X), \text{diam}(Y)) \leq \eta \text{ dist}(X,Y) \implies \dfrac{\min(\text{diam}(X), \text{diam}(Y))}{\text{dist}(X,Y)} \leq \eta 
\end{align}

\begin{figure}[H]
    \centering
    \subfloat[\scriptsize $\mathcal{H}_{\sqrt{d}}$ matrix at level $1$]{
        \includegraphics[scale=0.4]{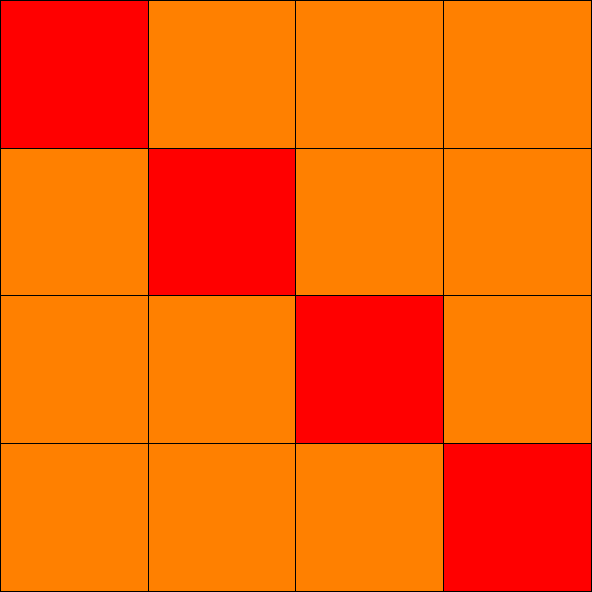}
        \label{fmm_mat_1}
        }\qquad \qquad 
    \subfloat[\scriptsize $\mathcal{H}_{\sqrt{d}}$ matrix at level $2$]{
        \includegraphics[scale=0.4]{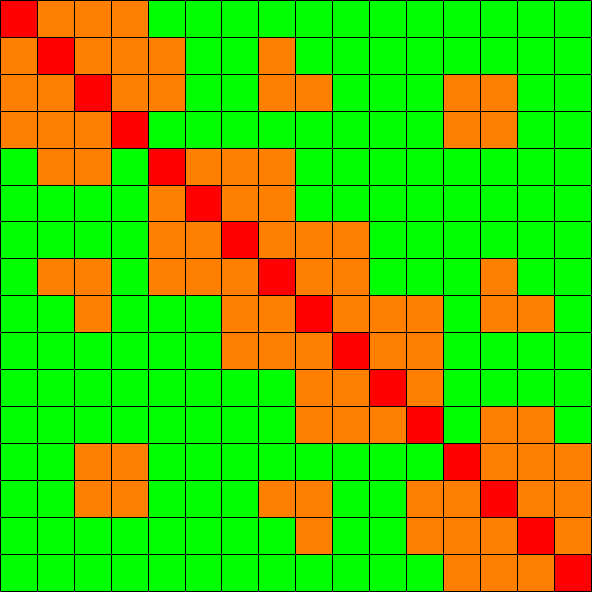}
        \label{fmm_mat_2}
        }\qquad \qquad 
    \subfloat[\scriptsize $\mathcal{H}_{\sqrt{d}}$ matrix at level $3$]{
        \includegraphics[scale=0.4]{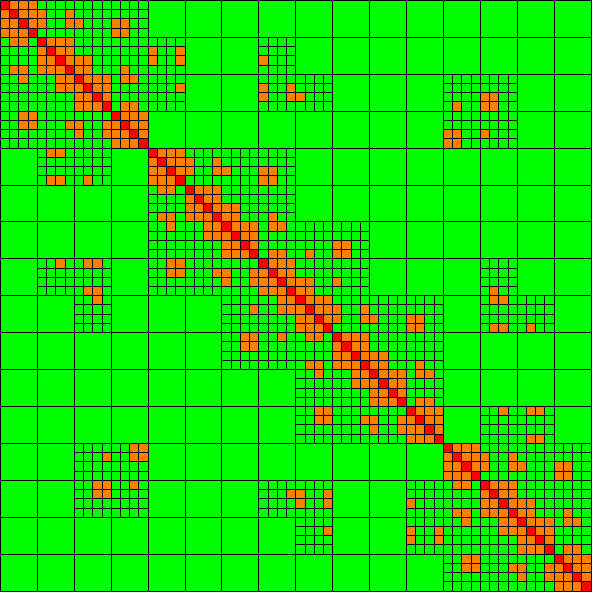}
        \label{fmm_mat_3}
        }\qquad \qquad 
    \subfloat{
        \begin{tikzpicture}
            [
            box/.style={rectangle,draw=black, minimum size=0.2cm},scale=0.2
            ]
            \node[box,fill=red,,font=\tiny,label=right:Self-interaction block,  anchor=west] at (-5,8){};
            \node[box,fill=orange,,font=\tiny,label=right:Neighbor block,  anchor=west] at (20,8){};
            \node[box,fill=green,,font=\tiny,label=right:Far-field compressed block,  anchor=west] at (40,8){};
        \end{tikzpicture}
    }
    \caption{$\mathcal{H}_{\sqrt{d}}$ in $2$D at the level $l$, where $l = 1, 2, 3$ (no compression at the level $1$).}
    \label{fig:FMM2D_matrix}
\end{figure}

 In this article, we set $\eta = \sqrt{d}$ in the admissibility condition \cref{H_admiss}, where $d$ is the underlying dimension. Specifically, $\eta = \sqrt{2}$ and $\sqrt{3}$ for $2$D and $3$D, respectively. \rk{\textbf{The $\mathcal{H}$ and $\mathcal{H}^2$ matrices discussed in this article are the standard $\mathcal{H}$ and $\mathcal{H}^2$ matrices, denoted as $\mathcal{H}_{\sqrt{d}}$ and $\mathcal{H}^2_{\sqrt{d}}$, respectively}. We use the number $``\sqrt{d}"$ in the subscript to indicate that the hierarchical representation is based on the standard/strong admissibility condition \cref{H_admiss} with $\eta = \sqrt{d}$. We denote the interaction list of a cluster $\mathcal{C}$ in $\mathcal{H}_{\sqrt{d}}$ and $\mathcal{H}^2_{\sqrt{d}}$ representations as $\mathcal{IL}_{\sqrt{d}}\bkt{\mathcal{C}}$, which contains the far-field clusters (at least one box away). For a cluster $\mathcal{C}$, the self cluster and its neighbors clusters collectively form the near-field list, denoted by $\mathcal{N}_{\sqrt{d}}\bkt{C}$ in $\mathcal{H}_{\sqrt{d}}$ and $\mathcal{H}^2_{\sqrt{d}}$ representations}. Note that the $\mathcal{H}^2_{\sqrt{d}}$ matrix is identical to the balanced FMM structure. The neighbors and the interaction list for $\eta = \sqrt{d}$ (in $2$D, $d=2$) at different levels of the tree are illustrated in \Cref{fig:FMM_interaction}. The $\mathcal{H}_{\sqrt{d}}$ matrix in $2$D at different levels is depicted in \Cref{fig:FMM2D_matrix}.
\begin{figure}[H]
\begin{center}
\subfloat[Level $1$ $(l=1)$]{
    \begin{tikzpicture}[scale=1.2]
        \fill [red] (1,0) rectangle (2,1);
        \fill [orange] (0,0) rectangle (1,1);
        \fill [orange] (1,1) rectangle (2,2);
        \fill [orange] (0,1) rectangle (1,2);
        \draw[step=1cm, black] (0, 0) grid (2, 2);
        \node[anchor=north] at (1.5,0.7) {\tiny $\mathcal{C}^{(1)}$};
    \end{tikzpicture}
    \label{fmm_interaction_1}
    }\qquad \qquad \qquad  \qquad 
    \subfloat[Level $2$ $(l=2)$]{
    \begin{tikzpicture}[scale=0.6]
        \draw[black] (0, 0) grid (4, 4);
        \fill [green] (0,0) rectangle (1,1);
        \fill [green] (0,1) rectangle (1,2);
        \fill [green] (0,2) rectangle (1,3);
        \fill [green] (0,3) rectangle (1,4);
        \fill [green] (1,3) rectangle (2,4);
        \fill [green] (2,3) rectangle (3,4);
        \fill [green] (3,3) rectangle (4,4);
        \fill [red] (2,1) rectangle (3,2);
        \fill [orange] (1,0) rectangle (2,1);
        \fill [orange] (1,1) rectangle (2,2);
        \fill [orange] (1,2) rectangle (2,3);
        \fill [orange] (2,2) rectangle (3,3);
        \fill [orange] (3,2) rectangle (4,3);
        \fill [orange] (2,0) rectangle (3,1);
        \fill [orange] (3,0) rectangle (4,1);
        \fill [orange] (3,1) rectangle (4,2);
        \node[anchor=north] at (2.5,1.8) {\tiny $\mathcal{C}^{(2)}$};
        \node[anchor=north] at (1.5,2.8) {};
        \node[anchor=north] at (1.5,1.8) {};
        \node[anchor=north] at (.5,1.8) {};
        \draw[black] (0, 0) grid (4, 4);
        \draw[line width=0.5mm,  black] (0, 0) rectangle (2, 2);
        \draw[line width=0.5mm,  black] (2, 2) rectangle (4, 4);
        \draw[line width=0.5mm,  black] (2, 0) rectangle (4, 2);
        \draw[line width=0.5mm,  black] (0, 2) rectangle (2, 4);
    \end{tikzpicture}
    \label{fmm_interaction_2}
    }\qquad \qquad \qquad  \qquad 
    \subfloat[Level $3$ $(l=3)$]{
    \begin{tikzpicture}[scale=0.6]
        \fill [green!50] (0,0) rectangle (1,1);
        \fill [green!50] (0,1) rectangle (1,2);
        \fill [green!50] (0,2) rectangle (1,3);
        \fill [green!50] (0,3) rectangle (1,4);
        \fill [green!50] (1,3) rectangle (2,4);
        \fill [green!50] (2,3) rectangle (3,4);
        \fill [green!50] (3,3) rectangle (4,4);
        \node[anchor=north] at (2.28,2) {\tiny $\mathcal{C}^{(3)}$};
        \node[anchor=north] at (1.5,2.8) {};
        \node[anchor=north] at (1.5,1.8) {};
        \node[anchor=north] at (.5,1.8) {};
        \draw[black] (0, 0) grid (4, 4);
        \draw[line width=0.5mm,  fill=green] (1, 0) rectangle (2, 1);
        \draw[line width=0.5mm,  fill=green] (1, 1) rectangle (2, 2);
        \draw[line width=0.5mm,  fill=green] (2, 2) rectangle (3, 3);
        \draw[line width=0.5mm,  fill=green] (3, 2) rectangle (4, 3);
        \draw[line width=0.5mm,  fill=green] (2, 0) rectangle (3, 1);
        \draw[line width=0.5mm,  fill=green] (3, 0) rectangle (4, 1);
        \draw[line width=0.5mm,  fill=green] (3, 1) rectangle (4, 2);
        \draw[line width=0.5mm,  fill=green] (1, 2) rectangle (2, 3);
        \fill [orange] (1.5,1) rectangle (3,2.5);
        \fill [red] (2,1.5) rectangle (2.5,2);
        \node[anchor=north] at (2.28,2.1) {\tiny $\mathcal{C}^{(3)}$};
        \draw[step=0.5cm, black] (1, 0) grid (4, 3);
        \draw[line width=0.5mm] (1, 0) rectangle (2, 1);
        \draw[line width=0.5mm] (1, 1) rectangle (2, 2);
        \draw[line width=0.5mm] (2, 2) rectangle (3, 3);
        \draw[line width=0.5mm] (3, 2) rectangle (4, 3);
        \draw[line width=0.5mm] (2, 0) rectangle (3, 1);
        \draw[line width=0.5mm] (3, 0) rectangle (4, 1);
        \draw[line width=0.5mm] (3, 1) rectangle (4, 2);
        \draw[line width=0.5mm] (1, 2) rectangle (2, 3);
    \end{tikzpicture}
    \label{fmm_interaction_3}
    }\qquad
    \subfloat{
        \begin{tikzpicture}
            [
            box/.style={rectangle,draw=black, minimum size=0.2cm},scale=0.2
            ]
            \node[box,fill=red,,font=\tiny,label=right:Cluster considered,  anchor=west] at (0,8){};
            \node[box,fill=orange,,font=\tiny,label=right:Neighbor cluster,  anchor=west] at (0,6){};
            \node[box,fill=green,,font=\tiny,label=right:Far-field cluster at current level,  anchor=west] at (40,8){};
            \node[box,fill=green!50,,font=\tiny,label=right:Far-field cluster at parent level,  anchor=west] at (40,6){};
        \end{tikzpicture}
    }
    \caption{The neighbors and the interaction list of cluster $\mathcal{C}^{(l)}$ $\bkt{\mathcal{IL}_{\sqrt{d}} \bkt{\mathcal{C}^{(l)}} \text{consists of the far-field clusters}}$ for $\eta = \sqrt{d}$ in \cref{H_admiss} at level $l$ (in $2$D, $d=2$). The superscript $``l"$ denotes the level of the uniform quad tree. The admissible clusters are the green colored clusters enclosed within a noticeable \textbf{black} border. The lighter shade colors represent the admissible clusters at the coarser levels (no admissible cluster at the level $1$).}
    \label{fig:FMM_interaction}
\end{center}
\end{figure}

\subsection{\texorpdfstring{$\mathcal{H}$}{H}-matrix based on Hackbusch's weak admissibility condition}
\rk{Hackbusch et al. \cite{hackbusch2004hierarchical} first introduced the notion of weak admissibility for one-dimensional problems, where they compress the sub-matrices corresponding to the adjacent intervals (clusters). This notion of weak admissibility doesn't look at the distance between the clusters; rather, it is defined as two non-overlapping clusters being admissible. To be more precise, in terms of matrix setting, this represents all the off-diagonal sub-matrices as low-rank. HODLR (non-nested) \cite{ambikasaran2013mathcal,ambikasaran2019hodlrlib}, HSS (nested)~\cite{xia2010fast,chandrasekaran2007fast} and HBS \cite{gillman2012direct} belong to the category of the hierarchical matrix based on Hackbusch's weak admissibility condition. 

The article \cite{hackbusch2004hierarchical} discusses the notion of weak admissibility only for one-dimensional problems. It does not discuss the notion of weak admissibility in higher dimensions ($d>1$) (though they mention in \cite{hackbusch2004hierarchical} that higher dimensional/multi-dimensional case will be considered later and refer to an article in the bibliography. However, to the best of our knowledge and searches, the article was never published or is available to the public).
A straightforward extension of the notion of weak admissibility discussed in~\cite{hackbusch2004hierarchical,hmatrix_book} to higher dimensions ($d>1$), i.e., compressing all the off-diagonal sub-matrices, will not result in a quasi-linear complexity hierarchical matrix algorithm. This is because, in higher dimensions, except for the vertex-sharing interaction sub-matrix, the rank of all other adjacent interactions sub-matrices grows with some positive power of $N$ \cite{ho2013hierarchical,kandappan2022hodlr2d,khan2022numerical}. Therefore, the notion of weak admissibility introduced by Hackbusch is fine for the one-dimensional problems but needs to be revised for the higher-dimensional problems}.

\subsection{\texorpdfstring{$\mathcal{H}$}{H}-matrix based on our \emph{weak admissibility} condition in higher dimensions} \label{hodlrdd_matrix}
\rk{The} article \cite{kandappan2022hodlr2d} by our group shows that for the $2$D Green's function, the rank of not just the far-field but also the \emph{vertex-sharing} interaction sub-matrices do not scale with any positive power of $N$ and introduces HODLR$2$D hierarchical representation based on it. In our recent article \cite{khan2022numerical}, we generalize this for any non-oscillatory translation invariant kernel function in $d$ dimensions and propose that the admissibility of \emph{far-field} and  \emph{vertex-sharing} clusters could be considered a way to extend the notion of \emph{weak admissibility} in higher dimensions. We also develop a hierarchical representation, namely HODLR$d$D \cite{khan2022numerical}, based on our proposed \emph{weak admissibility} condition. The HODLR$d$D is a non-nested algorithm, i.e., clusters' bases are non-nested with complexity $\mathcal{O} \bkt{pN \log(N)}$, where $p \in \mathcal{O} \bkt{\log(N) \log^d \bkt{\log(N)}}$. \rk{As mentioned earlier, we refer to HODLR$d$D as $\mathcal{H}_{*}$ in this article. We use the subscript $``*"$ to indicate that the hierarchical representation is based on our \emph{weak admissibility} condition in higher dimensions, i.e., the admissible clusters are \emph{far-field} and \emph{vertex-sharing} clusters}.

Let's discuss the interaction list of a cluster $\mathcal{C}^{(l)}$ at level $l$ of the $2^d$ tree, denoted as $\mathcal{IL}_{*} \bkt{\mathcal{C}^{(l)}}$. We consider the $2$D case for its simplicity in pictorial representation, but it is important to note that this article also considers the 3D case. We classify three different sets of clusters at the same level $l$: 
\\
(i) set of far-field/well-separated clusters (at least one cluster away from $\mathcal{C}^{(l)}$) denoted as $\mathcal{F} \bkt{\mathcal{C}^{(l)}}$, (ii) set of clusters share a vertex with $\mathcal{C}^{(l)}$ denoted as $\mathcal{V} \bkt{\mathcal{C}^{(l)}}$ and (iii) set of clusters share an edge with $\mathcal{C}^{(l)}$ denoted as $\mathcal{E} \bkt{\mathcal{C}^{(l)}}$. 

Let $child \bkt{\mathcal{C}^{(l)}}$, $siblings \bkt{\mathcal{C}^{(l)}}$ and $parent \bkt{\mathcal{C}^{(l)}}$ denote the child, siblings and parent of the cluster $\mathcal{C}^{(l)}$, respectively. We also define the clan set of the cluster $\mathcal{C}^{(l)}$ as follows:
\begin{align}
    clan \bkt{\mathcal{C}^{(l)}} = \Bigl \{ siblings \bkt{\mathcal{C}^{(l)}} \Bigr \} \bigcup \Bigl \{ child \bkt{P} : P \in \mathcal{E} \bkt{parent \bkt{\mathcal{C}^{(l)}}} \Bigr \}
\end{align}

For the \rk{$\mathcal{H}_{*}$ representation in $2$D}, the interaction list of the  cluster $\mathcal{C}^{(l)}$ is given by
\begin{align}
    \mathcal{IL}_{*} \bkt{\mathcal{C}^{(l)}} = clan \bkt{\mathcal{C}^{(l)}} \bigcap \bkt{\mathcal{V} \bkt{\mathcal{C}^{(l)}} \bigcup \mathcal{F} \bkt{\mathcal{C}^{(l)}}}
\end{align}
and the neighbors set of the  cluster $\mathcal{C}^{(l)}$ is given by $\mathcal{E} \bkt{\mathcal{C}^{(l)}}$.

The neighbors and the self cluster, i.e., the inadmissible clusters, collectively form the near-field list of $\mathcal{C}^{(l)}$, which is given by
\begin{align}
    \mathcal{N}_{*} \bkt{\mathcal{C}^{(l)}} = \mathcal{E} \bkt{\mathcal{C}^{(l)}} \bigcup \mathcal{C}^{(l)}
\end{align}

\begin{figure}[H]
    \centering
    \subfloat[\scriptsize $\mathcal{H}_{*}$ matrix at level $1$]{
        \includegraphics[scale=0.4]{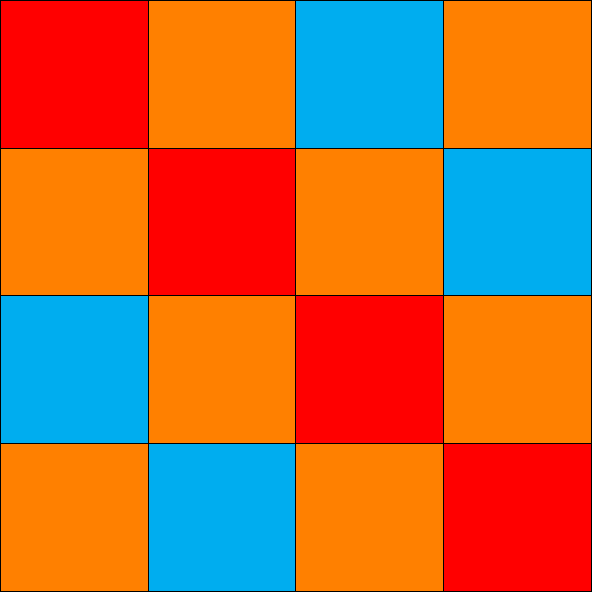}
        \label{hodlr2d_mat_1}
        }\qquad \qquad 
    \subfloat[\scriptsize $\mathcal{H}_{*}$ matrix at level $2$]{
        \includegraphics[scale=0.4]{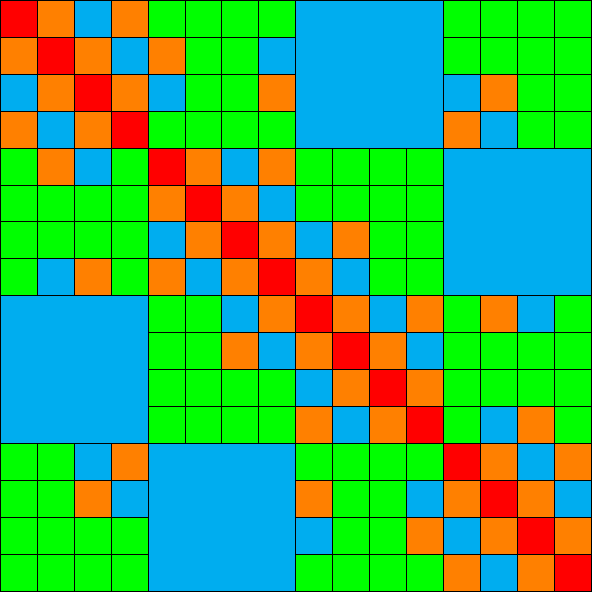}
        \label{hodlr2d_mat_2}
        }\qquad \qquad 
    \subfloat[\scriptsize $\mathcal{H}_{*}$ matrix at level $3$]{
        \includegraphics[scale=0.4]{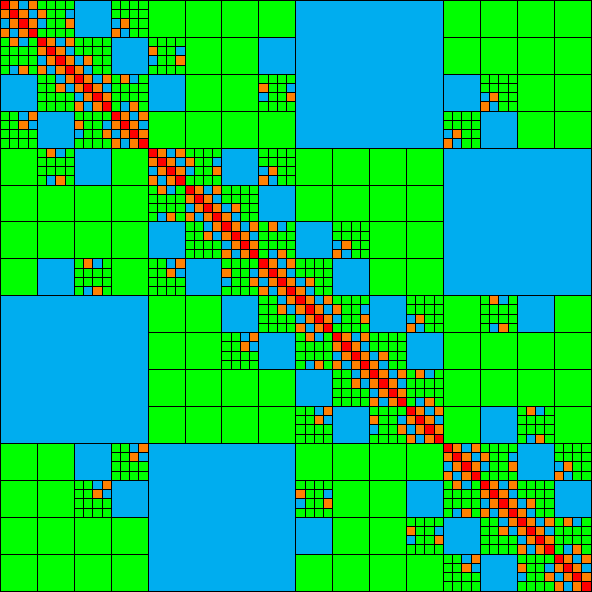}
        \label{hodlr2d_mat_3}
        }\qquad \qquad 
    \subfloat{
        \begin{tikzpicture}
            [
            box/.style={rectangle,draw=black, minimum size=0.2cm},scale=0.2
            ]
            \node[box,fill=red,,font=\tiny,label=right:Self-interaction block,  anchor=west] at (0,8){};
            \node[box,fill=orange,,font=\tiny,label=right:Neighbor block,  anchor=west] at (0,6){};
            \node[box,fill=cyan,,font=\tiny,label=right:Vertex-sharing compressed block,  anchor=west] at (30,8){};
            \node[box,fill=green,,font=\tiny,label=right:Far-field compressed block,  anchor=west] at (30,6){};
        \end{tikzpicture}
    }
    \caption{The $\mathcal{H}_{*}$ matrix in $2$D at the level $l$, where $l = 1,2,3$ (no far-field compression at the level $1$).}
    \label{fig:HODLRdD_matrix}
\end{figure}

The \Cref{fig:HODLRdD_interaction}, illustrates the neighbors and the interaction list at different levels of the tree for our \emph{weak admissibility} condition in $2$D. The $\mathcal{H}_{*}$ matrix in $2$D at different levels is depicted in \Cref{fig:HODLRdD_matrix}.

In article \cite{khan2022numerical}, we extend the notion of \emph{weak admissibility} condition to any dimension $d$ as follows: 

\boxed{\textit{\textbf{Weak admissibility condition in $d$ dimensions.}}} \label{weak_admis} \emph{Two clusters $X$ and $Y$ at the same level of the $2^d$ tree are admissible iff they share at the most a vertex, or equivalently, they do not share $d'$ hyper-surface ($d'>0$)}.

The interaction list and the near-field list for this \emph{weak admissibility} condition in $d$ dimensions are given in \Cref{tab:HOD_notation}.

\begin{figure}[H]
\begin{center}
\subfloat[Level $1$ $(l=1)$]{
    \begin{tikzpicture}[scale=1.2]
        \fill [red] (1,0) rectangle (2,1);
        \fill [orange] (0,0) rectangle (1,1);
        \fill [orange] (1,1) rectangle (2,2);
        \fill [cyan] (0,1) rectangle (1,2);
        \draw[step=1cm, black] (0, 0) grid (2, 2);
        \node[anchor=north] at (1.5,0.7) {\tiny $\mathcal{C}^{(1)}$};
        \draw[line width=0.5mm,  black] (0, 0) rectangle (2, 2);
    \end{tikzpicture}
    \label{hodlr2d_interaction_1}
    }\qquad \qquad \qquad \qquad 
    \subfloat[Level $2$ $(l=2)$]{
    \begin{tikzpicture}[scale=0.6]
        \draw[black] (0, 0) grid (4, 4);
        \fill [green] (0,0) rectangle (1,1);
        \fill [green] (0,1) rectangle (1,2);
        \fill [green] (2,3) rectangle (3,4);
        \fill [green] (3,3) rectangle (4,4);
        \fill [red] (2,1) rectangle (3,2);
        \fill [cyan] (1,0) rectangle (2,1);
        \fill [orange] (1,1) rectangle (2,2);
        \fill [orange] (2,2) rectangle (3,3);
        \fill [cyan] (3,2) rectangle (4,3);
        \fill [orange] (2,0) rectangle (3,1);
        \fill [cyan] (3,0) rectangle (4,1);
        \fill [orange] (3,1) rectangle (4,2);
        \node[anchor=north] at (2.5,1.8) {\tiny $\mathcal{C}^{(2)}$};
        \node[anchor=north] at (1.5,2.8) {};
        \node[anchor=north] at (1.5,1.8) {};
        \node[anchor=north] at (.5,1.8) {};
        \draw[black] (0, 0) grid (4, 4);
        \fill [cyan!50] (0,2) rectangle (2,4);
        \draw[line width=0.5mm,  black] (0, 0) rectangle (2, 2);
        \draw[line width=0.5mm,  black] (2, 2) rectangle (4, 4);
        \draw[line width=0.5mm,  black] (2, 0) rectangle (4, 2);
    \end{tikzpicture}
    \label{hodlr2d_interaction_2}
    }\qquad \qquad \qquad \qquad 
    \subfloat[Level $3$ $(l=3)$]{
    \begin{tikzpicture}[scale=0.6]
        \draw[black] (0, 0) grid (4, 4);
        \fill [green!50] (0,0) rectangle (1,1);
        \fill [green!50] (0,1) rectangle (1,2);
        \fill [green!50] (2,3) rectangle (3,4);
        \fill [green!50] (3,3) rectangle (4,4);
        \fill [cyan!50] (1,0) rectangle (2,1);
        \fill [cyan!50] (3,2) rectangle (4,3);
        \fill [cyan!50] (3,0) rectangle (4,1);
        \fill [red] (2,1) rectangle (3,2);
        \fill [green] (1,1) rectangle (2,2);
        \fill [orange] (1.5,1.5) rectangle (2,2);
        \fill [orange] (2.5,1.5) rectangle (3,2);
        \fill [cyan] (2.5,1) rectangle (3,1.5);
        \fill [orange] (2,1) rectangle (2.5,1.5);
        \fill [orange] (2,2) rectangle (2.5,2.5);
        \fill [green] (2,2) rectangle (3,3);
        \fill [green] (2,0) rectangle (3,1);
        \fill [green] (3,1) rectangle (4,2);
        \node[anchor=north] at (2.28,2.1) {\tiny $\mathcal{C}^{(3)}$};
        \node[anchor=north] at (1.5,2.8) {};
        \node[anchor=north] at (1.5,1.8) {};
        \node[anchor=north] at (.5,1.8) {};
        \draw[black] (0, 0) grid (4, 4);
        \draw[step=0.5cm, black] (3, 1) grid (4, 2);
        \draw[step=0.5cm, black] (2, 0) grid (3, 1);
        \draw[step=0.5cm, black] (2, 2) grid (3, 3);
        \draw[step=0.5cm, black] (2, 1) grid (3, 2);
        \draw[step=0.5cm, black] (1, 1) grid (2, 2);
        \fill [cyan!25] (0,2) rectangle (2,4);
        \fill [orange] (2,2) rectangle (2.5,2.5);
        \fill [cyan] (2.5,2) rectangle (3,2.5);
        \fill [cyan] (1.5,1) rectangle (2,1.5);
        \draw[line width=0.5mm,  black] (1, 1) rectangle (2, 2);
        \draw[line width=0.5mm,  black] (2, 2) rectangle (3, 3);
        \draw[line width=0.5mm,  black] (3, 1) rectangle (4, 2);
        \draw[line width=0.5mm,  black] (2, 0) rectangle (3, 1);
        \draw[black] (0, 2) rectangle (2, 4);
        \draw[step=1cm, black] (0, 0) grid (2, 2);
    \end{tikzpicture}
    \label{hodlr2d_interaction_3}
    }\qquad \qquad 
    \subfloat{
        \begin{tikzpicture}
            [
            box/.style={rectangle,draw=black, minimum size=0.2cm},scale=0.2
            ]
            \node[box,fill=red,,font=\tiny,label=right:Cluster considered,  anchor=west] at (0,8){};
            \node[box,fill=orange,,font=\tiny,label=right:Neighbor (edge-sharing) cluster,  anchor=west] at (0,6){};
            \node[box,fill=cyan,,font=\tiny,label=right:Vertex-sharing cluster at current level,  anchor=west] at (0,4){};
            \node[box,fill=cyan!50,,font=\tiny,label=right:Vertex-sharing cluster at parent level,  anchor=west] at (0,2){};
            \node[box,fill=cyan!25,,font=\tiny,label=right:Vertex-sharing cluster at grandparent level,  anchor=west] at (40,8){};
            \node[box,fill=green,,font=\tiny,label=right:Far-field cluster at current level,  anchor=west] at (40,6){};
            \node[box,fill=green!25,,font=\tiny,label=right:Far-field cluster at parent level,  anchor=west] at (40,4){};           
        \end{tikzpicture}
    }
    \caption{The neighbors and the interaction list of cluster $\mathcal{C}^{(l)}$ $\bkt{\mathcal{IL}_{*} \bkt{\mathcal{C}^{(l)}} \text{consists of vertex-sharing and far-field clusters}}$ for our weak admissibility condition in $2$D at level $l$. The superscript $``l"$ denotes the level of the uniform quad tree. The admissible clusters are the green and cyan colored clusters enclosed within a noticeable \textbf{black} border. The lighter shade colors represent the admissible clusters at the coarser levels.}
    \label{fig:HODLRdD_interaction}
\end{center}
\end{figure}

\rk{It is to be noted that the $\mathcal{H}_{*}$ representation does not satisfy the standard/strong admissibility condition \cref{H_admiss} because the distance \cref{usual_distance} between two vertex-sharing clusters (admissible clusters include vertex-sharing clusters) is zero as they share a common point. To be more precise, let $X$ and $Y$ be two clusters at the same level of the tree with $Y \in \mathcal{IL}_{*} \bkt{X}$. If $X$ and $Y$ share a vertex, then $\text{dist}(X,Y) = 0$. Given that the $\text{dist}(X,Y) = 0$, no real value of $\eta$ would generate the $\mathcal{H}_{*}$ representation from \cref{H_admiss}}. 

\textbf{Maximum size of interaction list $\bkt{\mathcal{IL}_{*} \bkt{\mathcal{C}}}$ and near-field list $\bkt{\mathcal{N}_{*} \bkt{\mathcal{C}}}$.} Let us find the maximum possible size of the interaction list and the near-field list (neighbors $+$ self) of a cluster $\mathcal{C}$ for the $\mathcal{H}_{*}$ representation. 

In $2$D, a square can have maximum $4$ edges. Therefore, the maximum size of the near-field list is 
\begin{align}
    \# \mathcal{N}_{*} \bkt{\mathcal{C}} = \bkt{ \# \text{edges} + \# \text{self} }= 4 + 1 =5
\end{align}
 Note that $\# \text{self}$ denotes the number of self-interaction or the cluster considered, which is $1$. \\
 The maximum size of the interaction list corresponding $\mathcal{C}$ is  
 \begin{align}
     \# \mathcal{IL}_{*} \bkt{\mathcal{C}} = \bkt{2^2 \times \# \mathcal{N}_{*} \bkt{\mathcal{C}}} - \# \mathcal{N}_{*} \bkt{\mathcal{C}} = (4 \times 5) - 5 = 15 
 \end{align}
 Out of a total of $15$ interactions, $3$ are vertex-sharing interactions, and $12$ are far-field interactions (\Cref{hodlr2d_interaction_3}). 

Similarly, in $3$D, a cube can have maximum $12$ edges and $6$ faces. Therefore, the maximum size of the near-field list is 
\begin{align}
    \# \mathcal{N}_{*} \bkt{\mathcal{C}} = \bkt{ \# \text{edges} + \# \text{faces} + \# \text{self} }= 12 + 6 + 1 = 19
\end{align}
The maximum size of the interaction list corresponding $\mathcal{C}$ is  
\begin{align}
  \# \mathcal{IL}_{*} \bkt{\mathcal{C}} = \bkt{2^3 \times \# \mathcal{N}_{*} \bkt{\mathcal{C}}} - \# \mathcal{N}_{*} \bkt{\mathcal{C}} = (8 \times 19) - 19 = 133   
\end{align}
Out of a total of $133$ interactions, $7$ are the vertex-sharing interactions, and $126$ are the far-field interactions. 

One can get a closed-form formula for the maximum possible size of the interaction list and the near-field list of a cluster in $d$ dimensions for the $\mathcal{H}_{*}$ representation. The maximum size of the near-field list is
\begin{align} \label{hodlrdd_nbd_list}
    \# \mathcal{N}_{*} \bkt{\mathcal{C}} = \bkt{3^d - 2^d}
\end{align}
and the maximum size of the interaction list is  
\begin{align} \label{hodlrdd_interaction_list}
    \# \mathcal{IL}_{*} \bkt{\mathcal{C}} = \bkt{2^d \times \# \mathcal{N}_{*} \bkt{\mathcal{C}}} - \# \mathcal{N}_{*} \bkt{\mathcal{C}} = \bkt{3^d - 2^d} \bkt{2^d -1} = 6^d - 4^d -3 ^d + 2^d
\end{align}
Out of a total of $(6^d - 4^d -3 ^d + 2^d)$ interactions, $(2^d -1)$ are the vertex-sharing interactions and remaining $(6^d - 4^d - 3^d +1)$ are the far-field interactions.

\begin{remark} \label{remark_3.1}
    The $\mathcal{H}_{*}$ matrix (\Cref{fig:HODLRdD_matrix}) has a couple of advantages compared to the $\mathcal{H}_{\sqrt{d}}$ matrix (\Cref{fig:FMM2D_matrix}), as given below.
    \begin{enumerate}
        \item In the $\mathcal{H}_{*}$ matrix, the maximum possible size of the interaction list corresponding to a cluster is $(6^d-4^d-3^d+2^d)$, whereas for the $\mathcal{H}_{\sqrt{d}}$ matrix, the maximum possible size of the interaction list is $(6^d-3^d)$. Therefore, the size of the interaction list of the $\mathcal{H}_{*}$ is slightly smaller than that of the $\mathcal{H}_{\sqrt{d}}$ matrix. This is visible from \Cref{hodlr2d_interaction_3} and \Cref{fmm_interaction_3}.
        \item $\mathcal{H}_{*}$ and $\mathcal{H}_{\sqrt{d}}$ matrices perform $(3^d -2^d)$ and $3^d$ dense matrix computations at leaf level, respectively. Therefore, $\mathcal{H}_{*}$ performs less dense matrix computations than the $\mathcal{H}_{\sqrt{d}}$ at the leaf level.
    \end{enumerate}
\end{remark}
\subsection{The Nested Cross Approximation (NCA) for \texorpdfstring{$\mathcal{H}^2$}{H2}-matrices}
The $\mathcal{H}$ matrix algorithm can be further accelerated using the nested form of the cluster bases, which leads to the $\mathcal{H}^2$ matrix algorithm. As we rely on the algebraic techniques in this article, this subsection briefly discusses algebraic ways called NCA to construct $\mathcal{H}^2$ matrix. Let $X$ and $Y$ be two clusters (hyper-cubes) which belong to the $2^d$ tree. Let $t^{X}$ and $s^{Y}$ be the index sets that store the indices of points $x (\in X)$ and $y (\in Y)$, respectively. 

\begin{align}
    t^X = \{i : x_i \in X \} \text{ and } s^Y = \{j : y_j \in Y \}
\end{align}
Let $K_{t^X, s^Y}$ be the matrix sub-block that captures the interaction between the clusters, whose $(i, j)^{th}$ entry is given by $K_{t^X, s^Y}(i,j) = K(t^X (i), s^Y (j))$.
If the cluster $Y$ is in the interaction list of the cluster $X$, i.e., $Y \in \mathcal{IL}_{\eta}\bkt{X}$ ($\mathcal{IL}_{\eta}\bkt{X}$ denotes the interaction list of $X$ generated from \cref{H_admiss} for a value $\eta >0$), then $K_{t^X, s^Y}$ is called admissible sub-block. The admissible sub-block $K_{t^X,s^Y}$ can be approximated by ACA \cite{aca,bebendorf2012constructing} as follows:
\begin{equation}
   K_{t^X, s^Y} \approx K_{t^X, s^Y}^{(p)} =  UV^* = K_{t^X,\sigma^Y} \bkt{K_{\tau^X,\sigma^Y}}^{-1} K_{\tau^X,s^Y}
\end{equation}
where $\tau^X \subset t^X$ and $\sigma^Y \subset s^Y$ are called the \emph{row} and \emph{column} pivots, respectively and the matrix $K_{t^X, s^Y}^{(p)}$ is the $p^{th}$ update of the original matrix $K_{t^X, s^Y}$. Here, we use the partially pivoted ACA algorithm \cite{aca}, where for a user-specified tolerance $\epsilon$, the iteration stops if the following condition is true
\begin{equation}
    \magn{u_p}_2 \magn{v_p}_2 \leq \epsilon \magn{K_{t^X, s^Y}^{(p)}}_F
\end{equation}
where $u_p$ and $v_p$ are the $p^{th}$ column vectors of the matrices $U$ and $V$, respectively. The ACA-based $\mathcal{H}$ matrix-like fast algorithm leads to a complexity of $\mathcal{O} \bkt{N \log(N)}$. One needs an algorithm that exploits the nestedness among the admissible clusters to achieve a lower complexity, so the NCA was proposed \cite{bebendorf2012constructing,zhao2019fast,gujjula2022new}. The validation of NCA, along with numerical error analysis, is presented in \cite{bebendorf2012constructing}. Let $\mathcal{A} \bkt{\mathcal{C}}$ denote the union of all the clusters forming the interaction list ($\mathcal{IL}_{\eta} \bkt{\mathcal{C}}$) at the same level and the clusters that form the interaction list with $\mathcal{C}$'s ancestors at the lower (coarser) level of the tree. The low-rank approximation via NCA takes the following form \cite{bebendorf2012constructing,zhao2019fast,gujjula2022new}

\begin{equation} \label{nca_equ}
    K_{t^X, s^Y} \approx \underbrace{K_{t^X, s^{X, i}} \bkt{K_{t^{X, i}, s^{X, i}}}^{-1}}_{U_X} \underbrace{K_{t^{X, i}, s^{Y, o}}}_{T_{X,Y}} \underbrace{\bkt{K_{t^{Y, o}, s^{Y, o}}}^{-1} K_{t^{Y, o}, s^Y}}_{V_Y^*}
\end{equation}
where $t^{X,i} \subset t^X$, $s^{X, i} \subset \mathcal{A}^{X,i}$ and $t^{Y,o} \subset \mathcal{A}^{Y,o}$, $s^{Y, o} \subset s^Y$. The sets $\mathcal{A}^{X,i}$ and $\mathcal{A}^{Y,o}$ are defined as
\begin{align}
    \mathcal{A}^{X,i} = \{ s^{X'} : X' \in \mathcal{A} \bkt{X} \} \text{ and } \mathcal{A}^{Y,o} = \{ t^{Y'} : Y' \in \mathcal{A} \bkt{Y} \} 
\end{align}

$t^{X,i}$ and $s^{X,i}$ are termed as \emph{incoming row pivots} and \emph{incoming column pivots} of $X$, respectively. $t^{Y,o}$ and $s^{Y,o}$ are termed as \emph{outgoing row pivots} and \emph{outgoing column pivots} of $Y$, respectively.
The $U_X$ and $V_Y^*$ in \Cref{nca_equ} are the column basis of $X$ and the row basis of $Y$, respectively. 

It is to be noted that if the matrix $K_{t^X, s^Y}$ is admissible, then so is $K_{t^Y, s^X}$ and the low-rank approximation via NCA takes the form
\begin{equation} 
    K_{t^Y, s^X} \approx \underbrace{K_{t^Y, s^{Y, i}} \bkt{K_{t^{Y, i}, s^{Y, i}}}^{-1}}_{U_Y} \underbrace{K_{t^{Y, i}, s^{X, o}}}_{T_{Y,X}} \underbrace{\bkt{K_{t^{X, o}, s^{X, o}}}^{-1} K_{t^{X, o}, s^X}}_{V_X^*}
\end{equation}
where $t^{Y,i} \subset t^Y$, $s^{Y, i} \subset \mathcal{A}^{Y,i}$ and $t^{X,o} \subset \mathcal{A}^{X,o}$, $s^{X, o} \subset s^X$.

Therefore, a cluster $X$ has four sets of pivots: (i) $t^{X,i} \subset t^X$, (ii) $s^{X,i} \subset \mathcal{A}^{X,i}$, (iii) $t^{X,o} \subset \mathcal{A}^{X,o}$ and (iv) $s^{X,o} \subset s^X$ with the bases $U_X$ and $V_X^*$.

\subsubsection{Various pivots selection strategies}
The efficiency (accuracy) of NCA depends on the strategy employed for selecting the pivots that represent the blocks at a specific level of the tree. 

Bebendorf et al. \cite{bebendorf2012constructing} first introduce projection-based algebraic NCA. They choose indices of points close to the predefined Chebyshev grids and perform ACA upon it to obtain the pivots, and this involves a top to bottom (T$2$B) traversal of the $2^d$ tree. The overall complexity to generate $\mathcal{H}^2$ matrix representation using this method is $\mathcal{O}\bkt{N \log(N)}$.

Zhao et al. \cite{zhao2019fast} avoid the geometrical projection to find the pivots and propose two purely algebraic NCAs based on different pivot selection strategies. They first discuss an $\mathcal{O} \bkt{N \log(N)}$ NCA algorithm to generate $\mathcal{H}^2$ matrix representation. In this algorithm, the pivots are obtained by traversing the $2^d$ tree from top to bottom (Algorithm 1, \cite{zhao2019fast}).\\
They also discuss an $\mathcal{O} \bkt{N}$ NCA algorithm, which is done through a two-step procedure. The first step involves a bottom to top (B$2$T) tree traversal to find the \emph{local pivots} (Algorithm 2, \cite{zhao2019fast}), and the second step involves a top to bottom (T$2$B) tree traversal to find the \emph{global pivots} (Algorithm 3, \cite{zhao2019fast}). The overall complexity of this algorithm is $\mathcal{O} \bkt{N}$.

Gujjula et al. \cite{gujjula2022new} propose a purely algebraic NCA of complexity $\mathcal{O} \bkt{N}$, which is similar to the Algorithm 2 of \cite{zhao2019fast}. It also involves a bottom to top (B$2$T) tree traversal to find the \emph{local pivots}. But after that, it does not traverse the tree top to bottom (T$2$B) like the $\mathcal{O} \bkt{N}$ NCA algorithm of \cite{zhao2019fast} for the \emph{global pivots}. They construct the operators based on the \emph{local pivots}. They show that to construct the standard $\mathcal{H}^2$ matrix representation \emph{local pivots} are enough. To be more precise, \cite{gujjula2022new} shows that for a purely algebraic construction of the standard $\mathcal{H}^2$ matrix, the bottom to top (B$2$T) tree traversal is sufficient; the subsequent top to bottom (T$2$B) tree traversal could be redundant without substantial change in the relative error. Therefore, the single bottom to top (B$2$T) tree traversal is sufficient, i.e., \emph{local pivots} are sufficient for a good approximation, and the second step of $\mathcal{O} \bkt{N}$ NCA algorithm of \cite{zhao2019fast}, which finds the \emph{global pivots} could be omitted without compromising on accuracy. Since \cite{gujjula2022new} follows a single tree traversal procedure than the two tree traversal procedure of \cite{zhao2019fast}, this makes the NCA reported in \cite{gujjula2022new} slightly computationally faster than the $\mathcal{O} \bkt{N}$ NCA algorithm of \cite{zhao2019fast}. The article \cite{gujjula2022new} compares their new NCA algorithm with the existing NCA algorithms (\cite{bebendorf2012constructing,zhao2019fast}) and validates their claim numerically for different kernels. We refer the readers to \cite{gujjula2022new} for more details.

\textbf{We denote the NCA with B$2$T pivot selection as B$2$T NCA and the NCA with T$2$B pivot selection as T$2$B NCA}. The standard $\mathcal{H}^2$/$\mathcal{H}^2_{\sqrt{d}}$ matrix can be formed using both B$2$T NCA and T$2$B NCA. 

\subsubsection{\texorpdfstring{$\mathcal{H}^2_{\sqrt{d}}$}{H2}(b): B\texorpdfstring{$2$}{2}T NCA to construct \texorpdfstring{$\mathcal{H}^2_{\sqrt{d}}$}{H2} matrix} \label{b2t_h2} 
In the B$2$T NCA, the pivots of a cluster at a parent level are obtained from the pivots at its child level. Therefore, one needs to traverse the $2^d$ uniform tree from the B$2$T direction. Once the pivots of all the clusters are available, we construct all the required operators (P$2$M/M$2$M, M$2$L, L$2$L/L$2$P) as described in \Cref{operator_construction}. The cost of initializing the $\mathcal{H}^2_{\sqrt{d}}$ matrix representation using the B$2$T NCA is $\mathcal{O} \bkt{N}$ \cite{zhao2019fast,gujjula2022new}. 

The $\mathcal{H}^2_{\sqrt{d}}$ MVP is performed by following the upward, transverse and downward tree traversal. Since the B$2$T NCA is applied to get the $\mathcal{H}^2_{\sqrt{d}}$ representation, we denote this algorithm as $\mathcal{H}^2_{\sqrt{d}}$(b). We use this algorithm many times for benchmarking with the proposed algorithms. Thus, for the convenience of the reader, we present it in the \Cref{b2t_h2matrix} by using similar notations to those used in this article. We also refer the reader to Algorithm $2$ of \cite{zhao2019fast}, and \cite{gujjula2022new} for more details.

\subsubsection{\texorpdfstring{$\mathcal{H}^2_{\sqrt{d}}$}{H2}(t): T\texorpdfstring{$2$}{2}B NCA to construct \texorpdfstring{$\mathcal{H}^2_{\sqrt{d}}$}{H2} matrix} \label{t2b_h2} 
In the T$2$B NCA, the pivots of a cluster at a child level are obtained from its own index set and the pivots at its parent level. Therefore, one needs to traverse the $2^d$ tree from the T$2$B direction. Once the pivots of all the clusters are available, we construct all the required operators as described in \Cref{operator_construction}. The cost of initializing the $\mathcal{H}^2_{\sqrt{d}}$ matrix representation using the T$2$B NCA is $\mathcal{O} \bkt{N \log(N)}$ \cite{zhao2019fast}. 

The $\mathcal{H}^2_{\sqrt{d}}$ MVP is performed by following the upward, transverse and downward tree traversal. Since we apply T$2$B NCA to get the $\mathcal{H}^2$ matrix, we denote this algorithm as $\mathcal{H}^2_{\sqrt{d}}$(t). We also use this algorithm for benchmarking with the proposed algorithms. Thus, for the convenience of the reader, we present it in the \Cref{t2b_h2matrix} by using similar notations to those used in this article. We refer the reader to Algorithm $1$ of \cite{zhao2019fast} for more details.

\begin{remark}
    It is worth noting that the search spaces of pivots in $\mathcal{H}^2_{\sqrt{d}}$(b) (\Cref{fig:B2T_pivot}) are smaller than that of $\mathcal{H}^2_{\sqrt{d}}$(t) (\Cref{fig:T2B_pivot}). The initialization cost of $\mathcal{H}^2_{\sqrt{d}}$(b) is $\mathcal{O} \bkt{N}$, while the initialization cost of $\mathcal{H}^2_{\sqrt{d}}$(t) is $\mathcal{O} \bkt{N \log(N)}$ (refer to \cite{zhao2019fast}). Therefore, in terms of initialization time, $\mathcal{H}^2_{\sqrt{d}}$(b) is a faster/better algorithm.
\end{remark}

\subsubsection{Construction of operators} \label{operator_construction}
Once all the pivots ( $t^{X,i}$, $s^{X,i}$, $t^{X,o}$ and $s^{X,o}$) corresponding to a cluster $X$ are available, the operators can be constructed as follows:

\begin{itemize}
    \item \textbf{L$2$P and P$2$M operators.} For a \emph{leaf} cluster $X$, the matrices $U_X$ and $V_X^*$ are given by
\end{itemize}
    \begin{equation} \label{L2P_operator}
        U_X = K_{t^X,s^{X,i}} \bkt{K_{t^{X, i}, s^{X,i}}}^{-1} 
    \end{equation}
The matrix $U_X$ is called \textbf{L$2$P} (local to particles) operator.

    \begin{equation} \label{P2M_operator}
        V_X^* = \bkt{K_{t^{X, o}, s^{X, o}}}^{-1} K_{t^{X, o}, s^X}
    \end{equation}
The matrix $V_X^*$ is called \textbf{P$2$M} (particles to multipole) operator.

\begin{itemize}
    \item \textbf{L$2$L and M$2$M operators.} For a \emph{non-leaf} cluster $ X = \displaystyle \bigcup_{c=1}^{2^d} X_c$, where $X_c$ is a child of $X$.
\end{itemize}
\begin{equation} \label{L2L_operator}
                U_{X} = \begin{bmatrix}
                    U_{X_{1}} & 0 & \hdots & 0\\
                    0 & U_{X_{2}} &  & 0\\
                    \vdots &  & \ddots & \vdots\\
                    0 & 0 & \hdots & U_{X_{2^{d}}}
                    \end{bmatrix}
                    \begin{bmatrix}
                    \Tilde{U}_{X_{1} X}\\
                    \Tilde{U}_{X_{2} X}\\
                    \vdots\\
                    \Tilde{U}_{X_{2^{d}} X}\\
                    \end{bmatrix}
\end{equation}
where
\begin{equation} 
    \Tilde{U}_{X_{c} X} = K_{t^{{X_c},i},s^{X,i}} \bkt{K_{t^{X, i}, s^{X,i}}}^{-1}, \qquad  1 \leq c \leq 2^d
\end{equation}
The matrices $\Tilde{U}_{X_{c} X}$ are called the \textbf{L$2$L} (local to local) operators.

\begin{equation} \label{M2M_operator}
                V_{X} = \begin{bmatrix}
                    V_{X_{1}} & 0 & \hdots & 0\\
                    0 & V_{X_{2}} &  & 0\\
                    \vdots &  & \ddots & \vdots\\
                    0 & 0 & \hdots & V_{X_{2^{d}}}
                    \end{bmatrix}
                    \begin{bmatrix}
                    \Tilde{V}^{*}_{X X_{1}}\\
                    \Tilde{V}^{*}_{X X_{2}}\\
                    \vdots\\
                    \Tilde{V}^{*}_{X X_{2^{d}}}\\
                    \end{bmatrix}
\end{equation}
where 
\begin{equation}
  \Tilde{V}_{X X_{c}} = \bkt{K_{t^{X, o}, s^{X, o}}}^{-1} K_{t^{X, o}, s^{{X_c},o}}, \qquad  1 \leq c \leq 2^d
\end{equation}
The matrices $\Tilde{V}^{*}_{X X_{c}}$ are called the \textbf{M$2$M} (multipole to multipole) operators. 

\begin{itemize}
    \item  \textbf{M$2$L operators.} For any cluster $X$, the \textbf{M$2$L} (multipole to local) operator is given by
\end{itemize}
    \begin{equation}
        T_{X,Y} = K_{t^{X,i}, s^{Y,o}}, \quad Y \in \mathcal{IL} \bkt{X}
    \end{equation}

\begin{remark}
    If the kernel matrix is symmetric, the P$2$M or M$2$M operators can be obtained simply by taking the transpose of L$2$P or L$2$L operators, respectively. 
\end{remark}

\section{Various pivots selection strategies on \emph{weak admissibility} condition} \label{b2t_vs_t2b}
\rk{The} NCAs \cite{bebendorf2012constructing, zhao2019fast, gujjula2022new} are primarily used to construct the nested hierarchical representations based on the \textbf{strong admissibility condition} \cref{H_admiss}, i.e., the admissible clusters are the far-field/well-separated clusters. To be more precise, all the existing literature on NCA discusses mainly strong admissibility-based $\mathcal{H}^2$ matrices. \rk{In this article, we want to construct nested bases algorithms based on our \emph{weak admissibility} condition in higher dimensions (\Cref{weak_admis}), where the admissible clusters are the far-field and the \emph{vertex-sharing} clusters. Therefore, exploring the effectiveness of B$2$T NCA and T$2$B NCA on this \emph{weak admissibility} condition is worthwhile}.

\subsection{B\texorpdfstring{$2$}{2}T NCA on our \emph{weak admissibility} condition} \label{nca_b2t_entire_hodlrdd}
\rk{The} B$2$T pivot selection selects the leaf clusters' global row and column indices and applies ACA to get the pivots. To obtain the pivots of a non-leaf cluster, ACA is applied to the pivots from the children-level clusters. This process recursively goes from bottom to top of the $2^d$ tree. Further explanation regarding the B$2$T NCA on our \rk{\emph{weak admissibility} condition in higher dimensions (\Cref{weak_admis}) is given below}.
\begin{itemize}
        \item If $X$ is a leaf cluster (childless), then construct the following four sets
    \end{itemize}
    \begin{align}
        \Tilde{t}^{X,i} := t^X \qquad \qquad \text{ and } \qquad \qquad \Tilde{s}^{X,i} := \bigcup_{Y \in \mathcal{IL}_{*}(X)} s^Y
    \end{align}
     \begin{align}
        \Tilde{t}^{X,o} := \bigcup_{Y \in \mathcal{IL}_{*}(X)} t^Y \qquad \qquad \text{ and } \qquad \qquad \Tilde{s}^{X,o} := s^X 
    \end{align}
    \begin{itemize}
        \item If $X$ is a non-leaf cluster, then construct the following four sets
    \end{itemize}
    \begin{align}
        \Tilde{t}^{X,i} := \bigcup_{X_c \in child(X)} t^{X_c, i} \qquad \qquad \text{ and } \qquad \qquad \Tilde{s}^{X,i} := \bigcup_{Y \in \mathcal{IL}_{*}(X)} \bigcup_{Y_c \in child(Y)} s^{Y_c,o}
    \end{align}
     \begin{align}
        \Tilde{t}^{X,o} := \bigcup_{Y \in \mathcal{IL}_{*}(X)} \bigcup_{Y_c \in child(Y)} t^{Y_c,i} \qquad \qquad \text{ and } \qquad \qquad \Tilde{s}^{X,o} := \bigcup_{X_c \in child(X)} s^{X_c,o} 
    \end{align}

To obtain the pivots $t^{X,i}$, $s^{X,i}$, $t^{X,o}$ and $s^{X,o}$, one needs to perform ACA \cite{aca}. We perform ACA on the matrix $K_{\Tilde{t}^{X,i}, \Tilde{s}^{X,i}}$ with user-given tolerance $\epsilon$. The sets $t^{X,i}$ and $s^{X,i}$ are the row and column pivots chosen by the ACA. The search spaces of the pivots for a particular cluster are illustrated in \Cref{fig:B2T_pivot_entire_interaction}. Similarly, perform ACA on the matrix $K_{\Tilde{t}^{X,o}, \Tilde{s}^{X,o}}$ to obtain the other two sets of pivots $t^{X,o}$ and $s^{X,o}$. Once we have the pivots, we construct the operators as in \Cref{operator_construction}. We denote this nested hierarchical representation as $\Tilde{K}_b$.

\begin{figure}[H]
\begin{center}
\captionsetup[subfloat]{labelformat=empty}
    \subfloat[]{
    \begin{tikzpicture}[scale=0.5]
        \draw[draw=black,fill=red] (-8,1) rectangle (-7,2);
        \pgfmathsetseed{58678246}
        \foreach \x in {1,...,200}
        {
		\pgfmathsetmacro{\x}{rand/2 -7.5}
		\pgfmathsetmacro{\y}{rand/2 +1.5}
		\fill[black]    (\x,\y) circle (0.015);
        };
        \pgfmathsetseed{58678246}
        \foreach \x in {1,...,80}
        {
		\pgfmathsetmacro{\x}{rand/2 -7.5}
		\pgfmathsetmacro{\y}{rand/2 +1.5}
		\fill[blue]    (\x,\y) circle (0.035);
        };
        \draw[black] (0, 0) grid (4, 4);
        \fill [red] (2,1) rectangle (3,2);
        \fill [green] (1,1) rectangle (2,2);
        \fill [white] (1.5,1.5) rectangle (2,2);
        \fill [white] (2.5,1.5) rectangle (3,2);
        \fill [white] (2.5,1) rectangle (3,1.5);
        \fill [white] (2,1) rectangle (2.5,1.5);
        \fill [white] (2,2) rectangle (2.5,2.5);
        \fill [green] (2,2) rectangle (3,3);
        \fill [green] (2,0) rectangle (3,1);
        \fill [green] (3,1) rectangle (4,2);
        \node[anchor=north] at (2.27,2.15) {\tiny $X$};
        \node[anchor=north] at (1.5,2.8) {};
        \node[anchor=north] at (1.5,1.8) {};
        \node[anchor=north] at (.5,1.8) {};
        \draw[black] (0, 0) grid (4, 4);
        \draw[step=0.5cm, black] (3, 1) grid (4, 2);
        \draw[step=0.5cm, black] (2, 0) grid (3, 1);
        \draw[step=0.5cm, black] (2, 2) grid (3, 3);
        \draw[step=0.5cm, black] (2, 1) grid (3, 2);
        \draw[step=0.5cm, black] (1, 1) grid (2, 2);
        \fill [white] (2,2) rectangle (2.5,2.5);
        \fill [white] (2.5,2) rectangle (3,2.5);
        \draw[line width=0.1mm,  black] (2.5, 2) rectangle (3, 2.5);
        \draw[line width=0.1mm,  black] (2, 2) rectangle (2.5, 2.5);
        \fill [white] (1.5,1) rectangle (2,1.5);
        \fill [white] (0,2) rectangle (2,4);
        \fill [cyan] (2.5,2) rectangle (3,2.5);
        \fill [cyan] (1.5,1) rectangle (2,1.5);
        \fill [cyan] (2.5,1) rectangle (3,1.5);
        \draw[line width=0.5mm,  black] (1, 1) rectangle (2, 2);
        \draw[line width=0.5mm,  black] (2, 2) rectangle (3, 3);
        \draw[line width=0.5mm,  black] (3, 1) rectangle (4, 2);
        \draw[line width=0.5mm,  black] (2, 0) rectangle (3, 1);
        \draw[black] (0, 2) rectangle (2, 4);
        \draw[step=1cm, black] (0, 0) grid (2, 2);

        \pgfmathsetseed{58678246}
        \foreach \x in {1,...,200}
        {
		\pgfmathsetmacro{\x}{rand/2+2.5}
		\pgfmathsetmacro{\y}{rand/2+0.5}
		\fill[black]    (\x,\y) circle (0.015);
        };

        \pgfmathsetseed{52213446}
        \foreach \x in {1,...,100}
        {
		\pgfmathsetmacro{\x}{rand/2+2.5}
		\pgfmathsetmacro{\y}{rand/4+2.75}
		\fill[black]    (\x,\y) circle (0.015);
        };

        \pgfmathsetseed{52213446}
        \foreach \x in {1,...,200}
        {
		\pgfmathsetmacro{\x}{rand/2+3.5}
		\pgfmathsetmacro{\y}{rand/2+1.5}
		\fill[black]    (\x,\y) circle (0.015);
        };

        \pgfmathsetseed{52213448}
        \foreach \x in {1,...,100}
        {
		\pgfmathsetmacro{\x}{rand/4+1.25}
		\pgfmathsetmacro{\y}{rand/2+1.5}
		\fill[black]    (\x,\y) circle (0.015);
        };

        \pgfmathsetseed{58678246}
        \foreach \x in {1,...,80}
        {
		\pgfmathsetmacro{\x}{rand/4+1.75}
		\pgfmathsetmacro{\y}{rand/4+1.25}
		\fill[black]    (\x,\y) circle (0.015);
        };

        \pgfmathsetseed{58678256}
        \foreach \x in {1,...,80}
        {
		\pgfmathsetmacro{\x}{rand/4+2.75}
		\pgfmathsetmacro{\y}{rand/4+1.25}
		\fill[black]    (\x,\y) circle (0.015);
        };

        \pgfmathsetseed{58678260}
        \foreach \x in {1,...,80}
        {
		\pgfmathsetmacro{\x}{rand/4+2.75}
		\pgfmathsetmacro{\y}{rand/4+2.25}
		\fill[black]    (\x,\y) circle (0.015);
        };

        \draw[->, line width=0.35mm] (2,1.7) -- (-7,1.7);
        \node at (-3,2.2) {Zoomed in view};

        \pgfmathsetseed{58678246}
        \foreach \x in {1,...,40}
        {
		\pgfmathsetmacro{\x}{rand/2+2.5}
		\pgfmathsetmacro{\y}{rand/2+0.5}
		\fill[blue]    (\x,\y) circle (0.035);
        };

        \pgfmathsetseed{58678246}
        \foreach \x in {1,...,20}
        {
		\pgfmathsetmacro{\x}{rand/2+2.5}
		\pgfmathsetmacro{\y}{rand/4+2.75}
		\fill[blue]    (\x,\y) circle (0.035);
        };

        \pgfmathsetseed{58678246}
        \foreach \x in {1,...,40}
        {
		\pgfmathsetmacro{\x}{rand/2+3.5}
		\pgfmathsetmacro{\y}{rand/2+1.5}
		\fill[blue]    (\x,\y) circle (0.035);
        };

        \pgfmathsetseed{58678246}
        \foreach \x in {1,...,20}
        {
		\pgfmathsetmacro{\x}{rand/4+1.25}
		\pgfmathsetmacro{\y}{rand/2+1.5}
		\fill[blue]    (\x,\y) circle (0.035);
        };

        \pgfmathsetseed{58678249}
        \foreach \x in {1,...,10}
        {
		\pgfmathsetmacro{\x}{rand/4+1.75}
		\pgfmathsetmacro{\y}{rand/4+1.25}
		\fill[blue]    (\x,\y) circle (0.035);
        };

        \pgfmathsetseed{58678248}
        \foreach \x in {1,...,10}
        {
		\pgfmathsetmacro{\x}{rand/4+2.75}
		\pgfmathsetmacro{\y}{rand/4+1.25}
		\fill[blue]    (\x,\y) circle (0.035);
        };

        \pgfmathsetseed{58678256}
        \foreach \x in {1,...,10}
        {
		\pgfmathsetmacro{\x}{rand/4+2.75}
		\pgfmathsetmacro{\y}{rand/4+2.25}
		\fill[blue]    (\x,\y) circle (0.035);
        };

    \end{tikzpicture}
    }\qquad \qquad
    \subfloat{
        \begin{tikzpicture}
            [
            box/.style={rectangle,draw=black, minimum size=0.25cm},scale=0.2
            ]
            \node[box,fill=red,,font=\tiny,label=right:Cluster considered,  anchor=west] at (-4,8){};
            \node[circle,fill=black,inner sep=0pt,minimum size=2pt,label=right:{Pivot}] (a) at (-4,6) {};
        \end{tikzpicture}
    }   
    \caption{Illustration of search spaces of the pivots in B$2$T NCA to construct a nested bases algorithm based on our weak admissibility condition in higher dimensions, which leads to poor approximation. The green and cyan-colored regions represent the interaction list of $X$, i.e., $\mathcal{IL}_{*} \bkt{X}$.}
    \label{fig:B2T_pivot_entire_interaction}
\end{center}
\end{figure}

We have shown in \cite{khan2022numerical} that for user-given tolerance $\epsilon$, the rank of interaction between two $d$ dimensional vertex-sharing hyper-cubes containing $N$ uniformly distributed particles scales $\mathcal{O} \bkt{\log(N) \log^d \bkt{\log(N)/\epsilon}}$. \rk{In our \emph{weak admissibility} condition in higher dimensions (\Cref{weak_admis}), the admissible clusters are the far-field and the vertex-sharing clusters, so the maximum rank scales as $\mathcal{O} \bkt{\log(N) \log^d \bkt{\log(N)}}$. If we apply B$2$T NCA on the \emph{weak admissibility} condition in higher dimensions}, the size of the sets $\Tilde{t}^{X,i}$, $\Tilde{s}^{X,i}$, $\Tilde{t}^{X,o}$ and $\Tilde{s}^{X,o}$ decrease as we traverse the tree from bottom to top. However, as we ascend from bottom to top of the tree, the actual size of the simulation cluster becomes larger. Consequently, the number of particles within a cluster $X$ increases and the actual vertex-sharing interaction rank grows as $\mathcal{O} \bkt{\log \bkt{\#(X)} \log^d \bkt{\log \bkt{\#(X)}}}$. Therefore, the pivots obtained from the reduced sets ($\Tilde{t}^{X,i}$, $\Tilde{s}^{X,i}$, $\Tilde{t}^{X,o}$ and $\Tilde{s}^{X,o}$) are not enough to get a good approximation of the vertex-sharing interaction matrix. We will see the effect on the rank of a vertex-sharing admissible cluster if we apply the B$2$T pivot selection. For our convenience, we assume that the $N$ particles are uniformly distributed inside the computational domain. Consider two leaf-level vertex-sharing admissible clusters containing $n_{max}$ particles. Let $R^{(l)}$ be the rank estimate of a vertex-sharing block interaction matrix at a level $l$ of the tree. At the leaf level, $\kappa$, ACA is applied to the global row and column index sets, and $R^{\bkt{\kappa}}$ scales as $\mathcal{O} \bkt{\log(n_{max}) \log^d \bkt{\log(n_{max})}}$. At a non-leaf level, two vertex-sharing admissible clusters (\Cref{fig:d2t_nl}) can be interpreted as $(2^d - 1)$ far-field clusters at the child level (\Cref{fig:d2t_far}), with a single vertex-sharing cluster at the child level (\Cref{fig:d2t_ver}) as in \Cref{fig:d2t_vs_t2d}. Considering $F^{(l)}_r$ as the rank of interaction between far-field (or well-separated) clusters and $V^{(l)}_r$ as the rank of interaction between vertex-sharing clusters at level $l$, the following recurrence relation governs a rough estimate of $R^{(l)}$ at the non-leaf level $l$.
\begin{equation} \label{b2t_useless_1}
    R^{(l)} \approx (2^d-1) F^{(l+1)}_r + V^{(l+1)}_r, \qquad l = \kappa-1:-1:1
\end{equation}
where $F^{(l)}_r \in \mathcal{O}(1)$ and $V^{(l)}_r \in \mathcal{O} \bkt{\log \bkt{R^{(l)}} \log^d \bkt{\log \bkt{R^{(l)}}}}$. 
From \cref{b2t_useless_1}, it is clear that as we traverse the tree from bottom to top, $V^{(l)}_r$ scales as a nested $``\log"$ expression, which yields a diminutive value. But the actual vertex-sharing rank grows as $\mathcal{O} \bkt{\log \bkt{N/2^{dl}} \log^d \bkt{\log \bkt{N/2^{dl}}}}$; the term $N/2^{dl}$ increase as $l$ decrease (as we go bottom to top the $l$ decrease). Thus, if we apply the B$2$T NCA to the \textbf{entire} interaction list, the rank's diminutive value is insufficient to get a good approximation of the vertex-sharing interaction block matrix. The interaction list of a cluster $\bkt{\mathcal{IL}_{*}(X)}$ may contain multiple vertex-sharing clusters, exacerbating the issue.
So, it is clear that picking the pivots in B$2$T traversal is insufficient to get a good nested hierarchical representation based on the \emph{weak admissibility} in higher dimensions.

\begin{figure}[H]
    \centering
    \subfloat[]{
    \begin{tikzpicture}[scale=0.3]
         \draw[draw=black,fill=red] (2,2) rectangle (4,4);
         \draw[draw=black,fill=cyan] (0,0) rectangle (2,2);
         \label{fig:d2t_nl}
    \end{tikzpicture}
    } \qquad  \qquad  \qquad 
    \subfloat[]{
    \begin{tikzpicture}[scale=0.3]
        \draw[black] (0, 0) grid (2, 2);
         \draw[black] (2, 2) grid (4, 4);
         \fill [cyan] (1,1) rectangle (2,2);
         \fill [green] (0,0) rectangle (1,1);
         \draw[line width=3mm,  black,thick] (0, 0) rectangle (1, 1);
         \fill [green] (1,0) rectangle (2,1);
         \draw[line width=3mm,  black,thick] (1, 0) rectangle (2, 1);
         \fill [green] (0,1) rectangle (1,2);
         \draw[line width=3mm,  black,thick] (0, 1) rectangle (1, 2);
         \fill [red] (2,2) rectangle (3,3);
         \fill [green] (2,3) rectangle (3,4);
         \draw[line width=3mm,  black,thick] (2, 3) rectangle (3, 4);
         \fill [green] (3,3) rectangle (4,4);
         \draw[line width=3mm,  black,thick] (3, 3) rectangle (4, 4);
         \fill [green] (3,2) rectangle (4,3);
         \draw[line width=3mm,  black,thick] (3, 2) rectangle (4, 3);
    \end{tikzpicture}
    } \qquad  \qquad  \qquad 
    \subfloat[]{
    \begin{tikzpicture}[scale=0.3]
         \fill [green] (0,0) rectangle (1,1);
         \draw[line width=3mm,  black,thick] (0, 0) rectangle (1, 1);
         \fill [green] (1,0) rectangle (2,1);
         \draw[line width=3mm,  black,thick] (1, 0) rectangle (2, 1);
         \fill [green] (0,1) rectangle (1,2);
         \draw[line width=3mm,  black,thick] (0, 1) rectangle (1, 2);
         \fill [green] (2,3) rectangle (3,4);
         \draw[line width=3mm,  black,thick] (2, 3) rectangle (3, 4);
         \fill [green] (3,3) rectangle (4,4);
         \draw[line width=3mm,  black,thick] (3, 3) rectangle (4, 4);
         \fill [green] (3,2) rectangle (4,3);
         \draw[line width=3mm,  black,thick] (3, 2) rectangle (4, 3);
         \label{fig:d2t_far}
    \end{tikzpicture}
    } \qquad  \qquad  \qquad 
    \subfloat[]{
    \begin{tikzpicture}[scale=0.3]
        \draw[black] (1, 1) grid (2, 2);
        \fill [cyan] (1,1) rectangle (2,2);
         \draw[black] (2, 2) grid (3, 3);
         \fill [red] (2,2) rectangle (3,3);
         \label{fig:d2t_ver}
    \end{tikzpicture}
    } \qquad 
    \caption{We can divide the vertex-sharing interaction like above in $2$D to get a rank estimate in B2T approach.}
    \label{fig:d2t_vs_t2d}
\end{figure}

Hence, we show that applying B$2$T NCA to the \textbf{entire} interaction list of a cluster $\bkt{\mathcal{IL}_{*} \bkt{X}}$ will not yield a \emph{well-approximated} nested hierarchical representation.

To verify it numerically, we consider $N=160000$ particles (same source and target) to be uniformly distributed inside the domain $[-1,1]^2$ and set the maximum particles at leaf-clusters $n_{max} = 400$. The kernel matrix $K \in \mathbb{R}^{N \times N}$ is formed using the kernel function $1/r$, where $r$ is the Euclidean distance between source and target. If $r=0$, we set the value of the kernel function as $1$. We construct the nested hierarchical representation, denoted as $\Tilde{K}_b$, based on the B$2$T pivot selection applied to the \textbf{entire} interaction list (search spaces of pivots are illustrated in \Cref{fig:B2T_pivot_entire_interaction}). The maximum rank, minimum rank, and average rank obtained from the NCA at each level of the quad tree are tabulated in \Cref{tab:d2t_matrix_approximation}. We also report the relative error in the matrix approximation, i.e.,$\magn{K - \Tilde{K}_b}_2 / \magn{K}_2$ with different NCA tolerance $(\epsilon)$. From the relative errors of \Cref{tab:d2t_matrix_approximation}, it is clear that applying the B$2$T pivot selection naively to the \textbf{entire} interaction list $\bkt{\mathcal{IL}_{*} \bkt{X}}$ does not result in a well-approximated nested hierarchical representation.

\begin{table}[H]
\centering
\resizebox{\textwidth}{!}{%
\begin{tabular}{|r|lll|l|lll|l|}
\hline
\multicolumn{1}{|l|}{\multirow{2}{*}{\makecell{Level of \\the tree}}} & \multicolumn{3}{l|}{Level wise rank in B$2$T ($\epsilon = 10^{-10}$)} & \multirow{2}{*}{\makecell{Relative error in Matrix\\ approximation ($2$-norm)}} & \multicolumn{3}{l|}{Level wise rank in B$2$T ($\epsilon = 10^{-12}$)} & \multirow{2}{*}{\makecell{Relative error in Matrix\\ approximation ($2$-norm)}} \\ \cline{2-4} \cline{6-8}
\multicolumn{1}{|l|}{} & \multicolumn{1}{l|}{Max. rank} & \multicolumn{1}{l|}{Min. Rank} & Avg. rank &  & \multicolumn{1}{l|}{Max. Rank} & \multicolumn{1}{l|}{Min. rank} & Avg. rank &  \\ \hline
1 & \multicolumn{1}{l|}{70} &  \multicolumn{1}{l|}{60} & 63    & \multirow{5}{*}{5.3697E-04} & \multicolumn{1}{l|}{96}  & \multicolumn{1}{l|}{80}  & 86 & \multirow{5}{*}{8.784E-05} \\ \cline{1-4} \cline{6-8}
2 & \multicolumn{1}{l|}{154} & \multicolumn{1}{l|}{73} & 113   &  & \multicolumn{1}{l|}{213} & \multicolumn{1}{l|}{108} & 158   & \\ \cline{1-4} \cline{6-8}
3 & \multicolumn{1}{l|}{155} & \multicolumn{1}{l|}{80} & 129   &  & \multicolumn{1}{l|}{222} & \multicolumn{1}{l|}{108} &  182  &  \\ \cline{1-4} \cline{6-8}
4 & \multicolumn{1}{l|}{147} & \multicolumn{1}{l|}{75} & 126   &  & \multicolumn{1}{l|}{202} & \multicolumn{1}{l|}{108} &  167  &  \\ \cline{1-4} \cline{6-8}
5 & \multicolumn{1}{l|}{139} & \multicolumn{1}{l|}{64} & 85    &  & \multicolumn{1}{l|}{181} & \multicolumn{1}{l|}{64} &  97   &  \\ \hline
\end{tabular}%
}
\caption{Level-wise rank and the relative error in the kernel matrix approximation, i.e., $\magn{K - \Tilde{K}_b}_2 / \magn{K}_2$. The relative errors indicate that the approximation is poor. We choose the kernel function $1/r$ with $N = 160000$ and set $n_{max} = 400$.}
\label{tab:d2t_matrix_approximation}
\end{table}

\subsection{T\texorpdfstring{$2$}{2}B NCA on our \emph{weak admissibility} condition} \label{nca_t2b_entire_hodlrdd}
\rk{In} contrast to B$2$T pivot selection, in the T$2$B pivot selection, the pivots of a cluster at a child level are obtained from its own index set and the pivots at its parent level. Therefore, one needs to traverse the $2^d$ tree from top to bottom direction to generate all the required sets of pivots for all the clusters. We discuss the T$2$B NCA on our \emph{weak admissibility} condition in higher dimensions (\Cref{weak_admis}) in detail as follows.
\begin{itemize}
        \item If $X$ has no parent, i.e., parent($X$) = NULL (parentless), then construct the following four sets
\end{itemize}
    \begin{align}
        \Tilde{t}^{X,i} := t^X \qquad \qquad \text{ and } \qquad \qquad \Tilde{s}^{X,i} := \bigcup_{Y \in \mathcal{IL}_{*}(X)} s^Y
    \end{align}
     \begin{align}
        \Tilde{t}^{X,o} := \bigcup_{Y \in \mathcal{IL}_{*}(X)} t^Y \qquad \qquad \text{ and } \qquad \qquad \Tilde{s}^{X,o} := s^X 
    \end{align}
    \begin{itemize}
        \item If $X$ has parent, i.e., parent($X$) $\neq$ NULL, then construct the following four sets
    \end{itemize}
    \begin{align}
        \Tilde{t}^{X,i} := t^X \qquad \qquad \text{ and } \qquad \qquad \Tilde{s}^{X,i} := \bigcup_{Y \in \mathcal{IL}_{*}(X)} s^{Y} \bigcup s^{parent(X),i}
    \end{align}
     \begin{align}
        \Tilde{t}^{X,o} := \bigcup_{Y \in \mathcal{IL}_{*}(X)} t^Y \bigcup t^{parent(X),o} \qquad \qquad \text{ and } \qquad \qquad \Tilde{s}^{X,o} := s^X 
    \end{align}

We perform ACA on the matrix $K_{\Tilde{t}^{X,i}, \Tilde{s}^{X,i}}$ with user-given tolerance $\epsilon$. The sets $t^{X,i}$ and $s^{X,i}$ are the row and column pivots chosen by the ACA. The search spaces of the pivots for a particular cluster are illustrated in \Cref{fig:T2B_pivot_entire_interaction}. Similarly, perform ACA on the matrix $K_{\Tilde{t}^{X,o}, \Tilde{s}^{X,o}}$ to obtain the other two sets of pivots $t^{X,o}$ and $s^{X,o}$. Once we have the pivots, we construct the operators as in \Cref{operator_construction}.

The T$2$B pivot selection does not suffer from the \textbf{reduced sets problem}, as it involves appending the parent level pivots with the admissible cluster's global index set. Hence, it possesses a larger search space from which ACA can choose the pivots, leading to a good approximation. We denote this nested hierarchical representation as $\Tilde{K}_t$.

\begin{figure}[H]
\begin{center}
\captionsetup[subfloat]{labelformat=empty}
    \subfloat[]{
    \begin{tikzpicture}[scale=0.5]
        \draw[draw=black,fill=red] (-8,1) rectangle (-7,2);
        \pgfmathsetseed{58678246}
        \foreach \x in {1,...,200}
        {
		\pgfmathsetmacro{\x}{rand/2 -7.5}
		\pgfmathsetmacro{\y}{rand/2 +1.5}
		\fill[black]    (\x,\y) circle (0.015);
        };
        \pgfmathsetseed{58678246}
        \foreach \x in {1,...,80}
        {
		\pgfmathsetmacro{\x}{rand/2 -7.5}
		\pgfmathsetmacro{\y}{rand/2 +1.5}
		\fill[blue]    (\x,\y) circle (0.035);
        };
        \draw[black] (0, 0) grid (4, 4);
        \fill [green!50] (0,0) rectangle (1,1);
        \fill [green!50] (0,1) rectangle (1,2);
        \fill [green!50] (2,3) rectangle (3,4);
        \fill [green!50] (3,3) rectangle (4,4);
        \fill [cyan!50] (1,0) rectangle (2,1);
        \fill [cyan] (1,2) rectangle (2,3);
        \fill [cyan!50] (3,2) rectangle (4,3);
        \fill [cyan!50] (3,0) rectangle (4,1);
        \fill [red] (2,1) rectangle (3,2);
        \fill [green] (1,1) rectangle (2,2);
        \fill [white] (1.5,1.5) rectangle (2,2);
        \fill [white] (2.5,1.5) rectangle (3,2);
        \fill [white] (2.5,1) rectangle (3,1.5);
        \fill [white] (2,1) rectangle (2.5,1.5);
        \fill [white] (2,2) rectangle (2.5,2.5);
        \fill [green] (2,2) rectangle (3,3);
        \fill [green] (2,0) rectangle (3,1);
        \fill [green] (3,1) rectangle (4,2);
        \node[anchor=north] at (2.27,2.15) {\tiny $X$};
        \node[anchor=north] at (1.5,2.8) {};
        \node[anchor=north] at (1.5,1.8) {};
        \node[anchor=north] at (.5,1.8) {};
        \draw[black] (0, 0) grid (4, 4);
        \draw[step=0.5cm, black] (3, 1) grid (4, 2);
        \draw[step=0.5cm, black] (2, 0) grid (3, 1);
        \draw[step=0.5cm, black] (2, 2) grid (3, 3);
        \draw[step=0.5cm, black] (2, 1) grid (3, 2);
        \draw[step=0.5cm, black] (1, 1) grid (2, 2);
        \fill [white] (2,2) rectangle (2.5,2.5);
        \fill [white] (2.5,2) rectangle (3,2.5);
        \draw[line width=0.1mm,  black] (2.5, 2) rectangle (3, 2.5);
        \draw[line width=0.1mm,  black] (2, 2) rectangle (2.5, 2.5);
        \fill [white] (1.5,1) rectangle (2,1.5);
        \fill [white] (0,2) rectangle (2,4);
        \fill [cyan] (2.5,2) rectangle (3,2.5);
        \fill [cyan] (1.5,1) rectangle (2,1.5);
        \fill [cyan] (2.5,1) rectangle (3,1.5);
        \draw[line width=0.5mm,  black] (1, 1) rectangle (2, 2);
        \draw[line width=0.5mm,  black] (2, 2) rectangle (3, 3);
        \draw[line width=0.5mm,  black] (3, 1) rectangle (4, 2);
        \draw[line width=0.5mm,  black] (2, 0) rectangle (3, 1);
        \draw[black] (0, 2) rectangle (2, 4);
        \draw[step=1cm, black] (0, 0) grid (2, 2);

        \pgfmathsetseed{58678246}
        \foreach \x in {1,...,200}
        {
		\pgfmathsetmacro{\x}{rand/2+2.5}
		\pgfmathsetmacro{\y}{rand/2+0.5}
		\fill[black]    (\x,\y) circle (0.015);
        };

        \pgfmathsetseed{58678226}
        \foreach \x in {1,...,200}
        {
		\pgfmathsetmacro{\x}{rand/2+1.5}
		\pgfmathsetmacro{\y}{rand/2+0.5}
		\fill[black]    (\x,\y) circle (0.015);
        };

        \pgfmathsetseed{58678229}
        \foreach \x in {1,...,300}
        {
		\pgfmathsetmacro{\x}{rand/2+0.5}
		\pgfmathsetmacro{\y}{rand+1}
		\fill[black]    (\x,\y) circle (0.015);
        };

        \pgfmathsetseed{58678226}
        \foreach \x in {1,...,200}
        {
		\pgfmathsetmacro{\x}{rand/2+3.5}
		\pgfmathsetmacro{\y}{rand/2+0.5}
		\fill[black]    (\x,\y) circle (0.015);
        };

        \pgfmathsetseed{58678231}
        \foreach \x in {1,...,200}
        {
		\pgfmathsetmacro{\x}{rand/2+2.5}
		\pgfmathsetmacro{\y}{rand/2+3.5}
		\fill[black]    (\x,\y) circle (0.015);
        };

        \pgfmathsetseed{58678235}
        \foreach \x in {1,...,300}
        {
		\pgfmathsetmacro{\x}{rand/2+3.5}
		\pgfmathsetmacro{\y}{rand+3}
		\fill[black]    (\x,\y) circle (0.015);
        };

        \pgfmathsetseed{52213446}
        \foreach \x in {1,...,100}
        {
		\pgfmathsetmacro{\x}{rand/2+2.5}
		\pgfmathsetmacro{\y}{rand/4+2.75}
		\fill[black]    (\x,\y) circle (0.015);
        };

        \pgfmathsetseed{52213446}
        \foreach \x in {1,...,200}
        {
		\pgfmathsetmacro{\x}{rand/2+3.5}
		\pgfmathsetmacro{\y}{rand/2+1.5}
		\fill[black]    (\x,\y) circle (0.015);
        };

        \pgfmathsetseed{52213448}
        \foreach \x in {1,...,100}
        {
		\pgfmathsetmacro{\x}{rand/4+1.25}
		\pgfmathsetmacro{\y}{rand/2+1.5}
		\fill[black]    (\x,\y) circle (0.015);
        };

        \pgfmathsetseed{58678246}
        \foreach \x in {1,...,80}
        {
		\pgfmathsetmacro{\x}{rand/4+1.75}
		\pgfmathsetmacro{\y}{rand/4+1.25}
		\fill[black]    (\x,\y) circle (0.015);
        };

        \pgfmathsetseed{58678256}
        \foreach \x in {1,...,80}
        {
		\pgfmathsetmacro{\x}{rand/4+2.75}
		\pgfmathsetmacro{\y}{rand/4+1.25}
		\fill[black]    (\x,\y) circle (0.015);
        };

        \pgfmathsetseed{58678260}
        \foreach \x in {1,...,80}
        {
		\pgfmathsetmacro{\x}{rand/4+2.75}
		\pgfmathsetmacro{\y}{rand/4+2.25}
		\fill[black]    (\x,\y) circle (0.015);
        };

        \draw[->, line width=0.35mm] (2,1.7) -- (-7,1.7);
        \node at (-3,2.2) {Zoomed in view};

        \pgfmathsetseed{58678246}
        \foreach \x in {1,...,40}
        {
		\pgfmathsetmacro{\x}{rand/2+2.5}
		\pgfmathsetmacro{\y}{rand/2+0.5}
		\fill[blue]    (\x,\y) circle (0.035);
        };

        \pgfmathsetseed{58678246}
        \foreach \x in {1,...,20}
        {
		\pgfmathsetmacro{\x}{rand/2+2.5}
		\pgfmathsetmacro{\y}{rand/4+2.75}
		\fill[blue]    (\x,\y) circle (0.035);
        };

        \pgfmathsetseed{58678246}
        \foreach \x in {1,...,40}
        {
		\pgfmathsetmacro{\x}{rand/2+3.5}
		\pgfmathsetmacro{\y}{rand/2+1.5}
		\fill[blue]    (\x,\y) circle (0.035);
        };

        \pgfmathsetseed{58678246}
        \foreach \x in {1,...,20}
        {
		\pgfmathsetmacro{\x}{rand/4+1.25}
		\pgfmathsetmacro{\y}{rand/2+1.5}
		\fill[blue]    (\x,\y) circle (0.035);
        };

        \pgfmathsetseed{58678249}
        \foreach \x in {1,...,10}
        {
		\pgfmathsetmacro{\x}{rand/4+1.75}
		\pgfmathsetmacro{\y}{rand/4+1.25}
		\fill[blue]    (\x,\y) circle (0.035);
        };

        \pgfmathsetseed{58678248}
        \foreach \x in {1,...,10}
        {
		\pgfmathsetmacro{\x}{rand/4+2.75}
		\pgfmathsetmacro{\y}{rand/4+1.25}
		\fill[blue]    (\x,\y) circle (0.035);
        };

        \pgfmathsetseed{58678256}
        \foreach \x in {1,...,10}
        {
		\pgfmathsetmacro{\x}{rand/4+2.75}
		\pgfmathsetmacro{\y}{rand/4+2.25}
		\fill[blue]    (\x,\y) circle (0.035);
        };

        \pgfmathsetseed{58678206}
        \foreach \x in {1,...,30}
        {
		\pgfmathsetmacro{\x}{rand/2+1.5}
		\pgfmathsetmacro{\y}{rand/2+0.5}
		\fill[black]    (\x,\y) circle (0.035);
        };

        \pgfmathsetseed{58678229}
        \foreach \x in {1,...,50}
        {
		\pgfmathsetmacro{\x}{rand/2+0.5}
		\pgfmathsetmacro{\y}{rand+1}
		\fill[black]    (\x,\y) circle (0.035);
        };

        \pgfmathsetseed{58678286}
        \foreach \x in {1,...,30}
        {
		\pgfmathsetmacro{\x}{rand/2+3.5}
		\pgfmathsetmacro{\y}{rand/2+0.5}
		\fill[black]    (\x,\y) circle (0.035);
        };

        \pgfmathsetseed{58678271}
        \foreach \x in {1,...,30}
        {
		\pgfmathsetmacro{\x}{rand/2+2.5}
		\pgfmathsetmacro{\y}{rand/2+3.5}
		\fill[black]    (\x,\y) circle (0.035);
        };

        \pgfmathsetseed{58678239}
        \foreach \x in {1,...,50}
        {
		\pgfmathsetmacro{\x}{rand/2+3.5}
		\pgfmathsetmacro{\y}{rand+3}
		\fill[black]    (\x,\y) circle (0.035);
        };

    \end{tikzpicture}
    }\qquad \qquad
    \subfloat{
        \begin{tikzpicture}
            [
            box/.style={rectangle,draw=black, minimum size=0.25cm},scale=0.2
            ]
            \node[box,fill=red,,font=\tiny,label=right:Cluster considered,  anchor=west] at (-4,8){};
            \node[circle,fill=black,inner sep=0pt,minimum size=2pt,label=right:{Pivot}] (a) at (-4,6) {};
        \end{tikzpicture}
    }   
    \caption{Illustration of search spaces of the pivots in T$2$B NCA to construct a nested bases algorithm based on our weak admissibility condition in higher dimensions, which leads to a good approximation. The lighter shade colors represent the interaction list at the parent level.}
    \label{fig:T2B_pivot_entire_interaction}
\end{center}
\end{figure}

To provide numerical evidence of the effectiveness of the T$2$B pivot selection, we develop the nested hierarchical representation $\Tilde{K}_t$ based on the T$2$B NCA and consider the last numerical example ($N=160000$ uniformly distributed particles and $1/r$ kernel) as in \Cref{nca_b2t_entire_hodlrdd}. We tabulate the maximum rank, minimum rank, and average rank obtained from the NCA at each level of the quad tree in \Cref{tab:t2d_matrix_approximation}. We also report the relative error in the matrix approximation, i.e.,$\magn{K - \Tilde{K}_t}_2 / \magn{K}_2$ with different NCA tolerance $(\epsilon)$. From the relative errors of \Cref{tab:t2d_matrix_approximation}, it is clear that applying the T$2$B pivot selection to the \textbf{entire} interaction list $\bkt{\mathcal{IL}_{*} \bkt{X}}$ leads to a well-approximated nested hierarchical representation, which is based on the \emph{weak admissibility} condition in higher dimensions.

\begin{table}[H]
\centering
\resizebox{\textwidth}{!}{%
\begin{tabular}{|r|lll|l|lll|l|}
\hline
\multicolumn{1}{|l|}{\multirow{2}{*}{\makecell{Level of \\the tree}}} & \multicolumn{3}{l|}{Level wise rank in T$2$B ($\epsilon = 10^{-10}$)} & \multirow{2}{*}{\makecell{Relative error in Matrix\\ approximation ($2$-norm)}} & \multicolumn{3}{l|}{Level wise rank in T$2$B ($\epsilon = 10^{-12}$)} & \multirow{2}{*}{\makecell{Relative error in Matrix\\ approximation ($2$-norm)}} \\ \cline{2-4} \cline{6-8}
\multicolumn{1}{|l|}{} & \multicolumn{1}{l|}{Max. rank} & \multicolumn{1}{l|}{Min. Rank} & Avg. rank &  & \multicolumn{1}{l|}{Max. Rank} & \multicolumn{1}{l|}{Min. rank} & Avg. rank &  \\ \hline
1 & \multicolumn{1}{l|}{102} &  \multicolumn{1}{l|}{92}   &  95      & \multirow{5}{*}{5.24526E-10} & \multicolumn{1}{l|}{127}  & \multicolumn{1}{l|}{114}  & 123 & \multirow{5}{*}{3.24098E-12} \\ \cline{1-4} \cline{6-8}
2 & \multicolumn{1}{l|}{234} & \multicolumn{1}{l|}{75}  &  153     &  & \multicolumn{1}{l|}{335} & \multicolumn{1}{l|}{142} & 218   & \\ \cline{1-4} \cline{6-8}
3 & \multicolumn{1}{l|}{222} & \multicolumn{1}{l|}{88}  &  168   &  & \multicolumn{1}{l|}{288} & \multicolumn{1}{l|}{106}   & 226   &  \\ \cline{1-4} \cline{6-8}
4 & \multicolumn{1}{l|}{187} & \multicolumn{1}{l|}{89}  &  143   &  & \multicolumn{1}{l|}{253} & \multicolumn{1}{l|}{119}   & 185   &  \\ \cline{1-4} \cline{6-8}
5 & \multicolumn{1}{l|}{150} & \multicolumn{1}{l|}{64}   &  88   &  & \multicolumn{1}{l|}{199} & \multicolumn{1}{l|}{64}     &  99   &  \\ \hline
\end{tabular}%
}
\caption{Level-wise rank and the relative error in the kernel matrix approximation, i.e., $\magn{K - \Tilde{K}_t}_2 / \magn{K}_2$. The relative errors indicate that the approximation is good. We choose the kernel function $1/r$ with $N = 160000$ and set $n_{max} = 400$.}
\label{tab:t2d_matrix_approximation}
\end{table}

Since $\Tilde{K}_t$ is a well-approximated nested hierarchical representation, it is possible to develop a fast MVP algorithm based on it, i.e., $\Tilde{K}_t \times \pmb{q}$. The MVP is performed by following the upward, transverse and downward tree traversal. We denote this algorithm as $\mathcal{H}^2_{*}$(t). The detailed algorithm can be found in \Cref{t2b_nHODt}.

Despite the good accuracy and quasi-linear complexity of the T$2$B NCA, the associated constant in the complexity estimate is very large, leading to a significant time consumption during the initialization of the representation. We demonstrate in \Cref{num_results} that the initialization cost of the $\mathcal{H}^2_{*}$(t) algorithm is high. However, there is room for improvement, leading to \textbf{better} and \textbf{efficient} algorithms, which is the main goal of this article. We discuss the proposed algorithms in the next section.

\begin{remark}
    From the \Cref{tab:d2t_matrix_approximation} and \Cref{tab:t2d_matrix_approximation}, it is evident that the average rank of both the B$2$T and T$2$B are almost the same at the leaf level (level $5$ of the tree). This is because, in the B$2$T pivot selection, the ACA is applied directly to the global index sets of the clusters. However, in the B$2$T pivot selection, as we traverse the tree upward, the pivots are selected from the \emph{reduced sets}, resulting in poor approximation. The substantial difference in the average rank between B$2$T and T$2$B pivot selection is apparent from the \Cref{tab:d2t_matrix_approximation} and \Cref{tab:t2d_matrix_approximation} at the coarser levels (level $1-4$ of the tree) of the tree.
\end{remark}

\section{Proposed algebraic hierarchical matrix algorithms for fast MVP} \label{proposed_algo} 
\rk{This section discusses two main algorithms in this article: efficient $\mathcal{H}^2_{*}$ and $\bkt{\mathcal{H}^2 + \mathcal{H}}_{*}$. The efficient $\mathcal{H}^2_{*}$ and the $\bkt{\mathcal{H}^2 + \mathcal{H}}_{*}$ are the \textbf{\emph{fully}} and the \textbf{\emph{semi/ partially}} nested bases algorithms, respectively}. Both of these proposed hierarchical matrix algorithms are based on our \emph{weak admissibility} condition in higher dimensions. 

We introduce some notations in \Cref{tab:HOD_notation}, which will be used to explain the proposed algorithms. 
 \begin{table}[H]
     \centering
     \resizebox{\textwidth}{!}{\begin{tabular}{|c|c|}
     \hline
        $\mathcal{C}$ &  Cluster of particles inside a hyper-cube at a level $l$ of the $2^d$ uniform tree. \\
    \hline
        $\mathcal{F} \bkt{{\mathcal{C}}}$  &  \makecell{Set of clusters that are in far-field of $\mathcal{C}$, i.e., hyper-cubes which are at least one hyper-cube away from $\mathcal{C}$ \\ or the clusters that are \textbf{well-separated} from $\mathcal{C}$.}  \\
    \hline
        ${HS}_{d'} \bkt{\mathcal{C}}$  &  \makecell{Set of clusters such that their corresponding hyper-cubes share $d'$ hyper-surface with the hyper-cube corresponding to $\mathcal{C}$ $(0 \leq d' \leq d-1)$. \\For example, ${HS}_{0} \bkt{\mathcal{C}}$ is the set of clusters such that their hyper-cubes share a vertex $(d'=0)$ with the hyper-cube corresponding to $\mathcal{C}$ \\and ${HS}_{1} \bkt{\mathcal{C}}$ is the set of clusters such that their hyper-cubes share an edge $(d'=1)$ with the hyper-cube corresponding to $\mathcal{C}$.}  \\
    \hline
        $child \bkt{\mathcal{C}}$  &  Set of clusters such that their hyper-cubes are children of the coarser hyper-cube corresponding to $\mathcal{C}$.  \\
    \hline
    $parent \bkt{\mathcal{C}}$  &  Set of cluster (at coarser level) with child as $\mathcal{C}$ (at finer level)  \\
    \hline
    $siblings \bkt{\mathcal{C}}$  &  $child\bkt{parent \bkt{\mathcal{C}}} \setminus \mathcal{C} $  \\
    \hline
    $clan \bkt{\mathcal{C}}$  &  $ \{ siblings \bkt{\mathcal{C}} \} \bigcup \{ child \bkt{P} : P \in \bigcup\limits_{d'=0}^{d-1} {HS}_{d'} \bkt{parent \bkt{\mathcal{C}}} \}$ \\
    \hline
    $\mathcal{IL}_{*} \bkt{\mathcal{C}}$  & \makecell{The interaction list of a cluster $\mathcal{C}$, based on our \emph{weak admissibility} condition in higher dimensions (refer to \Cref{weak_admis}), \\ is defined as $\mathcal{IL}_{*} \bkt{\mathcal{C}} = clan \bkt{\mathcal{C}} \bigcap \bkt{{HS}_{0} \bkt{\mathcal{C}} \bigcup \mathcal{F} \bkt{\mathcal{C}}}$. \rk{These are the admissible clusters for $\mathcal{H}^2_{*}$ and $\bkt{\mathcal{H}^2 + \mathcal{H}}_{*}$ representations.}}   \\
    \hline
    $\mathcal{N}_{*} \bkt{\mathcal{C}}$  & \makecell{The near-field (neighbor$+$self) list of a cluster $\mathcal{C}$ with our \emph{weak admissibility} condition in higher dimensions, \\ is defined as  $\mathcal{N}_{*} \bkt{\mathcal{C}} = \displaystyle  \bigcup_{d'=1}^{d-1} {HS}_{d'} \bkt{\mathcal{C}} \bigcup \mathcal{C}$. \rk{These are the inadmissible clusters for $\mathcal{H}^2_{*}$ and $\bkt{\mathcal{H}^2 + \mathcal{H}}_{*}$ representations.}}  \\
    \hline
    $n_{max}$  & Maximum number of particles at leaf-clusters.  \\
    \hline
    $\Tilde{K}$ & \makecell{$\Tilde{K}$ represent the hierarchical low-rank representation of the original kernel matrix $K$.} \\
    \hline
     \end{tabular}}
     \caption{Notations used to describe the $\mathcal{H}^2_{*}$ and $\bkt{\mathcal{H}^2 + \mathcal{H}}_{*}$  hierarchical matrix algorithm.}
     \label{tab:HOD_notation}
 \end{table}
 
\subsection{\texorpdfstring{$\mathcal{H}^2_{*}$(b$+$t)}{H2*}: The proposed efficient \texorpdfstring{$\mathcal{H}^2_{*}$}{H2*} matrix algorithm} \label{nHODLRdD}
\rk{The $\mathcal{H}^2_{*}$ is a hierarchical matrix algorithm with nested bases, which is based on our \emph{weak admissibility} condition in higher dimensions as described in \Cref{weak_admis}.

The admissibility criteria of $\mathcal{H}^2_{*}$ delineates that certain far-field (or well-separated) and vertex-sharing clusters are the admissible clusters, i.e., the interaction list of a cluster $X$ $\bkt{\mathcal{IL}_{*}(X)}$ consists of both the far-field and the vertex-sharing clusters}. In the \Cref{nca_b2t_entire_hodlrdd}, we have shown that if we naively apply the B$2$T NCA to the \textbf{entire} interaction list of a cluster $\bkt{\text{i.e., to } \mathcal{IL}_{*}(X)}$, it will not yield a \textbf{well-approximated} and \textbf{efficient} nested hierarchical representation because the vertex-sharing interaction rank grows polylogarithmically. Also, we know that B$2$T NCA works for the $\mathcal{H}^2_{\sqrt{d}}$ matrix (\Cref{b2t_h2}), where the admissible clusters are the far-field clusters. So, instead of applying the B$2$T pivot selection to the \textbf{entire} interaction list, if we can separate the far-field interaction from the vertex-sharing interaction, we can apply it specifically to the far-field interaction. Therefore, to construct an efficient $\mathcal{H}^2_{*}$ representation, we first partition the interaction list $\mathcal{IL}_{*} \bkt{X}$ into the \emph{far-field} interaction list denoted as $\mathcal{IL}_{far} \bkt{X}$ and the \emph{vertex-sharing} interaction list denoted as $\mathcal{IL}_{ver} \bkt{X}$, i.e., $\mathcal{IL}_{*} \bkt{X} = \mathcal{IL}_{far} \bkt{X} \bigcup \mathcal{IL}_{ver} \bkt{X}$ with $\mathcal{IL}_{far} \bkt{X} \bigcap \mathcal{IL}_{ver} \bkt{X} = \phi$. The visual representation of this partitioning corresponding to a cluster $X$ is illustrated in \Cref{fig:HODLRdD_interaction_split}. 

\begin{figure}[H]
\begin{center}
    \subfloat[\scriptsize Interaction list]{
    \begin{tikzpicture}[scale=0.5]
        \draw[black] (0, 0) grid (4, 4);
        \fill [green!50] (0,0) rectangle (1,1);
        \fill [green!50] (0,1) rectangle (1,2);
        \fill [green!50] (2,3) rectangle (3,4);
        \fill [green!50] (3,3) rectangle (4,4);
        \fill [cyan!50] (1,0) rectangle (2,1);
        \fill [cyan!50] (3,2) rectangle (4,3);
        \fill [cyan!50] (3,0) rectangle (4,1);
        \fill [red] (2,1) rectangle (3,2);
        \fill [green] (1,1) rectangle (2,2);
        \fill [white] (1.5,1.5) rectangle (2,2);
        \fill [white] (2.5,1.5) rectangle (3,2);
        \fill [cyan] (2.5,1) rectangle (3,1.5);
        \fill [white] (2,1) rectangle (2.5,1.5);
        \fill [white] (2,2) rectangle (2.5,2.5);
        \fill [green] (2,2) rectangle (3,3);
        \fill [green] (2,0) rectangle (3,1);
        \fill [green] (3,1) rectangle (4,2);
        \node[anchor=north] at (2.27,2.1) {\tiny $X$};
        \node[anchor=north] at (1.5,2.8) {};
        \node[anchor=north] at (1.5,1.8) {};
        \node[anchor=north] at (.5,1.8) {};
        \draw[black] (0, 0) grid (4, 4);
        \draw[step=0.5cm, black] (3, 1) grid (4, 2);
        \draw[step=0.5cm, black] (2, 0) grid (3, 1);
        \draw[step=0.5cm, black] (2, 2) grid (3, 3);
        \draw[step=0.5cm, black] (2, 1) grid (3, 2);
        \draw[step=0.5cm, black] (1, 1) grid (2, 2);
        \fill [cyan!25] (0,2) rectangle (2,4);
        \fill [white] (2,2) rectangle (2.5,2.5);
        \fill [cyan] (2.5,2) rectangle (3,2.5);
        \fill [cyan] (1.5,1) rectangle (2,1.5);
        \draw[line width=0.1mm,  black] (2.5, 2) rectangle (3, 2.5);
        \draw[line width=0.1mm,  black] (2, 2) rectangle (2.5, 2.5);
        \draw[line width=0.5mm,  black] (1, 1) rectangle (2, 2);
        \draw[line width=0.5mm,  black] (2, 2) rectangle (3, 3);
        \draw[line width=0.5mm,  black] (3, 1) rectangle (4, 2);
        \draw[line width=0.5mm,  black] (2, 0) rectangle (3, 1);
        \draw[black] (0, 2) rectangle (2, 4);
        \draw[step=1cm, black] (0, 0) grid (2, 2);
    \end{tikzpicture}
    }\qquad \qquad \qquad
    \subfloat[\scriptsize Far-field]{
    \begin{tikzpicture}[scale=0.5]
        \draw[black] (0, 0) grid (4, 4);
        \fill [green!50] (0,0) rectangle (1,1);
        \fill [green!50] (0,1) rectangle (1,2);
        \fill [green!50] (2,3) rectangle (3,4);
        \fill [green!50] (3,3) rectangle (4,4);
        \fill [green] (1,1) rectangle (2,2);
        \fill [white] (1.5,1.5) rectangle (2,2);
        \fill [white] (2.5,1.5) rectangle (3,2);
        \fill [white] (2.5,1) rectangle (3,1.5);
        \fill [white] (2,1) rectangle (2.5,1.5);
        \fill [white] (2,2) rectangle (2.5,2.5);
        \fill [green] (2,2) rectangle (3,3);
        \fill [green] (2,0) rectangle (3,1);
        \fill [green] (3,1) rectangle (4,2);
        \draw[step=0.5cm, black] (2, 0) grid (3, 1);
        \node[anchor=north] at (1.5,2.8) {};
        \node[anchor=north] at (1.5,1.8) {};
        \node[anchor=north] at (.5,1.8) {};
        \draw[black] (0, 0) grid (4, 4);
        \draw[step=0.5cm, black] (3, 1) grid (4, 2);
        \draw[step=0.5cm, black] (2, 0) grid (3, 1);
        \draw[step=0.5cm, black] (2, 2) grid (3, 3);
        \draw[step=0.5cm, black] (2, 1) grid (3, 2);
        \draw[step=0.5cm, black] (1, 1) grid (2, 2);
        \fill [white] (2,2) rectangle (2.5,2.5);
        \fill [white] (2.5,2) rectangle (3,2.5);
        \draw[line width=0.1mm,  black] (2.5, 2) rectangle (3, 2.5);
        \draw[line width=0.1mm,  black] (2, 2) rectangle (2.5, 2.5);
        \fill [white] (1.5,1) rectangle (2,1.5);
        \fill [white] (0,2) rectangle (2,4);
        \draw[line width=0.5mm,  black] (1, 1) rectangle (2, 2);
        \draw[line width=0.5mm,  black] (2, 2) rectangle (3, 3);
        \draw[line width=0.5mm,  black] (3, 1) rectangle (4, 2);
        \draw[line width=0.5mm,  black] (2, 0) rectangle (3, 1);
        \draw[black] (0, 2) rectangle (2, 4);
        \draw[step=1cm, black] (0, 0) grid (2, 2);
        \draw[step=0.5cm, black] (1, 1) grid (2, 2);
    \end{tikzpicture}
    }\qquad \qquad \qquad
    \subfloat[\scriptsize Vertex-sharing]{
    \begin{tikzpicture}[scale=0.5]
        \draw[black] (0, 0) grid (4, 4);
        \fill [cyan!50] (1,0) rectangle (2,1);
        \fill [cyan!50] (3,2) rectangle (4,3);
        \fill [cyan!50] (3,0) rectangle (4,1);
        \fill [cyan] (2.5,1) rectangle (3,1.5);
        \node[anchor=north] at (1.5,2.8) {};
        \node[anchor=north] at (1.5,1.8) {};
        \node[anchor=north] at (.5,1.8) {};
        \draw[black] (0, 0) grid (4, 4);
        \draw[step=0.5cm, black] (3, 1) grid (4, 2);
        \draw[step=0.5cm, black] (2, 0) grid (3, 1);
        \draw[step=0.5cm, black] (2, 2) grid (3, 3);
        \draw[step=0.5cm, black] (2, 1) grid (3, 2);
        \draw[step=0.5cm, black] (1, 1) grid (2, 2);
        \fill [cyan!25] (0,2) rectangle (2,4);
        \fill [cyan] (2.5,2) rectangle (3,2.5);
        \fill [cyan] (1.5,1) rectangle (2,1.5);
        \draw[line width=0.5mm,  black] (1, 1) rectangle (2, 2);
        \draw[line width=0.5mm,  black] (2, 2) rectangle (3, 3);
        \draw[line width=0.5mm,  black] (3, 1) rectangle (4, 2);
        \draw[line width=0.5mm,  black] (2, 0) rectangle (3, 1);
        \draw[black] (0, 2) rectangle (2, 4);
        \draw[step=1cm, black] (0, 0) grid (2, 2);
    \end{tikzpicture}
    }\qquad 
    \caption{Partitioning of the interaction list corresponding to cluster $X$ at level $3$ of the quad tree in $2$D. The interaction list of $X$ $\bkt{\mathcal{IL}_{*}(X)}$ consists of green and cyan colored clusters enclosed by a noticeable \textbf{black} border. The lighter shade colors represent the interaction list at the coarser levels of the tree (parent or grandparent of $X$).}
    \label{fig:HODLRdD_interaction_split}
\end{center}
\end{figure}

To initialize the efficient $\mathcal{H}^2_{*}$ representation, we apply the B$2$T NCA to the $\mathcal{IL}_{far} \bkt{X}$ and the T$2$B NCA to the $\mathcal{IL}_{ver} \bkt{X}$. Therefore, we need to independently traverse the $2^d$ tree twice in reverse directions with different interaction lists (B$2$T with far-field interaction and T$2$B with vertex-sharing interaction). Since we use both the B$2$T NCA (to far-field) and T$2$B NCA (to vertex-sharing), we abbreviate this algorithm as $\mathcal{H}^2_{*}$(b$+$t). Once we have all the required sets of pivots, we construct the P$2$M/M$2$M, M$2$L, L$2$L/L$2$P operators. After forming the $\mathcal{H}^2_{*}$(b$+$t) representation (or structure), we perform the MVP by following upward, transverse and downward tree traversal. We discuss our $\mathcal{H}^2_{*}$(b$+$t) hierarchical matrix algorithm for fast MVP in two steps:
\begin{enumerate}
    \item Initialization of $\mathcal{H}^2_{*}$(b$+$t) representation. (\Cref{nhdolr_init})
    \item Calculation of the potential (MVP). (\Cref{nhodlr_poten_cal})
\end{enumerate}

\subsubsection{Initialization of \texorpdfstring{$\mathcal{H}^2_{*}$(b$+$t)}{H2*} representation} \label{nhdolr_init} The initialization procedure is described in detail as follows:
\begin{enumerate}
    \item \emph{B$2$T NCA for far-field interaction.} \label{nhodlr_init_1}Since we partition the far-field and the vertex-sharing interactions associated with a cluster $X$, i.e., $\mathcal{IL}_{*} \bkt{X} = \mathcal{IL}_{far} \bkt{X} \bigcup \mathcal{IL}_{ver} \bkt{X}$, we can apply the B$2$T NCA only to the far-field interactions $\bkt{\mathcal{IL}_{far}(X)}$. The detailed procedure to obtain the pivots corresponding to a cluster $X$ is given below. 
    
\begin{itemize}
        \item If $X$ is a leaf cluster (childless), then construct the following four sets
    \end{itemize}
    \begin{align}
        \Tilde{t}^{X,i} := t^X \qquad \qquad \text{ and } \qquad \qquad \Tilde{s}^{X,i} := \bigcup_{Y \in \mathcal{IL}_{far}(X)} s^Y
    \end{align}
     \begin{align}
        \Tilde{t}^{X,o} := \bigcup_{Y \in \mathcal{IL}_{far}(X)} t^Y \qquad \qquad \text{ and } \qquad \qquad \Tilde{s}^{X,o} := s^X 
    \end{align}
    \begin{itemize}
        \item If $X$ is a non-leaf cluster, then construct the following four sets
    \end{itemize}
    \begin{align}
        \Tilde{t}^{X,i} := \bigcup_{X_c \in child(X)} t^{X_c, i} \qquad \qquad \text{ and } \qquad \qquad \Tilde{s}^{X,i} := \bigcup_{Y \in \mathcal{IL}_{far}(X)} \bigcup_{Y_c \in child(Y)} s^{Y_c,o}
    \end{align}
     \begin{align}
        \Tilde{t}^{X,o} := \bigcup_{Y \in \mathcal{IL}_{far}(X)} \bigcup_{Y_c \in child(Y)} t^{Y_c,i} \qquad \qquad \text{ and } \qquad \qquad \Tilde{s}^{X,o} := \bigcup_{X_c \in child(X)} s^{X_c,o} 
    \end{align} 

 \begin{figure}[H]
\begin{center}
\captionsetup[subfloat]{labelformat=empty}
    \subfloat[]{
    \begin{tikzpicture}[scale=0.5]
        \draw[draw=black,fill=red] (-8,1) rectangle (-7,2);
        \pgfmathsetseed{58678246}
        \foreach \x in {1,...,200}
        {
		\pgfmathsetmacro{\x}{rand/2 -7.5}
		\pgfmathsetmacro{\y}{rand/2 +1.5}
		\fill[black]    (\x,\y) circle (0.015);
        };
        \pgfmathsetseed{58678246}
        \foreach \x in {1,...,80}
        {
		\pgfmathsetmacro{\x}{rand/2 -7.5}
		\pgfmathsetmacro{\y}{rand/2 +1.5}
		\fill[blue]    (\x,\y) circle (0.035);
        };
        \draw[black] (0, 0) grid (4, 4);
        \fill [red] (2,1) rectangle (3,2);
        \fill [green] (1,1) rectangle (2,2);
        \fill [white] (1.5,1.5) rectangle (2,2);
        \fill [white] (2.5,1.5) rectangle (3,2);
        \fill [white] (2.5,1) rectangle (3,1.5);
        \fill [white] (2,1) rectangle (2.5,1.5);
        \fill [white] (2,2) rectangle (2.5,2.5);
        \fill [green] (2,2) rectangle (3,3);
        \fill [green] (2,0) rectangle (3,1);
        \fill [green] (3,1) rectangle (4,2);
        \node[anchor=north] at (2.27,2.15) {\tiny $X$};
        \node[anchor=north] at (1.5,2.8) {};
        \node[anchor=north] at (1.5,1.8) {};
        \node[anchor=north] at (.5,1.8) {};
        \draw[black] (0, 0) grid (4, 4);
        \draw[step=0.5cm, black] (3, 1) grid (4, 2);
        \draw[step=0.5cm, black] (2, 0) grid (3, 1);
        \draw[step=0.5cm, black] (2, 2) grid (3, 3);
        \draw[step=0.5cm, black] (2, 1) grid (3, 2);
        \draw[step=0.5cm, black] (1, 1) grid (2, 2);
        \fill [white] (2,2) rectangle (2.5,2.5);
        \fill [white] (2.5,2) rectangle (3,2.5);
        \draw[line width=0.1mm,  black] (2.5, 2) rectangle (3, 2.5);
        \draw[line width=0.1mm,  black] (2, 2) rectangle (2.5, 2.5);
        \fill [white] (1.5,1) rectangle (2,1.5);
        \fill [white] (0,2) rectangle (2,4);
        \draw[line width=0.5mm,  black] (1, 1) rectangle (2, 2);
        \draw[line width=0.5mm,  black] (2, 2) rectangle (3, 3);
        \draw[line width=0.5mm,  black] (3, 1) rectangle (4, 2);
        \draw[line width=0.5mm,  black] (2, 0) rectangle (3, 1);
        \draw[black] (0, 2) rectangle (2, 4);
        \draw[step=1cm, black] (0, 0) grid (2, 2);
        \draw[step=0.5cm, black] (1, 1) grid (2, 2);

        \pgfmathsetseed{58678246}
        \foreach \x in {1,...,200}
        {
		\pgfmathsetmacro{\x}{rand/2+2.5}
		\pgfmathsetmacro{\y}{rand/2+0.5}
		\fill[black]    (\x,\y) circle (0.015);
        };

        \pgfmathsetseed{52213446}
        \foreach \x in {1,...,100}
        {
		\pgfmathsetmacro{\x}{rand/2+2.5}
		\pgfmathsetmacro{\y}{rand/4+2.75}
		\fill[black]    (\x,\y) circle (0.015);
        };

        \pgfmathsetseed{52213446}
        \foreach \x in {1,...,200}
        {
		\pgfmathsetmacro{\x}{rand/2+3.5}
		\pgfmathsetmacro{\y}{rand/2+1.5}
		\fill[black]    (\x,\y) circle (0.015);
        };

        \pgfmathsetseed{52213446}
        \foreach \x in {1,...,100}
        {
		\pgfmathsetmacro{\x}{rand/4+1.25}
		\pgfmathsetmacro{\y}{rand/2+1.5}
		\fill[black]    (\x,\y) circle (0.015);
        };

        \draw[->, line width=0.35mm] (2,1.7) -- (-7,1.7);
        \node at (-3,2.2) {Zoomed in view};

        \pgfmathsetseed{58678246}
        \foreach \x in {1,...,40}
        {
		\pgfmathsetmacro{\x}{rand/2+2.5}
		\pgfmathsetmacro{\y}{rand/2+0.5}
		\fill[blue]    (\x,\y) circle (0.035);
        };

        \pgfmathsetseed{58678246}
        \foreach \x in {1,...,20}
        {
		\pgfmathsetmacro{\x}{rand/2+2.5}
		\pgfmathsetmacro{\y}{rand/4+2.75}
		\fill[blue]    (\x,\y) circle (0.035);
        };

        \pgfmathsetseed{58678246}
        \foreach \x in {1,...,40}
        {
		\pgfmathsetmacro{\x}{rand/2+3.5}
		\pgfmathsetmacro{\y}{rand/2+1.5}
		\fill[blue]    (\x,\y) circle (0.035);
        };

        \pgfmathsetseed{58678246}
        \foreach \x in {1,...,20}
        {
		\pgfmathsetmacro{\x}{rand/4+1.25}
		\pgfmathsetmacro{\y}{rand/2+1.5}
		\fill[blue]    (\x,\y) circle (0.035);
        };

    \end{tikzpicture}
    }\qquad \qquad
    \subfloat{
        \begin{tikzpicture}
            [
            box/.style={rectangle,draw=black, minimum size=0.25cm},scale=0.2
            ]
            \node[box,fill=red,,font=\tiny,label=right:Cluster considered,  anchor=west] at (-4,8){};
            \node[box,fill=green,,font=\tiny,label=right:Far-field interaction of the cluster,  anchor=west] at (-4,6){};
            \node[circle,fill=black,inner sep=0pt,minimum size=2pt,label=right:{Pivot}] (a) at (-4,4) {};
        \end{tikzpicture}
    }   
    \caption{Illustration of search spaces of the pivots in B$2$T NCA for $\mathcal{H}^2_{*}$(b$+$t) construction. Note that we consider only the far-field interaction here, i.e., $\mathcal{IL}_{far}(X) \subset \mathcal{IL}_{*}(X)$.}
    \label{fig:B2T_pivot_nHODLR}
\end{center}
\end{figure}

 We perform ACA on the matrix $K_{\Tilde{t}^{X,i}, \Tilde{s}^{X,i}}$ with user-given tolerance $\epsilon_{far}$. The sets $t^{X,i}$ and $s^{X,i}$ are the row and column pivots chosen by the ACA. The search spaces of the pivots for a particular cluster are illustrated in \Cref{fig:B2T_pivot_nHODLR}. Similarly, perform ACA on the matrix $K_{\Tilde{t}^{X,o}, \Tilde{s}^{X,o}}$ to obtain the other two sets of pivots $t^{X,o}$ and $s^{X,o}$. This process recursively goes from bottom to top of the $2^d$ tree and generates the pivots for all the clusters of the tree. As discussed earlier, the main components for constructing the operators are the four sets: $t^{X,i}$, $s^{X,i}$, $t^{X,o}$ and $s^{X,o}$. We construct different operators as described in \Cref{operator_construction}. Therefore, we get the following sets of operators corresponding to far-field interaction:
 \begin{enumerate}
     \item P$2$M $\bkt{V_{X}^{far}}$ / M$2$M $\bkt{\Tilde{V}_{X X_c}^{far}}$, \quad $X_c \in child \bkt{X}$.
     \item M$2$L $\bkt{T_{X,Y}^{far}}, \quad Y \in \mathcal{IL}_{far} \bkt{X}$.
     \item L$2$L $\bkt{\Tilde{U}_{X_c X}^{far}}$, \quad $X \in parent \bkt{X_c}$ / L$2$P $\bkt{U_{X}^{far}}$.
 \end{enumerate}

\item \emph{T$2$B NCA for vertex-sharing interaction.} As demonstrated in \Cref{nca_t2b_entire_hodlrdd}, the T$2$B NCA is well-suited for our hierarchical representation. However, instead of applying the T$2$B NCA to the entire interaction list of a cluster $X$ $(\mathcal{IL}_{*}(X))$, we employ it solely to the vertex-sharing interactions, i.e., only to the $\mathcal{IL}_{ver} \bkt{X}$. The detailed procedure for obtaining the pivots corresponding to a cluster $X$ is below. 

\begin{itemize}
        \item If $X$ has no parent, i.e., parent($X$) = NULL (parentless), then construct the following four sets
\end{itemize}
    \begin{align}
        \Tilde{t}^{X,i} := t^X \qquad \qquad \text{ and } \qquad \qquad \Tilde{s}^{X,i} := \bigcup_{Y \in \mathcal{IL}_{ver}(X)} s^Y
    \end{align}
     \begin{align}
        \Tilde{t}^{X,o} := \bigcup_{Y \in \mathcal{IL}_{ver}(X)} t^Y \qquad \qquad \text{ and } \qquad \qquad \Tilde{s}^{X,o} := s^X 
    \end{align}
    \begin{itemize}
        \item If $X$ has parent, i.e., parent($X$) $\neq$ NULL, then construct the following four sets
    \end{itemize}
    \begin{align}
        \Tilde{t}^{X,i} := t^X \qquad \qquad \text{ and } \qquad \qquad \Tilde{s}^{X,i} := \bigcup_{Y \in \mathcal{IL}_{ver}(X)} s^{Y} \bigcup s^{parent(X),i}
    \end{align}
     \begin{align}
        \Tilde{t}^{X,o} := \bigcup_{Y \in \mathcal{IL}_{ver}(X)} t^Y \bigcup t^{parent(X),o} \qquad \qquad \text{ and } \qquad \qquad \Tilde{s}^{X,o} := s^X 
    \end{align}

\begin{figure}[H]
\begin{center}
\captionsetup[subfloat]{labelformat=empty}
    \subfloat[]{
    \begin{tikzpicture}[scale=0.5]
        \draw[draw=black,fill=red] (-8,1) rectangle (-7,2);
        \pgfmathsetseed{58678246}
        \foreach \x in {1,...,200}
        {
		\pgfmathsetmacro{\x}{rand/2 -7.5}
		\pgfmathsetmacro{\y}{rand/2 +1.5}
		\fill[black]    (\x,\y) circle (0.015);
        };
        \pgfmathsetseed{58678246}
        \foreach \x in {1,...,80}
        {
		\pgfmathsetmacro{\x}{rand/2 -7.5}
		\pgfmathsetmacro{\y}{rand/2 +1.5}
		\fill[blue]    (\x,\y) circle (0.035);
        };
        \draw[black] (0, 0) grid (4, 4);
        \fill [cyan!50] (1,0) rectangle (2,1);
        \fill [cyan!50] (3,2) rectangle (4,3);
        \fill [cyan!50] (3,0) rectangle (4,1);
        \fill [red] (2,1.5) rectangle (2.5,2);
        \fill [cyan] (2.5,1) rectangle (3,1.5);
        \node[anchor=north] at (1.5,2.8) {};
        \node[anchor=north] at (1.5,1.8) {};
        \node[anchor=north] at (.5,1.8) {};
        \draw[black] (0, 0) grid (4, 4);
        \draw[step=0.5cm, black] (3, 1) grid (4, 2);
        \draw[step=0.5cm, black] (2, 0) grid (3, 1);
        \draw[step=0.5cm, black] (2, 2) grid (3, 3);
        \draw[step=0.5cm, black] (2, 1) grid (3, 2);
        \draw[step=0.5cm, black] (1, 1) grid (2, 2);
        \fill [white] (0,2) rectangle (2,4);
        \fill [cyan] (2.5,2) rectangle (3,2.5);
        \fill [cyan] (1.5,1) rectangle (2,1.5);
        \draw[line width=0.5mm,  black] (1, 1) rectangle (2, 2);
        \draw[line width=0.5mm,  black] (2, 2) rectangle (3, 3);
        \draw[line width=0.5mm,  black] (3, 1) rectangle (4, 2);
        \draw[line width=0.5mm,  black] (2, 0) rectangle (3, 1);
        \draw[black] (0, 2) rectangle (2, 4);
        \draw[step=1cm, black] (0, 0) grid (2, 2);
        \node[anchor=north] at (2.27,2.15) {\tiny $X$};

        \pgfmathsetseed{58678246}
        \foreach \x in {1,...,200}
        {
		\pgfmathsetmacro{\x}{rand/2+3.5}
		\pgfmathsetmacro{\y}{rand/2+0.5}
		\fill[black]    (\x,\y) circle (0.015);
        };

        \pgfmathsetseed{52213446}
        \foreach \x in {1,...,200}
        {
		\pgfmathsetmacro{\x}{rand/2+3.5}
		\pgfmathsetmacro{\y}{rand/2+2.5}
		\fill[black]    (\x,\y) circle (0.015);
        };

        \pgfmathsetseed{52213446}
        \foreach \x in {1,...,200}
        {
		\pgfmathsetmacro{\x}{rand/2+1.5}
		\pgfmathsetmacro{\y}{rand/2+0.5}
		\fill[black]    (\x,\y) circle (0.015);
        };

        \pgfmathsetseed{58678248}
        \foreach \x in {1,...,80}
        {
		\pgfmathsetmacro{\x}{rand/4+1.75}
		\pgfmathsetmacro{\y}{rand/4+1.25}
		\fill[blue]    (\x,\y) circle (0.015);
        };

        \pgfmathsetseed{58678249}
        \foreach \x in {1,...,80}
        {
		\pgfmathsetmacro{\x}{rand/4+2.75}
		\pgfmathsetmacro{\y}{rand/4+1.25}
		\fill[blue]    (\x,\y) circle (0.015);
        };

        \pgfmathsetseed{58678257}
        \foreach \x in {1,...,80}
        {
		\pgfmathsetmacro{\x}{rand/4+2.75}
		\pgfmathsetmacro{\y}{rand/4+2.25}
		\fill[blue]    (\x,\y) circle (0.015);
        };

        \draw[->, line width=0.35mm] (2,1.7) -- (-7,1.7);
        \node at (-3,2.2) {Zoomed in view};

        \pgfmathsetseed{58678246}
        \foreach \x in {1,...,40}
        {
		\pgfmathsetmacro{\x}{rand/2+3.5}
		\pgfmathsetmacro{\y}{rand/2+0.5}
		\fill[blue]    (\x,\y) circle (0.035);
        };

        \pgfmathsetseed{58678246}
        \foreach \x in {1,...,40}
        {
		\pgfmathsetmacro{\x}{rand/2+3.5}
		\pgfmathsetmacro{\y}{rand/2+2.5}
		\fill[blue]    (\x,\y) circle (0.035);
        };

        \pgfmathsetseed{58678246}
        \foreach \x in {1,...,40}
        {
		\pgfmathsetmacro{\x}{rand/2+1.5}
		\pgfmathsetmacro{\y}{rand/2+0.5}
		\fill[blue]    (\x,\y) circle (0.035);
        };

        \pgfmathsetseed{58678248}
        \foreach \x in {1,...,20}
        {
		\pgfmathsetmacro{\x}{rand/4+1.75}
		\pgfmathsetmacro{\y}{rand/4+1.25}
		\fill[blue]    (\x,\y) circle (0.035);
        };

        \pgfmathsetseed{58678249}
        \foreach \x in {1,...,20}
        {
		\pgfmathsetmacro{\x}{rand/4+2.75}
		\pgfmathsetmacro{\y}{rand/4+1.25}
		\fill[blue]    (\x,\y) circle (0.035);
        };

        \pgfmathsetseed{58678257}
        \foreach \x in {1,...,20}
        {
		\pgfmathsetmacro{\x}{rand/4+2.75}
		\pgfmathsetmacro{\y}{rand/4+2.25}
		\fill[blue]    (\x,\y) circle (0.035);
        };

    \end{tikzpicture}
    }\qquad \qquad
    \subfloat{
        \begin{tikzpicture}
            [
            box/.style={rectangle,draw=black, minimum size=0.25cm},scale=0.2
            ]
            \node[box,fill=red,,font=\tiny,label=right:Cluster considered,  anchor=west] at (-4,8){};
            \node[box,fill=cyan!50,,font=\tiny,label=right:Vertex-sharing interaction at parent level,  anchor=west] at (-4,6){};
            \node[box,fill=cyan,,font=\tiny,label=right:Vertex-sharing interaction of the cluster,  anchor=west] at (-4,4){};
            \node[circle,fill=black,inner sep=0pt,minimum size=2pt,label=right:{Pivot}] (a) at (-4,2) {};
        \end{tikzpicture}
    }   
    \caption{Illustration of search spaces of the pivots in T$2$B NCA for $\mathcal{H}^2_{*}$(b$+$t) construction.}
    \label{fig:T2B_pivot_nHODLR}
\end{center}
\end{figure}

 We perform ACA on the matrix $K_{\Tilde{t}^{X,i}, \Tilde{s}^{X,i}}$ with user-given tolerance $\epsilon_{ver}$. The sets $t^{X,i}$ and $s^{X,i}$ are the row and column pivots chosen by the ACA. The search spaces of the pivots for a particular cluster are illustrated in \Cref{fig:T2B_pivot_nHODLR}. Similarly, perform ACA on the matrix $K_{\Tilde{t}^{X,o}, \Tilde{s}^{X,o}}$ to obtain the other two sets of pivots $t^{X,o}$ and $s^{X,o}$. These four sets are the main components for constructing the operators. We construct the different operators corresponding to vertex-sharing interaction as described in \Cref{operator_construction}. Therefore, we get the following sets of operators corresponding to vertex-sharing interaction:
     \begin{enumerate}
         \item P$2$M $\bkt{V_{X}^{ver}}$ / M$2$M $\bkt{\Tilde{V}_{X X_c}^{ver}}$, \quad $X_c \in child \bkt{X}$.
         \item M$2$L $\bkt{T_{X,Y}^{ver}}, \quad Y \in \mathcal{IL}_{ver} \bkt{X}$.
         \item L$2$L $\bkt{\Tilde{U}_{X_c X}^{ver}}$, \quad $X \in parent \bkt{X_c}$ / L$2$P $\bkt{U_{X}^{ver}}$.
     \end{enumerate}
\end{enumerate}
\begin{remark}
    The B$2$T NCA and T$2$B NCA are performed completely independently with different interaction lists. There is no interdependency of pivots between these two methods.
\end{remark}
Therefore, after performing the B$2$T NCA (to the far-field interaction) and T$2$B NCA (to the vertex-sharing interaction) for all the clusters, we get two different sets of operators as given in \Cref{table:operators_nHODLR}.

\begin{table}[H]
      \centering
      \resizebox{6cm}{!}{%
        \begin{tabular}{l|l|l|}
        \cline{2-3}
        & Far-field  & Vertex-sharing \\ \hline
        \multicolumn{1}{|l|}{P$2$M Operators} & $V_{X}^{far}$ & $V_{X}^{ver}$ \\ \hline
        \multicolumn{1}{|l|}{M$2$M Operators} & $\Tilde{V}_{XX_c}^{far}$  & $\Tilde{V}_{X X_c}^{ver}$ \\ \hline
        \multicolumn{1}{|l|}{M$2$L Operators} & $T_{X,Y}^{far}$ & $T_{X,Y}^{ver}$ \\ \hline
        \multicolumn{1}{|l|}{L$2$L Operators} & $\Tilde{U}_{X_c X}^{far}$ & $\Tilde{U}_{X_c X}^{ver}$  \\ \hline
        \multicolumn{1}{|l|}{L$2$P Operators} & $U_{X}^{far}$ & $U_{X}^{ver}$ \\ \hline
        \end{tabular}
    }
    \caption{Two different sets of operators are readily available after performing B$2$T NCA and T$2$B NCA.}
    \label{table:operators_nHODLR}
\end{table}

\textbf{Pseudocodes of $\mathcal{H}^2_{*}$(b$+$t) initialization.}
The pseudocodes outlining the initialization of the $\mathcal{H}^2_{*}$(b$+$t) representation are presented here.\\ Initialization step I (\Cref{b2t}), Initialization step II (\Cref{t2b}).
\begin{algorithm}[H]
\scriptsize
\caption{B$2$T NCA for far-field interaction (Initialization step I)}
\label{b2t}
\begin{algorithmic}[1]
\Procedure{B2T-pivot-selection}{$\epsilon_{far}$}
    \For{\texttt{$l=\kappa : -1 : 2$}} 
        \Comment{Traverse the $2^d$ tree from bottom to top}
				\For{\texttt{$i=1:2^{dl}$}}
                    \State $X \leftarrow $ $i^{th}$ cluster at level $l$ of tree.
                    \State $I.clear()$
                    \State $J.clear()$
                    \If{$l \neq \kappa$} 
                        \State $I.insert \bkt{\bigcup t^{X_c, i}}$ \Comment{$X_c \in child(X)$, $c = 1 : 2^d$}
                        \For{$ Y \in $ \texttt{$\mathcal{IL}_{far}\bkt{{X}}$}}
						      \State $J.insert \bkt{\bigcup s^{Y_c,o}}$ \Comment{$Y_c \in child(Y)$, $c = 1 : 2^d$}
					    \EndFor
                    \Else
                    \State $I = t^X$ \Comment{$t^X$ = \text{global index set of } $X$}
                    \For{ $Y \in $  \texttt{$\mathcal{IL}_{far}\bkt{{X}}$}}
						      \State $J.insert \bkt{s^{Y}}$ \Comment{$s^Y$ = \text{global index set of } $Y$}
					\EndFor
                    \EndIf 
                \State $[t^{X,i} , s^{X,i}] = \text{ACA} \bkt{I, J, \epsilon_{far}}$ \Comment{Perform ACA on the matrix $K_{I,J}$ with tolerance $\epsilon_{far}$}
				\EndFor
	\EndFor
\EndProcedure
\Procedure{Construction-of-far-field-operators}{}
\State {One can get the other two sets of pivots $t^{X,o}$ and $s^{X,o}$ in a similar fashion. After that, all the operators corresponding to far-field are readily available (\Cref{operator_construction}).}
\EndProcedure
\end{algorithmic}
\end{algorithm}

\begin{algorithm}[H]
\scriptsize
\caption{T$2$B NCA for vertex-sharing interaction (Initialization step II)}
\label{t2b}
\begin{algorithmic}[1]
\Procedure{T2B-pivot-selection}{$\epsilon_{ver}$}
    \For{\texttt{$l=1:\kappa$}} 
        \Comment{Traverse the $2^d$ tree from top to bottom}
				\For{\texttt{$i=1:2^{dl}$}}
                    \State $X \leftarrow $ $i^{th}$ cluster at level $l$ of tree.
					\State $I =t^{X}$ \Comment{$t^X =$ \text{global index set of } $X$}
					\For{ $Y \in $  \texttt{$\mathcal{IL}_{ver}\bkt{{X}}$}}
						\State $J = s^Y$ \Comment{$s^Y =$ \text{global index set of } $Y$}
                        \If{$l > 1$} 
                            \State $J.insert \bkt{s^{parent\bkt{X},i}}$  \Comment{If column pivots of the parent exist, add it to $J$}
                        \EndIf 
					\EndFor
                \State $[t^{X,i} , s^{X,i}] = \text{ACA} \bkt{I, J, \epsilon_{ver}}$ \Comment{Perform ACA on the matrix $K_{I,J}$ with tolerance $\epsilon_{ver}$}
				\EndFor
	\EndFor
\EndProcedure
\Procedure{Construction-of-vertex-sharing-operators}{}
\State {One can get the other two sets of pivots $t^{X,o}$ and $s^{X,o}$ in a similar fashion. After that, all the operators corresponding to vertex-sharing interaction are readily available (\Cref{operator_construction}).}
\EndProcedure
\end{algorithmic}
\end{algorithm}

\subsubsection{Calculation of the potential (MVP)} \label{nhodlr_poten_cal} $\pmb{\Tilde{\phi}} = \Tilde{K} \pmb{q}$.  
All the required operators (\Cref{table:operators_nHODLR}) will be available after the initialization process. The final step is to discuss how to calculate the potential (MVP). Since we obtain two different sets of P$2$M, M$2$M, M$2$L, L$2$L, and L$2$P operators (\Cref{table:operators_nHODLR}), we use two different sets of vectors to keep track of the potential corresponding to the far-field and vertex-sharing interaction.  We independently compute the potentials corresponding to the far-field and vertex-sharing interactions by following upward, transverse and downward tree traversal. We denote the potential corresponding to the far-field interaction as \emph{far-field potential} and the potential corresponding to the vertex-sharing interaction as \emph{vertex-sharing potential}. The final potential is given by adding the far-field and vertex-sharing potentials to the near-field potential. Let the column vector $q_{X_i}^{\bkt{\kappa}}$ denotes the charge corresponding to the $i^{th}$ leaf cluster, i.e., $\pmb{q} = \bigg [q_{X_1}^{\bkt{\kappa}}; q_{X_2}^{\bkt{\kappa}}; \cdots; q_{X_{2^{d\kappa}}}^{\bkt{\kappa}} \bigg ]$ (MATLAB notation) and $X_i^{\bkt{l}}$ denotes the $i^{th}$ cluster at the level $l$ of the tree. Sometimes we omit the subscripts and superscripts from the cluster notation to improve the readability in the hope that the cluster id and the level of the tree can be understood easily from the context, i.e., instead of $X_i^{\bkt{l}}$ we use $X_i$ or $X^{(l)}$ or $X$. The procedure to compute the far-field and vertex-sharing potentials is as follows:
\begin{enumerate}
    \item \textbf{Upward traversal:} 
        \begin{itemize}
            \item Particles to multipole \emph{(P$2$M) at leaf level $\kappa$} : For all leaf clusters $X$, calculate \\ $v^{\bkt{\kappa}}_{X,far} = V_{X}^{far^*} q_X^{\bkt{\kappa}}$, \quad (Far-field) \\
            $v^{\bkt{\kappa}}_{X,ver} = V_{X}^{ver^*} q_X^{\bkt{\kappa}}$, \quad (Vertex-sharing)
            \item Multipole to multipole \emph{(M$2$M) at non-leaf level} : For all non-leaf $X$ clusters, calculate \\ $v^{\bkt{l}}_{X,far} = \dsum_{X_c \in child \bkt{X}}  \Tilde{V}_{X X_c}^{far^*} v_{X_c,far}^{\bkt{l+1}} , \quad \kappa-1 \geq l \geq 2$. (Far-field) \\
            $v^{\bkt{l}}_{X,ver} = \dsum_{X_c \in child \bkt{X}}  \Tilde{V}_{X X_c}^{ver^*} v_{X_c,ver}^{\bkt{l+1}} , \quad \kappa-1 \geq l \geq 1$. (Vertex-sharing)
        \end{itemize}
    \item  \textbf{Transverse traversal:}
        \begin{itemize}
            \item Multipole to local \emph{(M$2$L) at all levels and for all clusters} : For all cluster $X$, calculate \\
             $u^{\bkt{l}}_{X,far} = \dsum_{Y \in \mathcal{IL}_{far} \bkt{X}}  T_{X,Y}^{far} v_{Y,far}^{\bkt{l}} , \quad 2 \leq l \leq \kappa$. (Far-field) \\
            $u^{\bkt{l}}_{X,ver} = \dsum_{Y \in \mathcal{IL}_{ver} \bkt{X}}  T_{X,Y}^{ver} v_{Y,ver}^{\bkt{l}} , \quad 1 \leq l \leq \kappa$.  (Vertex-sharing)
        \end{itemize}
    \item  \textbf{Downward traversal:}
        \begin{itemize}
            \item Local to local \emph{(L$2$L) at non-leaf level} : For all non-leaf clusters $X$, calculate \\
             $u^{\bkt{l+1}}_{X_c,far} := u^{\bkt{l+1}}_{X_c,far} +    \Tilde{U}_{X_c X}^{far} u_{X,far}^{\bkt{l}} , \quad 2 \leq l \leq \kappa-1$ and $X \in parent \bkt{X_c} $. (Far-field) \\
            $u^{\bkt{l+1}}_{X_c,ver} := u^{\bkt{l+1}}_{X_c,ver} +    \Tilde{U}_{X_c X}^{ver} u_{X,ver}^{\bkt{l}} , \quad 1 \leq l \leq \kappa-1$ and $X \in parent \bkt{X_c} $. (Vertex-sharing) 
            \item Local to particles \emph{(L$2$P) at leaf level $\kappa$} : For all leaf clusters $X$, calculate \\ $\phi ^{\bkt{\kappa}}_{X,far} = U_{X}^{far} u^{\bkt{\kappa}}_{X,far}$, \quad (Far-field) \\
            $\phi ^{\bkt{\kappa}}_{X,ver} = U_{X}^{ver} u^{\bkt{\kappa}}_{X,ver}$, \quad (Vertex-sharing)
        \end{itemize}
\end{enumerate}

\textbf{Near-field potential and the total potential at leaf level.}
For each leaf cluster $X$, we add the near-field (neighbor$+$self) potential, which is a direct computation to the far-field and vertex-sharing potentials. Hence, the final computed potential of a leaf cluster $X$ is given by
\begin{align}
    \phi ^{\bkt{\kappa}}_{X} = \underbrace{\phi ^{\bkt{\kappa}}_{X,far}}_{ \text{Far-field potential (nested)}} + \underbrace{\phi ^{\bkt{\kappa}}_{X,ver}}_{\text{Vertex-sharing potential (nested)}} +   \underbrace{\dsum_{X' \in \mathcal{N}_{*} \bkt{X}}  K_{t^X,s^{X'}} q_{X'}^{\bkt{\kappa}}}_{\text{ Near-field potential (direct)}}
\end{align}
If $\kappa = 1$, then no far-field interaction exists so $\phi ^{\bkt{\kappa}}_{X,far}=0$.
We know that $X_i$ is the $i^{th}$ leaf cluster and $\phi ^{\bkt{\kappa}}_{X_i}$ represents the potential corresponding to it, $1 \leq i \leq 2^{d \kappa}$. 

\begin{figure}[H]
    \centering
    \captionsetup[subfloat]{labelformat=empty}
    \subfloat{
        \begin{tikzpicture}
        \node[anchor=west] at (0,-1.25){$\Tilde{\phi} = $};
        \filldraw  (1,0) to[out=260,in=100]  (1,-2.5) to [out=98,in=262] cycle;
        \end{tikzpicture}
    }
    \subfloat[]{
        \includegraphics[scale=0.25]{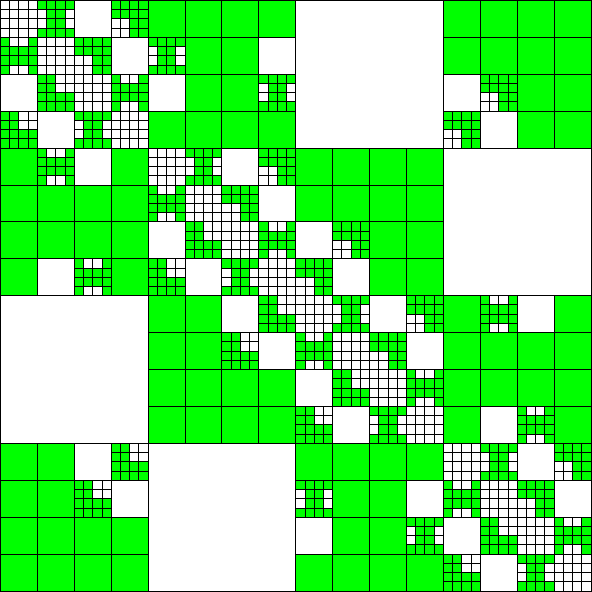}
        }\quad
    \subfloat[]{
    \begin{tikzpicture}
    \node at (-1.5,0) {};
    \node at (-1.5,1.1) {$+$};
    \end{tikzpicture}
    }\quad
    \subfloat[]{
        \includegraphics[scale=0.25]{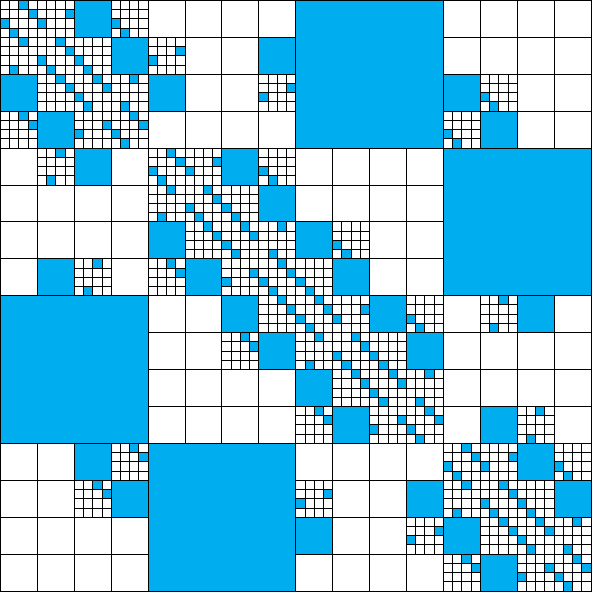}
        }\quad
    \subfloat[]{
        \begin{tikzpicture}
        \node at (-3,0) {};
        \node at (-3,1.1) {$+$};
        \end{tikzpicture}
    }\quad
    \subfloat[]{
        \includegraphics[scale=0.25]{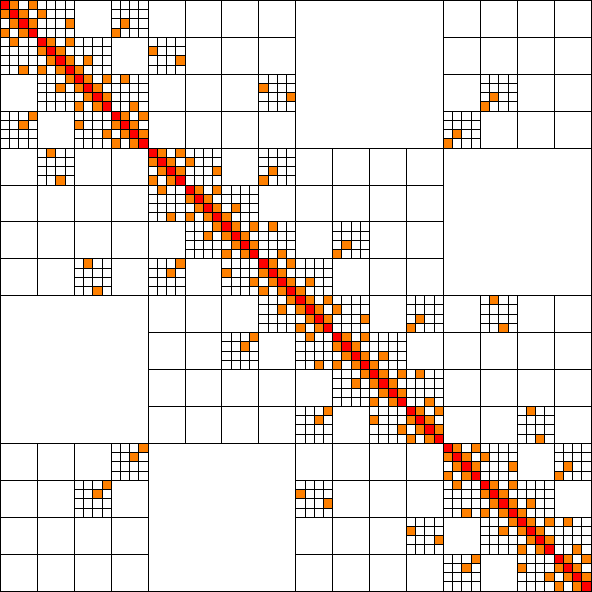}
        }
    \subfloat{
        \begin{tikzpicture}
        \filldraw (2.75,0) to[out=-80,in=80]  (2.75,-2.5) to [out=82,in=-82] cycle;
        \end{tikzpicture}
    }\quad
    \subfloat{
        \begin{tikzpicture}
        \node[anchor=west] at (-1,-2.75){$\times$};
        \draw[draw=black, fill=blue!40] (0, -4) rectangle (0.12, -1.5);
        \node[anchor=west] at (-0.15,-2.75){$\pmb{q}$};
        \end{tikzpicture}
    }\quad
    \caption{In the efficient $\mathcal{H}^2_{*}$/ $\mathcal{H}^2_{*}$(b$+$t) algorithm, the operators (P$2$M/M$2$M, M$2$L and L$2$L/L$2$P) corresponding to the far-field and the vertex-sharing interaction are constructed separately using B$2$T NCA and T$2$B NCA, respectively. After that, we calculate the far-field and vertex-sharing potentials independently by following Upward, Transverse and Downward tree traversal. We get the final potential by adding the far-field and vertex-sharing potentials to the near-field potential.}
    \label{fig:nHODLR2D_construction}
\end{figure}

Therefore, the computed potential is given by $\pmb{\Tilde{\phi}} =\bigg [\phi ^{\bkt{\kappa}}_{X_1};\phi ^{\bkt{\kappa}}_{X_2}; \cdots ;\phi ^{\bkt{\kappa}}_{X_{2^{d \kappa}}} \bigg ]$ (MATLAB notation). The schematic representation of the proposed $\mathcal{H}^2_{*}$(b$+$t) algorithm is given in \Cref{fig:nHODLR2D_construction}. 

\begin{remark}
    The difference between $\mathcal{H}^2_{*}$(b$+$t) and $\mathcal{H}^2_{*}$(t) (refer to \Cref{nca_t2b_entire_hodlrdd} and \Cref{t2b_nHODt}) lies in their construction techniques. (\Cref{fig:B2T_pivot_nHODLR}, \Cref{fig:T2B_pivot_nHODLR}) and \Cref{fig:T2B_pivot_entire_interaction} illustrate this. Since the search spaces for the pivots are small in $\mathcal{H}^2_{*}$(b$+$t), it takes less time to initialize the representation and is more efficient. This enhancement is illustrated in \Cref{num_results}.
\end{remark}

\subsubsection{Complexity analysis of \texorpdfstring{$\mathcal{H}^2_{*}$(b$+$t)}{H2*}}
We construct the $\mathcal{H}^2_{*}$(b$+$t) hierarchical representation using B$2$T NCA and T$2$B NCA. For a particular cluster $X$, its far-field interaction $\mathcal{IL}_{far} \bkt{X}$ and vertex-sharing interaction $\mathcal{IL}_{ver} \bkt{X}$ are compressed using the B$2$T NCA (\cref{b2t}) and T$2$B NCA (\cref{t2b}), respectively. 

\textbf{Time complexity.} The time complexity of the initialization steps and the potential calculation steps are given below
\begin{itemize}
    \item \emph{Far-field interaction compression: } Let the leaf cluster size be bounded by $p_1$ and we also assume $p_1 = \mathcal{O} \bkt{n_{max}}$. Let $c_{far}$ be the maximum of the far-field interaction list size of a cluster (in $2$D and $3$D for the $\mathcal{H}^2_{*}$ hierarchical representation $c_{far}=12$ and $126$, respectively, refer to \cref{hodlrdd_interaction_list}). For a non-leaf cluster, we choose the row indices from the row pivots of the children clusters and column indices from the column pivots of the children of the admissible far-field blocks. After that, the partially pivoted ACA \cite{aca} is applied to the row and column indices. The cost for applying ACA for a non-leaf cluster at level $l$ is bounded by $ \bkt{2^{d}p_1 + c_{far}2^{d}p_1}p_1^2$. By taking all the levels, the cost for the far-field compression is bounded by $\dsum_{l=2}^{\kappa} 2^{dl} \bkt{2^{d}p_1 + c_{far}2^{d}p_1}p_1^2 \approx \mathcal{O} \bkt{N}$. So, the complexity of \Cref{b2t} is $\mathcal{O} \bkt{N}$.
    \item \emph{Vertex-sharing interaction compression: } Let $c_{ver}$ be the maximum of the vertex-sharing interaction list size of a cluster (in $2$D and $3$D for the $\mathcal{H}^2_{*}$ hierarchical representation $c_{ver}=3$ and $7$, respectively, refer to \cref{hodlrdd_interaction_list}). For a cluster (with a parent), we choose the row indices from its index set and column indices from its index set, along with column pivots of the parent of the admissible vertex-sharing blocks. After that, the partially pivoted ACA \cite{aca} is applied to the row and column indices. Let $p_2$ be the maximum size of parent level pivots across the tree to be added, where $p_2$ scales at most polylogarithmically with $N$ \cite{khan2022numerical} and this is a loose upper bound of $p_2$. Hence, the cost for applying ACA for a particular cluster is bounded by $\bkt{\dfrac{N}{2^{dl}} + c_{ver} \dfrac{N}{2^{dl}} +p_2} p_2^2$. By taking all the levels, the cost for the vertex-sharing interaction compression is bounded by $\dsum_{l=1}^{\kappa} 2^{dl} \bkt{\dfrac{N}{2^{dl}} + c_{ver} \dfrac{N}{2^{dl}} +p_2} p_2^2 = \dsum_{l=1}^{\kappa} \bkt{N + c_{ver} N +p_2} p_2^2 \approx \mathcal{O} \bkt{p_2 N \log \bkt{N}}$. So, the complexity of \Cref{t2b} is $\mathcal{O} \bkt{p_2 N \log \bkt{N}}$.
\end{itemize}
    Thus, the overall time complexity to initialize the $\mathcal{H}^2_{*}$(b$+$t) representation scales as $\mathcal{O} \bkt{N + p_2 N \log \bkt{N}}$.
    
The time required for each step in the calculation of the far-field and vertex-sharing potentials (\Cref{nhodlr_poten_cal}) exhibits linear and quasi-linear scaling, respectively. Therefore, the overall time complexity for the potential calculation (MVP) is quasi-linear for non-oscillatory kernels. 

\textbf{Space complexity.}
For a cluster $X$, we store two different sets of operators as in \Cref{table:operators_nHODLR}. The total cost for storing all the P$2$M/M$2$M, M$2$L, and L$2$L/L$2$P  operators is $\mathcal{O} \bkt{Np_1^2+Np_2^2}$. The cost of storing the dense near-field operators at leaf level is $\mathcal{O} \bkt{N n_{max}^2}$. Therefore, the overall memory cost scales at most quasi-linearly.
\begin{remark}
    The quasi-linear complexity is an overestimation because the far-field part scales linearly. The $\mathcal{H}^2_{*}$(b$+$t) exhibits similar scaling when comparing its actual performance with the standard $\mathcal{H}^2$/$\mathcal{H}^2_{\sqrt{d}}$ algorithms \cite{zhao2019fast}. In fact, in the \Cref{num_results}, we show that $\mathcal{H}^2_{*}$(b$+$t) is competitive with the $\mathcal{H}^2_{\sqrt{d}}$ algorithms (\Cref{b2t_h2matrix}, \Cref{t2b_h2matrix}) with respect to memory and MVP time. The step-by-step costs of $\mathcal{H}^2_{*}$(b$+$t), along with a comparison of the same in algebraic $\mathcal{H}^2_{\sqrt{d}}$, are presented in \Cref{tab:complexity_com}.
\end{remark}
\begin{table}[H]
\centering
\resizebox{15cm}{!}{%
\begin{tabular}{|l|ll|ll|}
\hline
\multirow{2}{*}{\quad STEP} & \multicolumn{2}{l|}{ \qquad \qquad \qquad $2$D $(d=2)$} & \multicolumn{2}{l|}{ \qquad \qquad \qquad $3$D $(d=3)$} \\ \cline{2-5} 
 & \multicolumn{1}{l|}{$\mathcal{H}^2_{*}$(b$+$t) in $2$D} & $\mathcal{H}^2_{\sqrt{d}}$ in $2$D & \multicolumn{1}{l|}{$\mathcal{H}^2_{*}$(b$+$t) in $3$D} & $\mathcal{H}^2_{\sqrt{d}}$ in $3$D\\ \hline
P$2$M$+$M$2$M &  \multicolumn{1}{l|}{$\mathcal{O} \bkt{Np_1^2 + Np_2^2}$}       &  $\mathcal{O} \bkt{Np^2}$      & \multicolumn{1}{l|}{$\mathcal{O} \bkt{Np_1^2 + Np_2^2}$}    & $\mathcal{O} \bkt{Np^2}$ \\  
M$2$L & \multicolumn{1}{l|}{$\mathcal{O} \bkt{12 N p_1^2 + 3 N p_2^2}$}       &   $\mathcal{O} \bkt{27Np^2}$     & \multicolumn{1}{l|}{$\mathcal{O} \bkt{126 N p_1^2 + 7 N p_2^2}$}    & $\mathcal{O} \bkt{189 Np^2}$ \\ 
L$2$L$+$L$2$P & \multicolumn{1}{l|}{$\mathcal{O} \bkt{Np_1^2 + Np_2^2}$}       &   $\mathcal{O} \bkt{Np^2}$     & \multicolumn{1}{l|}{$\mathcal{O} \bkt{Np_1^2 + Np_2^2}$}    & $\mathcal{O} \bkt{Np^2}$ \\ 
Near-field & \multicolumn{1}{l|}{$\mathcal{O} \bkt{5N n_{max}^2}$}  &   $\mathcal{O} \bkt{9N n_{max}^2}$     & \multicolumn{1}{l|}{$\mathcal{O} \bkt{19 N n_{max}^2}$}    & $\mathcal{O} \bkt{27 N n_{max}^2}$ \\ \hline
\end{tabular}%
}
\caption{Comparison of complexities between $\mathcal{H}^2_{*}$(b$+$t) and $\mathcal{H}^2_{\sqrt{d}}$ for $d=2, 3$, where $p_1$ and $p$ scales as $\mathcal{O}(1)$ but $p_2$ scales at most polylogarithmically with $N$. The constants corresponding to Near-field and M$2$L steps are discussed in \cref{hodlrdd_nbd_list} and \cref{hodlrdd_interaction_list}.}
\label{tab:complexity_com}
\end{table}

\subsection{\texorpdfstring{$\bkt{\mathcal{H}^2 + \mathcal{H}}_{*}$}{H1.5}: The proposed semi-nested hierarchical matrix algorithm} \label{s-nHODLRdD}
\rk{The $\bkt{\mathcal{H}^2 + \mathcal{H}}_{*}$ is a hierarchical matrix algorithm with semi-nested (or partially nested) bases. It is also based on our \emph{weak admissibility} condition in higher dimensions (\Cref{weak_admis}}). We term this algorithm \emph{semi-nested} because we employ both the nested and non-nested bases in the representation. Therefore, $\bkt{\mathcal{H}^2 + \mathcal{H}}_{*}$ differs from the $\mathcal{H}^2_{*}$(b$+$t) algorithm in that it is not a \textbf{fully} nested algorithm. Let $X$ be a cluster at a particular level of the tree. We know that the interaction list of $X$ consists of far-field and vertex-sharing interactions; and it can be partitioned, i.e., $\mathcal{IL}_{*} \bkt{X} = \mathcal{IL}_{far} \bkt{X} \bigcup \mathcal{IL}_{ver} \bkt{X}$ (refer to \Cref{fig:HODLRdD_interaction_split}). To initialize the $\bkt{\mathcal{H}^2 + \mathcal{H}}_{*}$ representation, we apply the B$2$T NCA to $\mathcal{IL}_{far} \bkt{X}$ and partially pivoted ACA to $\mathcal{IL}_{ver} \bkt{X}$. It is to be noted that the far-field interaction compression routine is the same as the proposed $\mathcal{H}^2_{*}$(b$+$t) algorithm. The main difference between the $\mathcal{H}^2_{*}$(b$+$t) and the $\bkt{\mathcal{H}^2 + \mathcal{H}}_{*}$ representations is the vertex-sharing interaction compression routine. Once we have the $\bkt{\mathcal{H}^2 + \mathcal{H}}_{*}$ representation, i.e., all the required operators are available, we calculate the potential or MVP. The potential corresponding to the far-field interaction, i.e., the far-field potential, is obtained through upward, transverse and downward tree traversal. The potential corresponding to the vertex-sharing interaction, i.e., the vertex-sharing potential, is derived using $\mathcal{H}$ matrix-like algorithm (non-nested approach). We discuss our $\bkt{\mathcal{H}^2 + \mathcal{H}}_{*}$ hierarchical matrix algorithm for fast MVP in the following steps:
\begin{enumerate}
    \item Initialization of $\bkt{\mathcal{H}^2 + \mathcal{H}}_{*}$ representation. (\Cref{snhdolr_init})
    \item Calculation of the potential (MVP). (\Cref{snhodlr_poten_cal})
\end{enumerate}

\subsubsection{Initialization of \texorpdfstring{$\bkt{\mathcal{H}^2 + \mathcal{H}}_{*}$}{H1.5} representation} \label{snhdolr_init} The initialization procedure is described in detail as follows:
\begin{enumerate}
    \item \textit{B$2$T NCA for far-field interaction.} We apply the B$2$T NCA only to the $\mathcal{IL}_{far}(X)$. This step is analogous to the first step of $\mathcal{H}^2_{*}$(b$+$t) initialization (\Cref{nhodlr_init_1}). So, upon completion of the B$2$T NCA, we will get the following sets of operators:  
    \begin{enumerate}
     \item P$2$M $\bkt{V_{X}^{far}}$ / M$2$M $\bkt{\Tilde{V}_{X X_c}^{far}}$, \quad $X_c \in child \bkt{X}$.
     \item M$2$L $\bkt{T_{X,Y}^{far}}, \quad Y \in \mathcal{IL}_{far} \bkt{X}$.
     \item L$2$L $\bkt{\Tilde{U}_{X_c X}^{far}}$, \quad $X \in parent \bkt{X_c}$ / L$2$P $\bkt{U_{X}^{far}}$.
 \end{enumerate}
    \item  \textit{ACA-based compression for vertex-sharing interaction.} We compress the vertex-sharing interaction using partially pivoted ACA \cite{aca}. For each pair of clusters $X$ and $Y \in \mathcal{IL}_{ver} \bkt{X}$ at all the levels of the $2^d$ tree, the vertex-sharing interaction matrix $\bkt{K_{t^X, s^Y}}$ is compressed by partially pivoted ACA with user-given tolerance $\epsilon_{ver}$ as follows:
    \begin{align}
        K_{t^X, s^Y} \approx U V^*,\quad \text{$t^X$ and $s^Y$ are the global index sets of the clusters $X$ and $Y$, respectively}.
    \end{align}
\end{enumerate}

Therefore, after performing the B$2$T NCA (to the far-field) and ACA (to the vertex-sharing) for all the clusters, we get two different sets of operators as given in \Cref{table:operators_snHODLR}.

\begin{table}[H]
    \centering
    \subfloat[\scriptsize Far-field operators]{        
    \resizebox{4cm}{!}{%
        \begin{tabular}{l|l|}
        \cline{2-2}
                                  & Far-field \\ \hline
        \multicolumn{1}{|l|}{P$2$M Operators}    & $V_{X}^{far}$ \\ \hline
        \multicolumn{1}{|l|}{M$2$M Operators}    & $\Tilde{V}_{X X_c}^{far}$  \\ \hline
        \multicolumn{1}{|l|}{M$2$L Operators}    & $T_{X,Y}^{far}$ \\ \hline
        \multicolumn{1}{|l|}{L$2$L Operators}    & $\Tilde{U}_{X_c X}^{far}$ \\ \hline
        \multicolumn{1}{|l|}{L$2$P Operators}    & $U_{X}^{far}$ \\ \hline
        \end{tabular}
        }}  \qquad \qquad%
    \subfloat[\scriptsize Vertex-sharing operators]{        \resizebox{3.5cm}{!}{%
        \begin{tabular}{|l|}
        \hline
        \scriptsize Vertex-sharing \\ \hline
        $U$ \\ \hline
        $V$ \\ \hline
        \end{tabular}
        }}  \qquad
        \caption{Different operators after performing the B$2$T NCA on far-field and ACA on vertex-sharing interaction.}
        \label{table:operators_snHODLR}
\end{table}

\textbf{Pseudocodes of \texorpdfstring{$\bkt{\mathcal{H}^2 + \mathcal{H}}_{*}$}{H1.5} initialization.}
The pseudocodes for initializing $\bkt{\mathcal{H}^2 + \mathcal{H}}_{*}$ representation are given here.\\
Initialization step I (\Cref{b2t}), Initialization step II (\Cref{pure_aca})
 \begin{algorithm}[H]
 \scriptsize
	\caption{Compression of the vertex-sharing interaction using only ACA (Initialization step II)}\label{pure_aca}
	\begin{algorithmic}[1]
		\Procedure{VERTEX-SHARING-interaction-compression}{$\epsilon_{ver}$}
			\For{\texttt{$l=1:\kappa$}} 
				\For{\texttt{$i=1:2^{dl}$}}
                    \State $X \leftarrow $ $i^{th}$ cluster at level $l$ of tree.
					\State $I = t^X$ \Comment{$t^X = $ \text{ global index set of } $X$}
					\For{$Y \in \mathcal{IL}_{ver} \bkt{X}$}
						\State $J = s^Y$ \Comment{$s^Y = $ \text{ global index set of } $Y$}
                          \State $[U, V] = \text{ACA}(I, J, \epsilon_{ver})$ \Comment{Perform ACA on the matrix $K_{I,J}$, $K_{I,J} \approx U V^{*}$}
					\EndFor
				\EndFor
			\EndFor
		\EndProcedure
	\end{algorithmic}
\end{algorithm}
\subsubsection{Calculation of the potential (MVP)} \label{snhodlr_poten_cal} $\pmb{\Tilde{\phi}} = \Tilde{K} \pmb{q}$. After the initialization process, all the required operators (\Cref{table:operators_snHODLR}) will be available, and the final step is to discuss the calculation of the potential (MVP). We compute the far-field potential using upward, transverse and downward tree traversal, which is similar to the far-field potential computation of the $\mathcal{H}^2_{*}$(b$+$t) algorithm. The vertex-sharing potential is calculated using a single tree traversal in a manner analogous to the $\mathcal{H}$ matrix-like algorithms. Thus, the $\bkt{\mathcal{H}^2 + \mathcal{H}}_{*}$ algorithm can be thought of as a cross of $\mathcal{H}^2$ and $\mathcal{H}$ matrix-like algorithms. The final potential is given by adding the far-field, vertex-sharing and near-field potentials. 

The procedure to compute the far-field potential is as follows:

\begin{enumerate}
    \item \textbf{Upward traversal:} 
        \begin{itemize}
            \item Particles to multipole \emph{(P$2$M) at leaf level $\kappa$} : For all leaf clusters $X$, calculate \\ $v^{\bkt{\kappa}}_{X,far} = V_{X}^{far^*} q_X^{\bkt{\kappa}}$, \quad (Far-field)
            \item Multipole to multipole \emph{(M$2$M) at non-leaf level} : For all non-leaf $X$ clusters, calculate \\ $v^{\bkt{l}}_{X,far} = \dsum_{X_c \in child \bkt{X}}  \Tilde{V}_{X X_c}^{far^*} v_{X_c,far}^{\bkt{l+1}} , \quad \kappa-1 \geq l \geq 2$. (Far-field)
        \end{itemize}
    \item  \textbf{Transverse traversal:}
        \begin{itemize}
            \item Multipole to local \emph{(M$2$L) at all levels and for all clusters} : For all cluster $X$, calculate \\
             $u^{\bkt{l}}_{X,far} = \dsum_{Y \in \mathcal{IL}_{far} \bkt{X}}  T_{X,Y}^{far} v_{Y,far}^{\bkt{l}} , \quad 2 \leq l \leq \kappa$. (Far-field) 
        \end{itemize}
    \item  \textbf{Downward traversal:}
        \begin{itemize}
            \item Local to local \emph{(L$2$L) at non-leaf level} : For all non-leaf clusters $X$, calculate \\
             $u^{\bkt{l+1}}_{X_c,far} := u^{\bkt{l+1}}_{X_c,far} +    \Tilde{U}_{X_c X}^{far} u_{X,far}^{\bkt{l}} , \quad 2 \leq l \leq \kappa-1$ and $X \in parent \bkt{X_c} $. (Far-field) 
            \item Local to particles \emph{(L$2$P) at leaf level $\kappa$} : For all leaf clusters $X$, calculate \\ $\phi ^{\bkt{\kappa}}_{X,far} = U_{X}^{far} u^{\bkt{\kappa}}_{X,far}$, \quad (Far-field)
        \end{itemize}
\end{enumerate}

The calculation of the vertex-sharing potential $\pmb{\phi}_{ver}$ is executed in a non-nested manner, as given in \Cref{s-nHODLRdD_ver}.

\textbf{Pseudocode of calculating the vertex-sharing potential in \texorpdfstring{$\bkt{\mathcal{H}^2 + \mathcal{H}}_{*}$}{H1.5} algorithm.} The pseudocode outlining the calculation of the vertex-sharing potential ($\pmb{\phi}_{ver}$) is presented here.
 \begin{algorithm}[H]
 \scriptsize
	\caption{Compute the vertex-sharing potential in $\bkt{\mathcal{H}^2 + \mathcal{H}}_{*}$}\label{s-nHODLRdD_ver}
	\begin{algorithmic}[1]
		\Procedure{VERTEX-SHARING potential}{$\pmb{\phi}_{ver}$}
			\For{\texttt{$l=1:\kappa$}} 
				\For{\texttt{$i=1:2^{dl}$}}
                    \State $X \leftarrow $ $i^{th}$ cluster at level $l$ of tree.
					\For{$Y \in \mathcal{IL}_{ver} \bkt{X}$}. 
                        \State $\pmb{\phi}_{ver} \bkt{X} := \pmb{\phi}_{ver} \bkt{X} + U \bkt{V^* q_Y^{\bkt{l}}}$ \Comment{Low-rank MVP, MATLAB notation is applied to $\pmb{\phi}_{ver}$}
					\EndFor
				\EndFor
			\EndFor
            \State \textbf{return} $\pmb{\phi}_{ver}$
		\EndProcedure
	\end{algorithmic}
\end{algorithm}

\textbf{Near-field potential and the total potential at leaf level.}
For each leaf cluster $X$, we add the near-field (neighbor$+$self) potential, which is a direct computation with the far-field and vertex-sharing potential. Hence, the final computed potential of a leaf cluster is given by
\begin{align}
    \phi ^{\bkt{\kappa}}_{X} = \underbrace{\phi ^{\bkt{\kappa}}_{X,far}}_{ \text{Far-field potential (nested)}} + \underbrace{\pmb{\phi}_{ver} \bkt{X}}_{\text{Vertex-sharing potential (non-nested)}} +   \underbrace{\dsum_{X' \in \mathcal{N}_{*} \bkt{X}}  K_{t^X,s^{X'}} q_{X'}^{\bkt{\kappa}}}_{\text{ Near-field potential (direct)}}
\end{align}
If $\kappa = 1$, then no far-field interaction exists so $\phi ^{\bkt{\kappa}}_{X,far}=0$.
We know that $X_i$ is the $i^{th}$ leaf cluster and $\phi ^{\bkt{\kappa}}_{X_i}$ represents the potential corresponding to it, $1 \leq i \leq 2^{d \kappa}$. 

\begin{figure}[H]
    \centering
    \captionsetup[subfloat]{labelformat=empty}
    \subfloat{
        \begin{tikzpicture}
        \node[anchor=west] at (0,-1.25){$\Tilde{\phi} = $};
        \filldraw  (1,0) to[out=260,in=100]  (1,-2.5) to [out=98,in=262] cycle;
        \end{tikzpicture}
    }
    \subfloat[]{
        \includegraphics[scale=0.25]{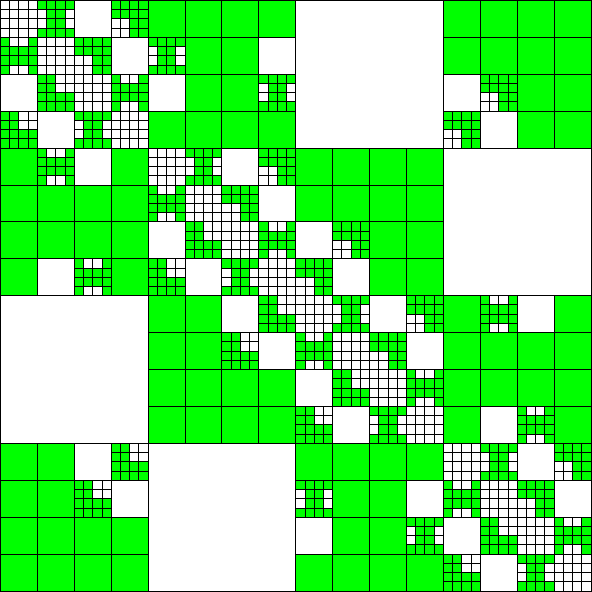}
        }\quad
    \subfloat[]{
    \begin{tikzpicture}
    \node at (-1.5,0) {};
    \node at (-1.5,1.1) {$+$};
    \end{tikzpicture}
    }\quad
    \subfloat[]{
        \includegraphics[scale=0.25]{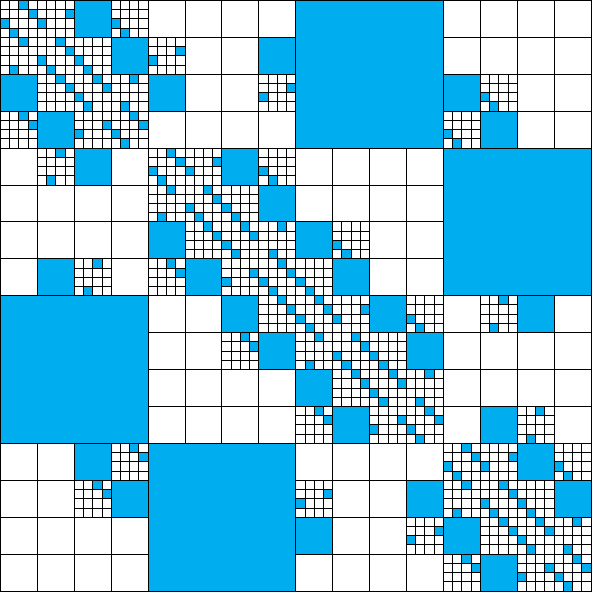}
        }\quad
    \subfloat[]{
        \begin{tikzpicture}
        \node at (-3,0) {};
        \node at (-3,1.1) {$+$};
        \end{tikzpicture}
    }\quad
    \subfloat[]{
        \includegraphics[scale=0.25]{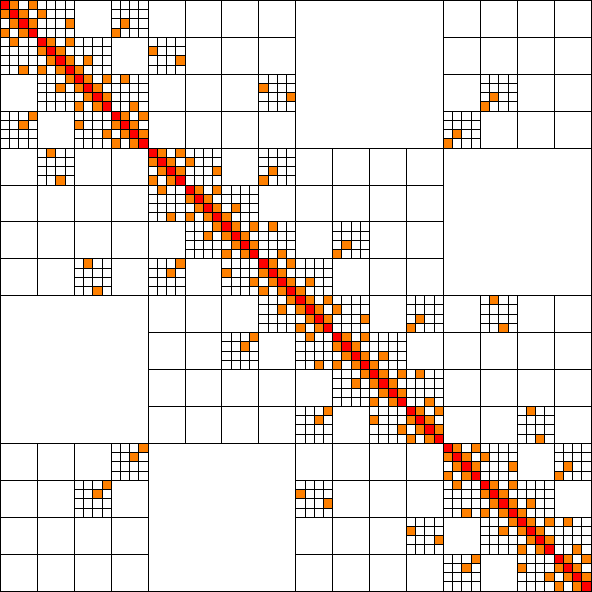}
        }
    \subfloat{
        \begin{tikzpicture}
        \filldraw (2.75,0) to[out=-80,in=80]  (2.75,-2.5) to [out=82,in=-82] cycle;
        \end{tikzpicture}
    }\quad
    \subfloat{
        \begin{tikzpicture}
        \node[anchor=west] at (-1,-2.75){$\times$};
        \draw[draw=black, fill=blue!40] (0, -4) rectangle (0.12, -1.5);
        \node[anchor=west] at (-0.15,-2.75){$\pmb{q}$};
        \end{tikzpicture}
    }\quad
    \caption{In the $\bkt{\mathcal{H}^2 + \mathcal{H}}_{*}$ algorithm, the operators (P$2$M/M$2$M, M$2$L and L$2$L/L$2$P) corresponding to the far-field interaction are constructed using B$2$T NCA and the operators $(U, V)$ corresponding to vertex-sharing interaction are constructed using partially pivoted ACA. After that, we calculate the far-field potential by following the Upward, Transverse and Downward tree traversal. The vertex-sharing potential is calculated independently following a single tree traversal. We get the final potential by adding the far-field, vertex-sharing and near-field potentials.}
    \label{fig:snHODLR2D_construction}
\end{figure}

Therefore, the computed potential is given by $\pmb{\Tilde{\phi}} =\bigg [\phi ^{\bkt{\kappa}}_{X_1};\phi ^{\bkt{\kappa}}_{X_2}; \cdots ;\phi ^{\bkt{\kappa}}_{X_{2^{d \kappa}}} \bigg ]$ (MATLAB notation). The visual representation of the $\bkt{\mathcal{H}^2 + \mathcal{H}}_{*}$ algorithm is given in \Cref{fig:snHODLR2D_construction}.

\subsubsection{Complexity analysis of \texorpdfstring{$\bkt{\mathcal{H}^2 + \mathcal{H}}_{*}$}{H1.5}}
We construct the $\bkt{\mathcal{H}^2 + \mathcal{H}}_{*}$ hierarchical representation using B$2$T NCA and the ACA. 

\textbf{Time complexity.} The time complexity of the initialization steps and the potential calculation steps are given below
\begin{itemize}
    \item \emph{Far-field interaction compression: }  The far-field compression cost of this algorithm is the same as the $\mathcal{H}^2_{*}$(b$+$t) algorithm. Let the leaf cluster size be bounded by $p_1$, and we also assume $p_1 = \mathcal{O} \bkt{n_{max}}$. Let $c_{far}$ be the maximum of the far-field interaction list size (in $2$D and $3$D for this hierarchical representation $c_{far}=12$ and $126$, respectively). Considering all the levels, the cost is bounded by $\dsum_{l=2}^{\kappa} 2^{dl} \bkt{2^{d}p_1 + c_{far}2^{d}p_1}p_1^2 \approx \mathcal{O} \bkt{N}$.
    \item \emph{Vertex-sharing interaction compression: } Let $c_{ver}$ be the maximum of the vertex-sharing interaction list size of a cluster (in $2$D and $3$D for this hierarchical representation $c_{ver}=3$ and $7$, respectively). The partially pivoted ACA \cite{aca} is applied to a cluster's global row and column indices. Let $p_3$ be the maximum rank (vertex-sharing blocks) and $p_3 \in \mathcal{O} \bkt{\log(N) \log^d \bkt{\log(N)}}$ \cite{khan2022numerical}. Then the cost for applying ACA for a particular cluster is bounded by $\bkt{\dfrac{N}{2^{dl}} + c_{ver} \dfrac{N}{2^{dl}}} p_3^2$. By taking all the levels, the cost for the vertex-sharing interaction compression is bounded by $\dsum_{l=1}^{\kappa} 2^{dl} \bkt{\dfrac{N}{2^{dl}} + c_{ver} \dfrac{N}{2^{dl}}} p_3^2 = \dsum_{l=1}^{\kappa} \bkt{N + c_{ver}N} p_3^2 \approx \mathcal{O} \bkt{p_3 N \log \bkt{N}}$.
\end{itemize}
    Therefore, the overall time complexity to initialize the $\bkt{\mathcal{H}^2 + \mathcal{H}}_{*}$ hierarchical representation scales asymptotically $\mathcal{O} \bkt{N + p_3 N \log \bkt{N}}$, $p_3 \in \mathcal{O} \bkt{\log(N) \log^d \bkt{\log(N)}}$ \cite{khan2022numerical}.
    
The time required to compute the far-field and vertex-sharing potentials calculation exhibits a scaling of $\mathcal{O} \bkt{N}$ and $\mathcal{O} \bkt{p_3N \log (N)}$, respectively. Therefore, the overall time complexity for potential calculation (MVP) scales quasi-linearly for the non-oscillatory kernels. 

\textbf{Space complexity.} The cost of storing the far-field P$2$M/M$2$M, M$2$L, and L$2$L/L$2$P operators is $\mathcal{O} (N)$, same as in the $\mathcal{H}^2_{*}$(b$+$t) algorithm. The cost of storing the vertex-sharing $U$ and $V$ operators is roughly $\mathcal{O} \bkt{p_3 N \log(N)}$. So, the overall storage cost of the $\bkt{\mathcal{H}^2 + \mathcal{H}}_{*}$ algorithm scales quasi-linearly.
\begin{remark}
It is noteworthy that the compression routines of the far-field and vertex-sharing interactions are independent and operate separately for both the proposed algorithms. Hence, one has the flexibility to opt for different values of tolerances, i.e., $\epsilon_{far}$ and $\epsilon_{ver}$, to achieve improved compression and, consequently, enhanced relative error control. However, we set $\epsilon_{far} = \epsilon_{ver} = \epsilon$ throughout this article.
\end{remark}
\begin{remark}
Despite having the same asymptotic complexities as the $\mathcal{H}_{*}$ algorithm \cite{khan2022numerical}, the proposed $\bkt{\mathcal{H}^2 + \mathcal{H}}_{*}$ algorithm demonstrates improved time and storage efficiency due to the use of nested bases for partial interaction list. This enhancement is illustrated in numerical results (\Cref{num_results}).
\end{remark}
\begin{remark}
Note that in $1$D ($d=1$), $\mathcal{H}^2_{*}$ $\equiv$ HSS, and $\bkt{\mathcal{H}^2 + \mathcal{H}}_{*}$ $\equiv$ HODLR \cite{ambikasaran2019hodlrlib} $\equiv$ $\mathcal{H}_{*}$.
\end{remark}

\section{Numerical results} \label{num_results}
This section presents various numerical experiments to demonstrate the performance of the proposed algorithms. We also compare the proposed algorithms with other related existing \textbf{algebraic} fast MVP algorithms and present various benchmarks. Notably, all the algorithms are developed using NCA/ACA-based compression.

For convenience, we discuss a few abbreviations in \Cref{tab:nex_notation} corresponding to the different algebraic fast algorithms discussed in this article. The notation b$+$t in the parenthesis signifies the utilization of both the B$2$T NCA and T$2$B NCA. On the other hand, b or t denotes the exclusive application of either B$2$T NCA or T$2$B NCA, respectively.

We perform the following experiments for different algorithms in $2$D and $3$D, report their performance, scalability, and plot various benchmarks.
\begin{enumerate}[label=\textnormal{(\arabic*)}]
    \item Fast MVP for the kernel matrix arises from different kernel functions. We choose the single layer Laplacian $\bkt{\log(r)}$ in $2$D. In $3$D, we select single layer Laplacian $\bkt{1/r}$, Matérn covariance kernel $\bkt{\exp{(-r)}}$ and Helmholtz kernel $\bkt{\exp{(ir)}/r}$.
    \item Fast iterative solver (GMRES) for the Fredholm integral equation of the second kind. 
    \item Fast iterative solver (GMRES) for RBF interpolation.
\end{enumerate}

 \begin{table}[H]
      \centering
     \resizebox{\textwidth}{!}{
    \begin{tabular}{|l|l|l|l|l|}
    \hline
     \textbf{Abbreviation} & \textbf{Description of the fast MVP algorithm} & \textbf{Admissible clusters} & \textbf{Compression technique(s)} &  \textbf{Type of bases} \\ \hline
     $\mathcal{H}^2_{*}$(b$+$t)        & \makecell{The proposed \textbf{efficient} $\mathcal{H}^2_{*}$ algorithm.\\ (\Cref{nHODLRdD})}   & Far-field and Vertex-sharing & B$2$T NCA and T$2$B NCA & Nested  \\ \hline
     $\bkt{\mathcal{H}^2 + \mathcal{H}}_{*}$             & \makecell{The proposed $\bkt{\mathcal{H}^2 + \mathcal{H}}_{*}$ algorithm.\\ (\Cref{s-nHODLRdD})}& Far-field and Vertex-sharing & B$2$T NCA and ACA & Semi-nested \\ \hline
     $\mathcal{H}^2_{*}$(t)            & \makecell{The nested algorithm, where T$2$B NCA is used \\ upon the \textbf{entire} interaction list $\bkt{\mathcal{IL}_*(X)}$. \\(\Cref{nca_t2b_entire_hodlrdd}/ \Cref{t2b_nHODt})} & Far-field and Vertex-sharing & T$2$B NCA  & Nested  \\ \hline
     $\mathcal{H}^2_{\sqrt{d}}$(b) & \makecell{The standard algebraic $\mathcal{H}^2$ matrix algorithm, \\ where B$2$T NCA is used. \\(\Cref{b2t_h2}/ \Cref{b2t_h2matrix})} & Far-field & B$2$T NCA & Nested  \\ \hline
     $\mathcal{H}^2_{\sqrt{d}}$(t) &  \makecell{The standard algebraic $\mathcal{H}^2$ matrix algorithm, \\ where T$2$B NCA is used. \\(\Cref{t2b_h2}/ \Cref{t2b_h2matrix})} & Far-field & T$2$B NCA & Nested  \\ \hline
     $\mathcal{H}_{*}$                & \makecell{The HODLR$d$D algorithm in \cite{kandappan2022hodlr2d, khan2022numerical}.\\ \rk{In this article, we call it $\mathcal{H}_{*}$}.} & Far-field and Vertex-sharing & ACA & Non-nested  \\ \hline
     $\mathcal{H}_{\sqrt{d}}$      & \makecell{The standard algebraic $\mathcal{H}$ matrix algorithm.\\ (\Cref{strong_admiss})} & Far-field & ACA & Non-nested \\ \hline
    \end{tabular}}
     \caption{Abbreviations of the algebraic fast MVP algorithms used in this section for various numerical experiments.}
     \label{tab:nex_notation}
 \end{table}

 It is worth noting that the \emph{initialization time} of the $\mathcal{O}(N)$ NCA discussed in \cite{zhao2019fast, gujjula2022new} is faster than that of the $\mathcal{O}(N \log(N))$ NCA discussed in \cite{bebendorf2012constructing}. However, the storage and MVP time scales are similar. The article \cite{gujjula2022new} also verifies it numerically. Therefore, we believe that it would be sufficient to compare the proposed algorithms with $\mathcal{H}^2_{\sqrt{d}}$(b) and $\mathcal{H}^2_{\sqrt{d}}$(t) from \cite{zhao2019fast}.

 All the algorithms are developed in a similar fashion in $\texttt{C++}$ and performed on an Intel Xeon Gold $2.5$GHz processor with $8$ OpenMP threads within the same environment configuration. We want to emphasize that high-performance implementation is not the goal of this article. The central theme here is to maintain consistency in implementing all the fast algorithms, thereby allowing for meaningful comparisons.
 
 Further, we introduce some other notations in \Cref{tab:app_notations}, which are used in this section.
  \begin{table}[H]
     \centering
     \resizebox{\textwidth}{!}{\begin{tabular}{|c|c|}
     \hline
        $N$ &  Total number of particles (we consider the same source and target) in the domain \\
    \hline
    $\Tilde{K}$ & \makecell{$\Tilde{K}$ represent the hierarchical low-rank representation of the original kernel matrix $K$.} \\
    \hline
    Memory (Mem.) & \makecell{The total memory (in GB) needed to store the hierarchical low-rank representation ($\Tilde{K}$). \\This refers to the total memory for storing all the operators, including the near-field operators.} \\
    \hline
        Initialization time $(t_{init})$ &  \makecell{The time taken (in seconds) to create the hierarchical low-rank representation $\Tilde{K}$ (Initialization routine). \\ This refers to the overall execution time, excluding the time spent on the matrix-vector product operation.}  \\
    \hline
    Fast MVP time $(t_{MVP})$  & \makecell{The time taken (in seconds) to compute the MVP using the hierarchical low-rank representation,\\ i.e., time for \textbf{fast} MVP $(\Tilde{K} \pmb{q})$.}   \\
    \hline
    Total time $(t_{total})$  & \makecell{Total time taken (in seconds), i.e., $t_{total} = t_{init} + t_{mvp}$.}   \\
    \hline
    Direct MVP time $(t_{direct})$  & \makecell{The time taken (in seconds) to compute the \textbf{direct} matrix-vector product, i.e., time for naive MVP.}   \\
    \hline
    Solution time $(t_{sol})$ & \makecell{Total time (in seconds) to solve the system $K \pmb{\sigma} = \pmb{f}$ using the \textbf{fast} iterative solver GMRES \cite{saad1986gmres}. \\ In the GMRES routine, the MVP is accelerated using the hierarchical low-rank representation ($\Tilde{K}$)}\\
    \hline
    $n_{max}$ & \makecell{Maximum number of particles in leaf-clusters. \\ We set $n_{max} = 100$ and $125$ for the experiments in $2$D and $3$D, respectively.} \\
    \hline
    Relative error in MVP $({RE}_{MVP})$ & \makecell{Let $\pmb{q}$ be a random column vector and the exact matrix-vector product $\pmb{\phi}=K \pmb{q}$ (exact up to round-off).  \\ The computed column vector $\Tilde{\pmb{\phi}}$ is given by $\Tilde{\pmb{\phi}} = \Tilde{K} \pmb{q}$. \\ The $2$-norm relative error in the \mvp is given by $\magn{\Tilde{\pmb{\phi}} - \pmb{\phi}}_2 / \magn{\pmb{\phi}}_2$ } \\
    \hline
    Relative error in solution $({RE}_{sol})$ & \makecell{Let $\pmb{q}$ be a random column vector and the exact \mvp $\pmb{b} = K \pmb{q}$ (exact up to round-off). \\ We set $\pmb{b}$ as RHS and solve the system $K \pmb{\lambda} = \pmb{b}$ using the iterative solver, GMRES. \\ The $2$-norm relative error in the solution of the system is given by $\magn{\pmb{\lambda} - \pmb{q}}_2 / \magn{\pmb{q}}_2$ } \\
    \hline
    NCA/ACA tolerance $(\epsilon)$ & \makecell{The NCA/ACA tolerance. We set $ \epsilon = \epsilon_{ver} = \epsilon_{far}$ for all the experiments.} \\
    \hline
    GMRES stopping condition & \makecell{We use GMRES \cite{saad1986gmres} without restart, and GMRES routine will terminate if the relative residual is less than $\epsilon_{GMRES}$. \\ The \mvp part in the GMRES routine is accelerated using various \textbf{fast} MVP algorithms (\Cref{tab:nex_notation}). \\ We set $\epsilon_{GMRES} = 10^{-12}$ and $10^{-10}$ for the experiments in $2$D and $3$D, respectively.} \\
    \hline
    $\#iter$ & \makecell{The number of iterations of the \textbf{fast} iterative solver (GMRES) to reach the stopping condition.} \\
    \hline
     \end{tabular}}
     \caption{Notations used in this section for various numerical experiments.}
     \label{tab:app_notations}
 \end{table}

\subsection{Experiments in two dimensions \texorpdfstring{$(d=2)$}{d=2}}  In this subsection, we perform various experiments in $2$D, and we consider the same source and target for all the experiments.

\subsubsection{Fast MVP in \texorpdfstring{$2$}{2}D}
Let us consider $N$ uniformly distributed particles with location at $\{\pmb{x}_i \}_{i=1}^N$ in the square $[-1,1]^2$. We choose the kernel function as the single layer Laplacian in $2$D. The $(i,j)^{th}$ entry of the kernel matrix $K \in \mathbb{R}^{N \times N}$ arises from the kernel function $\log(r)$ is given by
\begin{equation}
K(i,j) = \begin{cases}
            \log(r) = \log \bkt{\magn{\pmb{x}_i - \pmb{x}_j}_2} & \text{ if } i \neq j\\
             0 & \text{ otherwise}
            \end{cases}    
\end{equation}
We compare the proposed algorithms with the $\mathcal{H}^2_{\sqrt{d}}$(b), $\mathcal{H}_{*}$ and $\mathcal{H}_{\sqrt{d}}$ algorithms. Note that the initialization time of $\mathcal{H}^2_{\sqrt{d}}$(b) is better than $\mathcal{H}^2_{\sqrt{d}}$(t), and the scaling of memory and MVP time are the same for both (refer to \cite{zhao2019fast}). Hence, it is sufficient to compare the proposed algorithms with $\mathcal{H}^2_{\sqrt{d}}$(b), i.e., the better one in terms of initialization time. We choose the non-nested algorithms, $\mathcal{H}_{*}$ and $\mathcal{H}_{\sqrt{d}}$, to demonstrate the improvement over memory (Mem.) and MVP time $(t_{MVP})$ achieved in the nested algorithms.

We randomly select $5$ different column vectors $(\pmb{q})$ and perform the fast MVPs $(\Tilde{K} \pmb{q})$ using $\mathcal{H}^2_{*}$(b$+$t), $\bkt{\mathcal{H}^2 + \mathcal{H}}_{*}$, $\mathcal{H}^2_{\sqrt{d}}$(b), $\mathcal{H}_{*}$ and $\mathcal{H}_{\sqrt{d}}$ algorithms. The average memory (Mem.), initialization time $(t_{init})$, MVP time $(t_{MVP})$, and relative error in MVP $({RE}_{MVP})$ with respect to different tolerances $(\epsilon = 10^{-08}, 10^{-10}$, $10^{-12})$ are tabulated in \Cref{table:smvp2d_lap_1}, \Cref{table:smvp2d_lap_2}, and \Cref{table:smvp2d_lap_3}. We also report the direct MVP time $(t_{direct})$. 
            
\begin{table}[H]
      \centering
      \resizebox{12cm}{!}{%
        \begin{tabular}{l|l|l|l|l|l|l|}
        \cline{2-7}
                                                   & \scriptsize $N$ & \scriptsize Mem. (GB) & \scriptsize $t_{init}$ (s) & \scriptsize $t_{MVP}$ (s) &\scriptsize $t_{direct}$ (s) & \scriptsize Rel. error $({RE}_{MVP})$ \\ \hline \cline{2-7}
        \multicolumn{1}{|l|}{\multirow{4}{*}{\rotatebox[origin=c]{90}{\scriptsize $\mathcal{H}^2_{*}$(b$+$t)}}}   & 102400   &  0.42   &  2.94   & 0.01   & 63.05    & 1.82E-08   \\  
        \multicolumn{1}{|l|}{}                     &  409600  &   1.72     &  12.66     &  0.04   & 1006.78     &  8.04E-08   \\  
        \multicolumn{1}{|l|}{}                     &  1638400 &   6.90     &   61.19    &  0.16   & 16116.5     &  2.01E-07   \\  
        \multicolumn{1}{|l|}{}                     &  6553600  &  27.61    &  296.31    &  0.62   & 257863.7     & 1.03E-06    \\ \hline \hline
        \multicolumn{1}{|l|}{\multirow{4}{*}{\rotatebox[origin=c]{90}{\scriptsize $\bkt{\mathcal{H}^2 + \mathcal{H}}_{*}$}}} &  102400  &  0.66   &  1.83    & 0.01    &  63.05   & 3.91E-08   \\  
        \multicolumn{1}{|l|}{}                     &  409600  &  3.05       &  8.17     &  0.07   & 1006.78     &  4.89E-08   \\  
        \multicolumn{1}{|l|}{}                     &  1638400 &  13.55      &  38.98    &  0.38   & 16116.5     &  5.36E-08   \\  
        \multicolumn{1}{|l|}{}                     &  6553600  & 59.70      &  158.44   &  1.70   & 257863.7    & 1.12E-07    \\ \hline \hline
        \multicolumn{1}{|l|}{\multirow{4}{*}{\rotatebox[origin=c]{90}{\scriptsize $\mathcal{H}^2_{\sqrt{d}}$(b)}}}   & 102400   &  0.56    & 2.56   &  0.02    &  63.05   & 1.36E-08   \\   
        \multicolumn{1}{|l|}{}                     &  409600 &  2.32   & 10.16       &  0.08    &  1006.78   &  1.29E-08    \\  
        \multicolumn{1}{|l|}{}                     &  1638400 & 9.31   &  43.23      &  0.28    &  16116.5   &  8.36E-08    \\  
        \multicolumn{1}{|l|}{}                     &  6553600 & 37.30    & 172.28     &  1.11     & 257863.7    & 3.22E-07     \\ \hline \hline
        \multicolumn{1}{|l|}{\multirow{4}{*}{\rotatebox[origin=c]{90}{\scriptsize $\mathcal{H}_{*}$}}} & 102400  & 1.02    & 3.05    &  0.03   & 63.05    & 8.14E-08  \\    
        \multicolumn{1}{|l|}{}                     &  409600  & 4.92        & 12.81     & 0.15    & 1006.78       &  1.24E-07   \\  
        \multicolumn{1}{|l|}{}                     &  1638400 & 23.22       & 65.12     & 0.72    & 16116.5       &  8.59E-07   \\  
        \multicolumn{1}{|l|}{}                     &  6553600  & 108.23     & 313.25     & 3.59    &  257863.7    &  1.51E-07   \\ \hline \hline
        \multicolumn{1}{|l|}{\multirow{4}{*}{\rotatebox[origin=c]{90}{\scriptsize $\mathcal{H}_{\sqrt{d}}$}}}   & 102400  &  1.26   &  2.92    & 0.04    &  63.05   &  2.07E-09  \\    
        \multicolumn{1}{|l|}{}                     &  409600  &  6.22      &  12.71    &  0.23   &  1006.78    &   3.18E-08  \\  
        \multicolumn{1}{|l|}{}                     &  1638400 &  29.58     &  64.42    &  1.09   &  16116.5    &  6.73E-08   \\  
        \multicolumn{1}{|l|}{}                     &  6553600  & 136.44    &  308.83    & 5.19    & 257863.7     & 5.35E-08    \\ \hline 
        \end{tabular}
}
    \caption{Comparison between the proposed algorithms and other \textbf{algebraic} algorithms as mentioned in \cref{tab:nex_notation} for the single layer Laplacian in $2$D with tolerance $\epsilon = 10^{-08}$. Rel. error refers to $2$-norm relative error in MVP, i.e., ${RE}_{MVP}$ (refer to \Cref{tab:app_notations}).}
    \label{table:smvp2d_lap_1}
\end{table}

\begin{table}[H]
      \centering
      \resizebox{12cm}{!}{%
        \begin{tabular}{l|l|l|l|l|l|l|}
        \cline{2-7}
                                                   & \scriptsize $N$ & \scriptsize Mem. (GB) & \scriptsize $t_{init}$ (s) & \scriptsize $t_{MVP}$ (s) &\scriptsize $t_{direct}$ (s) & \scriptsize Rel. error $({RE}_{MVP})$ \\ \hline \cline{2-7}
        \multicolumn{1}{|l|}{\multirow{4}{*}{\rotatebox[origin=c]{90}{\scriptsize $\mathcal{H}^2_{*}$(b$+$t)}}}   & 102400   &  0.49   & 3.96  & 0.01  &  63.05    & 7.67E-10   \\  
        \multicolumn{1}{|l|}{}                     &  409600  &  2.02      &  15.93     & 0.05    &  1006.78    &  2.40E-09   \\  
        \multicolumn{1}{|l|}{}                     &  1638400 &  8.18      &  74.03     & 0.21    &  16116.5    &  1.37E-09   \\  
        \multicolumn{1}{|l|}{}                     &  6553600  & 32.63     &  392.58    & 0.86    &  257863.7   &  4.81E-09    \\ \hline \hline
        \multicolumn{1}{|l|}{\multirow{4}{*}{\rotatebox[origin=c]{90}{\scriptsize $\bkt{\mathcal{H}^2 + \mathcal{H}}_{*}$}}} &  102400  &  0.77   & 2.61    & 0.02    &  63.05  & 3.41E-11   \\  
        \multicolumn{1}{|l|}{}                     &  409600  &  3.57      &  12.35    &  0.11   & 1006.78     &  7.58E-11   \\  
        \multicolumn{1}{|l|}{}                     &  1638400 &  16.25     &  50.02    &  0.47   & 16116.5     &  1.44E-10   \\  
        \multicolumn{1}{|l|}{}                     &  6553600  & 71.02     &  228.27   &  2.11   & 257863.7    &  8.14E-10    \\ \hline \hline
        \multicolumn{1}{|l|}{\multirow{4}{*}{\rotatebox[origin=c]{90}{\scriptsize $\mathcal{H}^2_{\sqrt{d}}$(b)}}}   & 102400  &  0.62   &  2.71   & 0.02    & 63.05   & 9.37E-11   \\   
        \multicolumn{1}{|l|}{}                     &  409600  &  2.46      & 12.52     & 0.08    & 1006.78     &  3.53E-10   \\  
        \multicolumn{1}{|l|}{}                     &  1638400 &  10.28     & 51.66     & 0.35    & 16116.5     &  3.84E-10   \\  
        \multicolumn{1}{|l|}{}                     &  6553600  & 41.20     & 232.39    & 1.41    & 257863.7    &  2.17E-09   \\ \hline \hline
        \multicolumn{1}{|l|}{\multirow{4}{*}{\rotatebox[origin=c]{90}{\scriptsize $\mathcal{H}_{*}$}}} & 102400  &  1.21   &  4.04   &  0.05   &  63.05   & 6.26E-10  \\    
        \multicolumn{1}{|l|}{}                     &  409600  &  5.88      & 16.44     &  0.21   & 1006.78     &  1.99E-09   \\  
        \multicolumn{1}{|l|}{}                     &  1638400 &  28.55     & 80.12     &  0.99   & 16116.5     &  1.40E-09   \\  
        \multicolumn{1}{|l|}{}                     &  6553600  & 131.8     & 398.17    &  5.11   & 257863.7     & 4.60E-09    \\ \hline \hline
        \multicolumn{1}{|l|}{\multirow{4}{*}{\rotatebox[origin=c]{90}{\scriptsize $\mathcal{H}_{\sqrt{d}}$}}}   & 102400  & 1.47    &  4.01   &  0.06   & 63.05    & 1.98E-11  \\    
        \multicolumn{1}{|l|}{}                     &  409600  &  7.34      & 16.87     &  0.23   & 1006.78     & 2.08E-11    \\  
        \multicolumn{1}{|l|}{}                     &  1638400 &  35.09     & 82.23     &  1.21   & 16116.5     & 8.64E-11    \\  
        \multicolumn{1}{|l|}{}                     &  6553600  & 162.62    & 395.71    &  6.08   & 257863.7     & 5.31E-10     \\ \hline 
        \end{tabular}
}
    \caption{Comparison between the proposed algorithms and other \textbf{algebraic} algorithms as mentioned in \cref{tab:nex_notation} for the single layer Laplacian in $2$D with tolerance $\epsilon = 10^{-10}$. Rel. error refers to $2$-norm relative error in MVP, i.e., ${RE}_{MVP}$ (refer to \Cref{tab:app_notations}).}
    \label{table:smvp2d_lap_2}
\end{table}

\begin{table}[H]
      \centering
      \resizebox{12cm}{!}{%
        \begin{tabular}{l|l|l|l|l|l|l|}
        \cline{2-7}
                                                   & \scriptsize $N$ & \scriptsize Mem. (GB) & \scriptsize $t_{init}$ (s) & \scriptsize $t_{MVP}$ (s) &\scriptsize $t_{direct}$ (s) & \scriptsize Rel. error $({RE}_{MVP})$ \\ \hline \cline{2-7}
        \multicolumn{1}{|l|}{\multirow{4}{*}{\rotatebox[origin=c]{90}{\scriptsize $\mathcal{H}^2_{*}$(b$+$t)}}}   & 102400   &  0.56   &  4.94   & 0.02   &  63.05   &  1.63E-12  \\  
        \multicolumn{1}{|l|}{}                     &  409600  &  2.33      & 22.33     &  0.06    &  1006.78    &  1.06E-11   \\  
        \multicolumn{1}{|l|}{}                     &  1638400 &  9.44      & 102.08     & 0.25    &  16116.5    &  3.93E-11   \\  
        \multicolumn{1}{|l|}{}                     &  6553600  & 37.75     & 495.43     & 1.01     &  257863.7    &  1.15E-10   \\ \hline \hline
        \multicolumn{1}{|l|}{\multirow{4}{*}{\rotatebox[origin=c]{90}{\scriptsize $\bkt{\mathcal{H}^2 + \mathcal{H}}_{*}$}}}   & 102400  &  0.88  &  3.36   &  0.04   & 63.05    & 7.86E-13   \\   
        \multicolumn{1}{|l|}{}                     &  409600  &  4.12      &  15.65    &  0.2   &   1006.78    & 1.36E-12    \\  
        \multicolumn{1}{|l|}{}                     &  1638400 &  18.89     &  68.01    &  1.01   &  16116.5    & 1.36E-11    \\  
        \multicolumn{1}{|l|}{}                     &  6553600  & 86.49     &  288.72   &  4.50   &  257863.7   & 2.49E-11    \\ \hline \hline
        \multicolumn{1}{|l|}{\multirow{4}{*}{\rotatebox[origin=c]{90}{\scriptsize $\mathcal{H}^2_{\sqrt{d}}$(b)}}} &  102400  & 0.69    &  3.25   & 0.02   &  63.05   &  7.09E-12  \\  
        \multicolumn{1}{|l|}{}                     &  409600  &  2.82        &  15.46    & 0.09    & 1006.78     & 2.93E-11    \\  
        \multicolumn{1}{|l|}{}                     &  1638400 &  11.21       &  69.47    & 0.39    & 16116.5     &  3.34E-11   \\  
        \multicolumn{1}{|l|}{}                     &  6553600  & 45.99       &  294.3    &  1.64   & 257863.7     &  1.13E-11   \\ \hline \hline
        \multicolumn{1}{|l|}{\multirow{4}{*}{\rotatebox[origin=c]{90}{\scriptsize $\mathcal{H}_{*}$}}} & 102400  &  1.41  & 3.57    &  0.08   & 63.05   &  1.76E-11  \\    
        \multicolumn{1}{|l|}{}                     &  409600  &  6.91      &  21.19    &  0.41    &  1006.78    &  2.63E-11   \\  
        \multicolumn{1}{|l|}{}                     &  1638400 &  33.07     &  109.33    & 2.05    & 16116.5     &  7.37E-11   \\  
        \multicolumn{1}{|l|}{}                     &  6553600  &   -     &   -   & -    &  -    &  -   \\ \hline \hline
        \multicolumn{1}{|l|}{\multirow{4}{*}{\rotatebox[origin=c]{90}{\scriptsize $\mathcal{H}_{\sqrt{d}}$}}}   & 102400  & 1.68   & 3.66    &  0.10   &  63.05   & 1.35E-13  \\    
        \multicolumn{1}{|l|}{}                     &  409600  &  8.41        &  26.84    &  0.42   &   1006.78   & 1.43E-13    \\  
        \multicolumn{1}{|l|}{}                     &  1638400 & 40.47        &  121.10    & 2.32    &  16116.5    & 6.69E-12     \\  
        \multicolumn{1}{|l|}{}                     &  6553600  &  -      &   -   &   -  &    -  &    - \\ \hline
        \end{tabular}
}
    \caption{Comparison between the proposed algorithms and other \textbf{algebraic} algorithms as mentioned in \cref{tab:nex_notation} for the single layer Laplacian in $2$D with tolerance $\epsilon = 10^{-12}$. Rel. error refers to $2$-norm relative error in MVP, i.e., ${RE}_{MVP}$ (refer to \Cref{tab:app_notations}).}
    \label{table:smvp2d_lap_3}
\end{table}

The job either exceeds the specified wall time or is out of the memory corresponding to the cells with $``-"$ in the tables.

From the \Cref{table:smvp2d_lap_1}, \Cref{table:smvp2d_lap_2}, and \Cref{table:smvp2d_lap_3}, the following conclusions can be drawn.

\begin{enumerate}
    \item The scaling of the memory (Mem.) and MVP time $(t_{MVP})$ of the $\mathcal{H}^2_{*}$(b$+$t) algorithm is similar to that of $\mathcal{H}^2_{\sqrt{d}}$(b). Notably, $\mathcal{H}^2_{*}$(b$+$t) shows slightly less storage usage and improved MVP time.
    \item In $2$D, the initialization time $(t_{init})$ of $\mathcal{H}^2_{*}$(b$+$t) is slightly higher than that of $\mathcal{H}^2_{\sqrt{d}}$(b). This is due to the fact that in $\mathcal{H}^2_{*}$(b$+$t), the \emph{vertex-sharing} interactions are compressed using T$2$B NCA, which scales quasi-linearly, while in $\mathcal{H}^2_{\sqrt{d}}$(b), the initialization time scales roughly $\mathcal{O}(N)$. Overall, the initialization time of the $\mathcal{H}^2_{*}$(b$+$t) scales quasi-linearly.
    \item $\mathcal{H}^2_{*}$(b$+$t) and $\mathcal{H}^2_{\sqrt{d}}$(b) are the nested versions of $\mathcal{H}_{*}$ and $\mathcal{H}_{\sqrt{d}}$ algorithms, respectively. We can see that the nested versions require much less storage and improved MVP time than the non-nested ones.
    \item The initialization time of $\bkt{\mathcal{H}^2 + \mathcal{H}}_{*}$ is the fastest and scales similar to the initialization time of $\mathcal{H}^2_{\sqrt{d}}$(b).
    \item The \emph{semi-nested} $\bkt{\mathcal{H}^2 + \mathcal{H}}_{*}$ is better than the $\mathcal{H}_{*}$ and $\mathcal{H}_{\sqrt{d}}$ matrix algorithms due to the use of semi/partially nested bases. Thus, $\bkt{\mathcal{H}^2 + \mathcal{H}}_{*}$ demonstrates a nice balance between nested and non-nested algorithms in $2$D.
\end{enumerate}

We arrange the above fast algorithms in ascending order with respect to memory usage and MVP time.\\
\textbf{Memory:} $\mathcal{H}^2_{*}$(b$+$t) < $\mathcal{H}^2_{\sqrt{d}}$(b) < $\bkt{\mathcal{H}^2 + \mathcal{H}}_{*}$ < $\mathcal{H}_{*}$ < $\mathcal{H}_{\sqrt{d}}$ \\ 
\textbf{MVP time:} $\mathcal{H}^2_{*}$(b$+$t) < $\mathcal{H}^2_{\sqrt{d}}$(b) < $\bkt{\mathcal{H}^2 + \mathcal{H}}_{*}$ < $\mathcal{H}_{*}$ < $\mathcal{H}_{\sqrt{d}}$

Therefore, though we do not have a strict theoretical error bound, the above numerical experiment demonstrates that the proposed $\mathcal{H}^2_{*}$(b$+$t) algorithm shows slightly less memory usage and improved MVP time than the $\mathcal{H}^2_{\sqrt{d}}$ matrix algorithm with comparable accuracy.

\subsubsection{Fast MVP accelerated GMRES for integral equation in \texorpdfstring{$2$}{2}D} 
\begin{figure}[H]
    \centering
    \subfloat[]{\includegraphics[height=3.8cm, width=5.5cm]{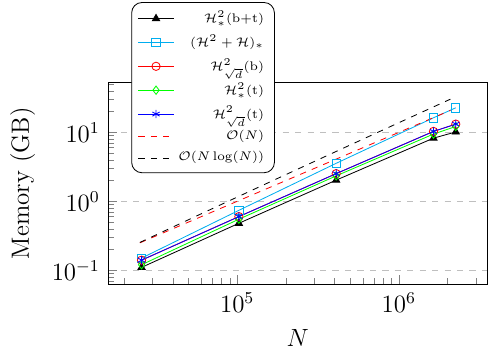}\label{green2d_mem}}%
    \subfloat[]{\includegraphics[height=3.8cm, width=5.5cm]{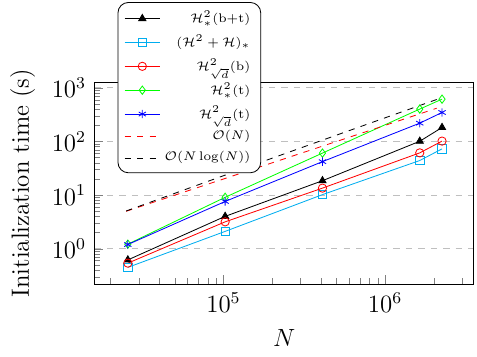}\label{green2d_assem}}%
    \subfloat[]{\includegraphics[height=3.8cm, width=5.5cm]{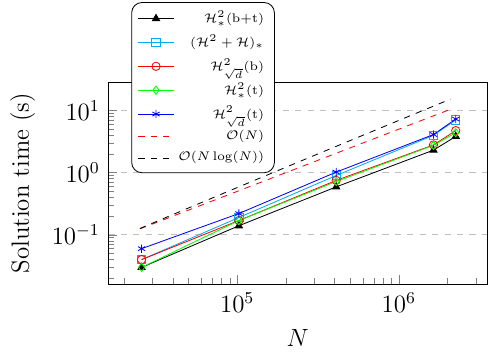}\label{green2d_sol}}%
    \caption{Plots of Memory, Initialization time and Solution time of the fast MVP algorithms.}
    \label{fig:green2d}
\end{figure}

In this experiment, Our goal is to compare the performance of the proposed algorithms with various nested algorithms from \Cref{tab:nex_notation}. Consider the Fredholm integral equation of the second kind over $C = [-1,1]^2 \subset \mathbb{R}^2$, which is given by
\begin{equation} \label{ieq2d_1}
    \sigma(\pmb{x})+\int_{C} F(\pmb{x},\pmb{y})\sigma(\pmb{y})d\pmb{y} = f(\pmb{x}) \qquad \pmb{x}, \pmb{y} \in C
\end{equation}
with $F(\pmb{x},\pmb{y}) = -\dfrac{1}{2\pi} \log{\bkt{\magn{\pmb{x}-\pmb{y}}_2}}$ (Green's function for Laplace equation in $2$D). We follow a piece-wise constant collocation method with collocation points on a uniform grid in $C = [-1,1]^2$ to discretize \cref{ieq2d_1} and set the RHS $\pmb{f}$ as described in the $11^{th}$ cell of \Cref{tab:app_notations}. Therefore, the following linear system can be obtained
\begin{equation} \label{ieq2d_2}
    K \pmb{\sigma} = \pmb{f}
\end{equation}
The \Cref{ieq2d_2} is solved using \textbf{fast} GMRES, i.e., the MVP part of GMRES is accelerated using the $\mathcal{H}^2_{*}$(b$+$t), $\bkt{\mathcal{H}^2 + \mathcal{H}}_{*}$, $\mathcal{H}^2_{\sqrt{d}}$(b), $\mathcal{H}^2_{*}$(t) and $\mathcal{H}^2_{\sqrt{d}}$(t) algorithms. We set $\epsilon = 10^{-10}$ and $\epsilon_{GMRES} = 10^{-12}$ and report the memory, initialization time, solution time and number of iterations in \Cref{table:green2d_int_eq}. The relative error in solution $(RE_{sol})$ is of order $10^{-12}$ in all cases. We also plot the memory (\Cref{green2d_mem}), initialization time (\Cref{green2d_assem}) and solution time (\Cref{green2d_sol}) in \Cref{fig:green2d}.

\begin{table}[H]
      \centering
      \resizebox{\textwidth}{!}{%
      \setlength\extrarowheight{0.9pt}
\begin{tabular}{|l|lllll|lllll|lllll|l|}
\hline
\multirow{2}{*}{N} & \multicolumn{5}{l|}{\qquad \qquad \qquad Memory (GB)}                                                                                           & \multicolumn{5}{l|}{\qquad \qquad \quad Initialization time (s)}                                                                                & \multicolumn{5}{l|}{\qquad \qquad \quad Solution time (s)}                                                                                      & \multirow{2}{*}{$\#iter$} \\ \cline{2-16}
                   & \multicolumn{1}{l|}{\scriptsize $\mathcal{H}^2_{*}$(b$+$t)} & \multicolumn{1}{l|}{\scriptsize $\bkt{\mathcal{H}^2 + \mathcal{H}}_{*}$} & \multicolumn{1}{l|}{\scriptsize $\mathcal{H}^2_{\sqrt{d}}$(b)} & \multicolumn{1}{l|}{\scriptsize $\mathcal{H}^2_{*}$(t)} & \scriptsize $\mathcal{H}^2_{\sqrt{d}}$(t) & \multicolumn{1}{l|}{\scriptsize $\mathcal{H}^2_{*}$(b$+$t)} & \multicolumn{1}{l|}{\scriptsize $\bkt{\mathcal{H}^2 + \mathcal{H}}_{*}$} & \multicolumn{1}{l|}{\scriptsize $\mathcal{H}^2_{\sqrt{d}}$(b)} & \multicolumn{1}{l|}{\scriptsize $\mathcal{H}^2_{*}$(t)} & \scriptsize $\mathcal{H}^2_{\sqrt{d}}$(t) & \multicolumn{1}{l|}{\scriptsize $\mathcal{H}^2_{*}$(b$+$t)} & \multicolumn{1}{l|}{\scriptsize $\bkt{\mathcal{H}^2 + \mathcal{H}}_{*}$} & \multicolumn{1}{l|}{\scriptsize $\mathcal{H}^2_{\sqrt{d}}$(b)} & \multicolumn{1}{l|}{\scriptsize $\mathcal{H}^2_{*}$(t)} & \scriptsize $\mathcal{H}^2_{\sqrt{d}}$(t) &                       \\ \hline 
                   25600& \multicolumn{1}{l|}{0.11}     & \multicolumn{1}{l|}{0.15}       & \multicolumn{1}{l|}{0.14}       & \multicolumn{1}{l|}{0.12}        &  0.14    & \multicolumn{1}{l|}{0.62}     & \multicolumn{1}{l|}{0.45}       & \multicolumn{1}{l|}{0.54}      & \multicolumn{1}{l|}{1.22}             &  1.2    & \multicolumn{1}{l|}{0.03}     & \multicolumn{1}{l|}{0.04}       & \multicolumn{1}{l|}{0.04}      & \multicolumn{1}{l|}{0.03}        &   0.06    &   6                    \\ \hline
                   102400& \multicolumn{1}{l|}{0.48}     & \multicolumn{1}{l|}{0.74}       & \multicolumn{1}{l|}{0.60}      & \multicolumn{1}{l|}{0.55}        &  0.60   & \multicolumn{1}{l|}{4.02}     & \multicolumn{1}{l|}{2.11}       & \multicolumn{1}{l|}{3.18}      & \multicolumn{1}{l|}{9.11}             &  7.625  & \multicolumn{1}{l|}{0.14}     & \multicolumn{1}{l|}{0.19}       & \multicolumn{1}{l|}{0.17}      & \multicolumn{1}{l|}{0.17}        &  0.22     &   6                    \\ \hline
                   409600& \multicolumn{1}{l|}{2.03}     & \multicolumn{1}{l|}{3.56}       & \multicolumn{1}{l|}{2.52}      & \multicolumn{1}{l|}{2.34}        &  2.52   & \multicolumn{1}{l|}{18.57}     & \multicolumn{1}{l|}{10.22}       & \multicolumn{1}{l|}{13.51}      & \multicolumn{1}{l|}{60.79}          & 42.11   & \multicolumn{1}{l|}{0.59}     & \multicolumn{1}{l|}{0.90}       & \multicolumn{1}{l|}{0.75}      & \multicolumn{1}{l|}{0.72}        &  1.02     &   6                    \\ \hline
                   1638400& \multicolumn{1}{l|}{8.25}     & \multicolumn{1}{l|}{16.29}       & \multicolumn{1}{l|}{10.33}       & \multicolumn{1}{l|}{9.45}     &   10.33 & \multicolumn{1}{l|}{101.05}     & \multicolumn{1}{l|}{44.2}       & \multicolumn{1}{l|}{61.39}      & \multicolumn{1}{l|}{402.29}          & 219.66  & \multicolumn{1}{l|}{2.29}     & \multicolumn{1}{l|}{3.98}       & \multicolumn{1}{l|}{2.79}      & \multicolumn{1}{l|}{2.75}        &  4.10     &   6                    \\ \hline
                   2250000& \multicolumn{1}{l|}{10.15}     & \multicolumn{1}{l|}{22.74}       & \multicolumn{1}{l|}{13.18}      & \multicolumn{1}{l|}{12.19}   &  13.17  & \multicolumn{1}{l|}{182.42}     & \multicolumn{1}{l|}{72.85}       & \multicolumn{1}{l|}{101.67}      & \multicolumn{1}{l|}{615.54}        &  347.74 & \multicolumn{1}{l|}{3.82}     & \multicolumn{1}{l|}{7.02}       & \multicolumn{1}{l|}{4.74}      & \multicolumn{1}{l|}{4.39}        & 7.27       &   6                    \\ \hline
\end{tabular}
}
    \caption{Performance of the proposed and the nested algorithms as mentioned in \Cref{tab:nex_notation}. We set $\epsilon = 10^{-10}$ and $\epsilon_{GMRES} = 10^{-12}$.}
    \label{table:green2d_int_eq}
\end{table}

Let us arrange the fast algorithms in ascending order with respect to memory usage and solution time.\\
\textbf{Memory:} $\mathcal{H}^2_{*}$(b$+$t) < $\mathcal{H}^2_{*}$(t) < $\mathcal{H}^2_{\sqrt{d}}$(b) $\leq$ $\mathcal{H}^2_{\sqrt{d}}$(t) < $\bkt{\mathcal{H}^2 + \mathcal{H}}_{*}$ \\ 
\textbf{Solution time:} $\mathcal{H}^2_{*}$(b$+$t) < $\mathcal{H}^2_{*}$(t) < $\mathcal{H}^2_{\sqrt{d}}$(b) < $\mathcal{H}^2_{\sqrt{d}}$(t) < $\bkt{\mathcal{H}^2 + \mathcal{H}}_{*}$

\subsubsection{Fast MVP accelerated GMRES for RBF interpolation in \texorpdfstring{$2$}{2}D}
\begin{figure}[H]
    \centering
    \subfloat[]{\includegraphics[height=3.8cm, width=5.5cm]{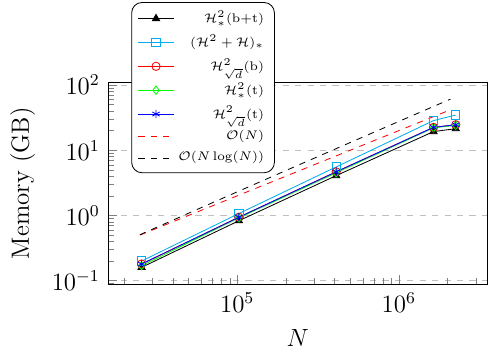}\label{rbf2d_a_r_mem}}%
    \subfloat[]{\includegraphics[height=3.8cm, width=5.5cm]{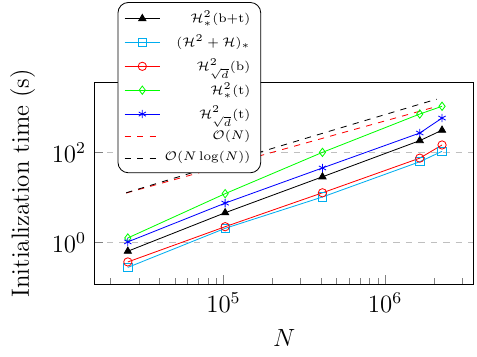}\label{rbf2d_a_r_assem}}%
    \subfloat[]{\includegraphics[height=3.8cm, width=5.5cm]{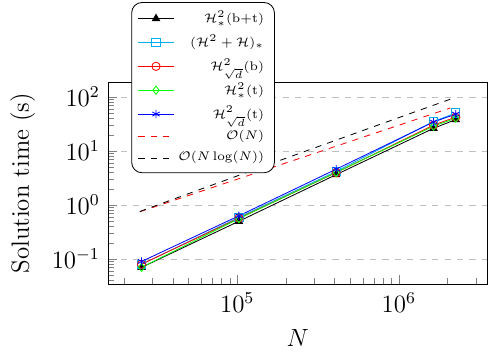}\label{rbf2d_a_r_sol}}%
    \caption{Plots of Memory, Initialization time and Solution time of the fast MVP algorithms.}
    \label{fig:rbf2d_a_r}
\end{figure}

Let the location of the particles $\{\pmb{x}_i \}_{i=1}^N$ be the $\sqrt{N} \times \sqrt{N}$ Chebyshev grid on the domain $[-1,1]^2$. We consider the Chebyshev distribution of particles to study the performance of the nested algorithms over slightly non-uniformly distributed particles. However, we use the uniform $2^d$ tree (quad tree in $2$D) as described in \Cref{tree_construction}. Let us consider the following radial basis function
\begin{equation}
 G(r) = \begin{cases}
            a/r & \text{ if } r \geq a\\
             r/a & \text{ if } r < a
            \end{cases}   
\end{equation}
Below is the dense linear system generated by the radial basis function $G$
\begin{equation} \label{rbf2d_eq2}
    \alpha \lambda_i + \dsum_{j=1, j \neq i}^N G \bkt{\magn{\pmb{x}_i - \pmb{x}_j}_2} \lambda_j = b_i, \qquad i = 1,2,\dots, N.
\end{equation}
We set $a=0.0001$ and $\alpha=N^{1/4}$. By setting $\pmb{b}$ as described in the $11^{th}$ cell of \Cref{tab:app_notations}, the \Cref{rbf2d_eq2} can be written in the form 
\begin{equation} \label{rbf2d_eq3}
    K \pmb{\lambda} = \pmb{b}
\end{equation}
The \Cref{rbf2d_eq3} is solved using \textbf{fast} GMRES. We set $\epsilon = 10^{-10}$ and $\epsilon_{GMRES} = 10^{-12}$ and report the memory, initialization time and solution time in \Cref{table:rbf2d_a_r_a} for all the nested algorithms. The relative error in solution $(RE_{sol})$ is of order $10^{-9}$ in all cases. We also plot the memory (\Cref{rbf2d_a_r_mem}), initialization time (\Cref{rbf2d_a_r_assem}) and solution time (\Cref{rbf2d_a_r_sol}) in \Cref{fig:rbf2d_a_r}.

\begin{table}[H]
      \centering
      \resizebox{\textwidth}{!}{%
      \setlength\extrarowheight{0.9pt}
\begin{tabular}{|l|lllll|lllll|lllll|l|}
\hline
\multirow{2}{*}{N} & \multicolumn{5}{l|}{\qquad \qquad \qquad Memory (GB)}                                                                                           & \multicolumn{5}{l|}{\qquad \qquad \quad Initialization time (s)}                                                                                & \multicolumn{5}{l|}{\qquad \qquad \quad Solution time (s)}                                                                                      & \multirow{2}{*}{$\#iter$} \\ \cline{2-16}
                   & \multicolumn{1}{l|}{\scriptsize $\mathcal{H}^2_{*}$(b$+$t)} & \multicolumn{1}{l|}{\scriptsize $\bkt{\mathcal{H}^2 + \mathcal{H}}_{*}$} & \multicolumn{1}{l|}{\scriptsize $\mathcal{H}^2_{\sqrt{d}}$(b)} & \multicolumn{1}{l|}{\scriptsize $\mathcal{H}^2_{*}$(t)} & \scriptsize $\mathcal{H}^2_{\sqrt{d}}$(t) & \multicolumn{1}{l|}{\scriptsize $\mathcal{H}^2_{*}$(b$+$t)} & \multicolumn{1}{l|}{\scriptsize $\bkt{\mathcal{H}^2 + \mathcal{H}}_{*}$} & \multicolumn{1}{l|}{\scriptsize $\mathcal{H}^2_{\sqrt{d}}$(b)} & \multicolumn{1}{l|}{\scriptsize $\mathcal{H}^2_{*}$(t)} & \scriptsize $\mathcal{H}^2_{\sqrt{d}}$(t) & \multicolumn{1}{l|}{\scriptsize $\mathcal{H}^2_{*}$(b$+$t)} & \multicolumn{1}{l|}{\scriptsize $\bkt{\mathcal{H}^2 + \mathcal{H}}_{*}$} & \multicolumn{1}{l|}{\scriptsize $\mathcal{H}^2_{\sqrt{d}}$(b)} & \multicolumn{1}{l|}{\scriptsize $\mathcal{H}^2_{*}$(t)} & \scriptsize $\mathcal{H}^2_{\sqrt{d}}$(t) &                       \\ \hline 
                   25600& \multicolumn{1}{l|}{0.16}     & \multicolumn{1}{l|}{0.20}       & \multicolumn{1}{l|}{0.18}      & \multicolumn{1}{l|}{0.17}        &    0.18     & \multicolumn{1}{l|}{0.64}     & \multicolumn{1}{l|}{0.28}       & \multicolumn{1}{l|}{0.37}      & \multicolumn{1}{l|}{1.26}        &  1.03       & \multicolumn{1}{l|}{0.07}     & \multicolumn{1}{l|}{0.08}       & \multicolumn{1}{l|}{0.08}      & \multicolumn{1}{l|}{0.07}        &   0.09    &    9                  \\ \hline
                   102400& \multicolumn{1}{l|}{0.83}     & \multicolumn{1}{l|}{2.06}       & \multicolumn{1}{l|}{0.94}      & \multicolumn{1}{l|}{0.92}        &  0.94     & \multicolumn{1}{l|}{4.57}     & \multicolumn{1}{l|}{2.06}       & \multicolumn{1}{l|}{2.24}      & \multicolumn{1}{l|}{12.12}        & 7.48       & \multicolumn{1}{l|}{0.50}     & \multicolumn{1}{l|}{0.59}       & \multicolumn{1}{l|}{0.57}      & \multicolumn{1}{l|}{0.56}        &  0.63     &  13                     \\ \hline
                   409600& \multicolumn{1}{l|}{4.08}     & \multicolumn{1}{l|}{5.55}       & \multicolumn{1}{l|}{4.67}      & \multicolumn{1}{l|}{4.52}        &  4.67     & \multicolumn{1}{l|}{28.68}     & \multicolumn{1}{l|}{10.22}       & \multicolumn{1}{l|}{12.67}      & \multicolumn{1}{l|}{100.06}      &  45.09    & \multicolumn{1}{l|}{3.65}     & \multicolumn{1}{l|}{4.18}       & \multicolumn{1}{l|}{3.98}      & \multicolumn{1}{l|}{4.01}        &  4.60     &  19                     \\ \hline
                   1638400& \multicolumn{1}{l|}{19.31}     & \multicolumn{1}{l|}{28.47}       & \multicolumn{1}{l|}{22.43}      & \multicolumn{1}{l|}{22.09}   &  22.41    & \multicolumn{1}{l|}{182.27}     & \multicolumn{1}{l|}{61.95}       & \multicolumn{1}{l|}{74.50}      & \multicolumn{1}{l|}{709.76}        &  268.86 & \multicolumn{1}{l|}{26.27}     & \multicolumn{1}{l|}{35.09}      & \multicolumn{1}{l|}{30.11}      & \multicolumn{1}{l|}{28.78}        & 34.42      &  29                     \\ \hline
                   2250000& \multicolumn{1}{l|}{21.43}     & \multicolumn{1}{l|}{34.64}       & \multicolumn{1}{l|}{24.61}      & \multicolumn{1}{l|}{24.12}   & 24.60      & \multicolumn{1}{l|}{310.58}     & \multicolumn{1}{l|}{105.48}       & \multicolumn{1}{l|}{145.85}      & \multicolumn{1}{l|}{1054.3}      & 577.35  & \multicolumn{1}{l|}{38.15}     & \multicolumn{1}{l|}{50.57}      & \multicolumn{1}{l|}{42.32}      & \multicolumn{1}{l|}{40.55}        & 48.12      &   33                    \\ \hline
\end{tabular}
}
    \caption{Performance of the proposed and the nested algorithms as mentioned in \Cref{tab:nex_notation}. We set $\epsilon = 10^{-10}$ and $\epsilon_{GMRES} = 10^{-12}$.}
    \label{table:rbf2d_a_r_a}
\end{table}

Let us arrange the fast algorithms in ascending order with respect to memory usage and solution time.\\
\textbf{Memory:} $\mathcal{H}^2_{*}$(b$+$t) < $\mathcal{H}^2_{*}$(t) < $\mathcal{H}^2_{\sqrt{d}}$(b) $\leq$ $\mathcal{H}^2_{\sqrt{d}}$(t) < $\bkt{\mathcal{H}^2 + \mathcal{H}}_{*}$ \\ 
\textbf{Solution time:} $\mathcal{H}^2_{*}$(b$+$t) < $\mathcal{H}^2_{*}$(t) < $\mathcal{H}^2_{\sqrt{d}}$(b) < $\mathcal{H}^2_{\sqrt{d}}$(t) < $\bkt{\mathcal{H}^2 + \mathcal{H}}_{*}$

\emph{\textbf{Summary of the experiments in $2$D.}} Here is a summary of the results obtained from all the experiments performed in $2$D.
\begin{enumerate}
    \item In all the experiments, the proposed $\mathcal{H}^2_{*}$(b$+$t) outperforms all the other algebraic algorithms (\Cref{tab:nex_notation}) in terms of memory and MVP/Solution time.
    \item The initialization time of the $\mathcal{H}^2_{*}$(t) algorithm (\Cref{nca_t2b_entire_hodlrdd}) is the highest and the proposed $\mathcal{H}^2_{*}$(b$+$t) performs way better than the $\mathcal{H}^2_{*}$(t) algorithm. 
    \item The $\bkt{\mathcal{H}^2 + \mathcal{H}}_{*}$ is the fastest among the other algorithms in terms of the initialization time and it performs better than the non-nested algorithms.
    \item The initialization time of $\mathcal{H}^2_{\sqrt{d}}$(b) is better than the $\mathcal{H}^2_{\sqrt{d}}$(t), and the scaling of the memory and MVP time are almost the same, which was shown in \cite{zhao2019fast} also.
\end{enumerate}

\subsection{Experiments in three dimensions \texorpdfstring{$(d=3)$}{d=3}} In this subsection, we perform various experiments in $3$D. We consider the same source and target for all the experiments. 
\subsubsection{Fast MVP in \texorpdfstring{$3$}{3}D}
Let us consider $N$ uniformly distributed particles with location at $\{\pmb{x}_i \}_{i=1}^N$ in the cube $[-1,1]^3$. We choose the kernel function as the single layer Laplacian in $3$D. The $(i,j)^{th}$ entry of the kernel matrix $K \in \mathbb{R}^{N \times N}$ arises from the kernel function $1/r$ is given by
\begin{equation}
K(i,j) = \begin{cases}
            \dfrac{1}{r} = \dfrac{1}{\magn{\pmb{x}_i - \pmb{x}_j}_2} & \text{ if } i \neq j\\
             0 & \text{ otherwise}
            \end{cases}    
\end{equation}
We compare the proposed algorithms with the $\mathcal{H}^2_{\sqrt{d}}$(b) (the better $\mathcal{H}^2$ matrix algorithm in terms of initialization time), $\mathcal{H}_{*}$ and $\mathcal{H}_{\sqrt{d}}$ algorithms. Again, we choose the non-nested algorithms, $\mathcal{H}_{*}$ and $\mathcal{H}_{\sqrt{d}}$, to demonstrate the improvement over memory (Mem.) and MVP time $(t_{MVP})$ achieved in the nested algorithms.

We randomly select $5$ different column vectors $(\pmb{q})$ and perform the fast MVPs $(\Tilde{K} \pmb{q})$ using $\mathcal{H}^2_{*}$(b$+$t), $\bkt{\mathcal{H}^2 + \mathcal{H}}_{*}$, $\mathcal{H}^2_{\sqrt{d}}$(b), $\mathcal{H}_{*}$ and $\mathcal{H}_{\sqrt{d}}$ algorithms. The average memory (Mem.), initialization time $(t_{init})$, MVP time $(t_{MVP})$, and relative error in MVP $({RE}_{MVP})$ with respect to different tolerances $(\epsilon = 10^{-4}, 10^{-6}, 10^{-8})$ are tabulated in \Cref{table:smvp3d_lap_1}, \Cref{table:smvp3d_lap_2}, and \Cref{table:smvp3d_lap_3}. We also report the direct MVP time $(t_{direct})$. 

\begin{table}[H]
      \centering
      \resizebox{12cm}{!}{%
        \begin{tabular}{l|l|l|l|l|l|l|}
        \cline{2-7}
                                                   & \scriptsize $N$ & \scriptsize Mem. (GB) & \scriptsize $t_{init}$ (s) & \scriptsize $t_{MVP}$ (s) &\scriptsize $t_{direct}$ (s) & \scriptsize Rel. error $({RE}_{MVP})$ \\ \hline \cline{2-7}
        \multicolumn{1}{|l|}{\multirow{4}{*}{\rotatebox[origin=c]{90}{\scriptsize $\mathcal{H}^2_{*}$(b$+$t)}}}   & 64000  &  0.91  &  2.76  & 0.03   &  13.18  &  1.94E-04 \\  
        \multicolumn{1}{|l|}{}                     &  512000  &  9.18      &  51.78     &  0.41   &  851.35      &  3.85E-04   \\  
        \multicolumn{1}{|l|}{}                     &  1000000 &  21.36     &  131.65    &  1.03   &  3829.64     & 7.15E-04    \\  
        \multicolumn{1}{|l|}{}                     &  1728000 &  39.07     &  237.65    &  2.25   &  11575.3     & 1.08E-03    \\ \hline \hline
        \multicolumn{1}{|l|}{\multirow{4}{*}{\rotatebox[origin=c]{90}{\scriptsize $\bkt{\mathcal{H}^2 + \mathcal{H}}_{*}$}}} &  64000  &   1.04 &  1.56  &  0.04    &  13.18    &   8.41E-05 \\  
        \multicolumn{1}{|l|}{}                     & 512000   &   10.68     &  37.81     & 0.57    &  851.35      & 1.24E-04    \\ 
        \multicolumn{1}{|l|}{}                     & 1000000  &   25.41     &  117.32    & 1.83    &  3829.64      & 1.38E-04    \\ 
        \multicolumn{1}{|l|}{}                     & 1728000  &   46.71     &  210.51    & 3.12    &  11575.3      & 1.67E-03    \\ \hline \hline
        \multicolumn{1}{|l|}{\multirow{4}{*}{\rotatebox[origin=c]{90}{\scriptsize $\mathcal{H}^2_{\sqrt{d}}$(b)}}}   & 64000  &  1.16    &  4.11   &  0.04  &  13.18  &  7.57E-05  \\ 
        \multicolumn{1}{|l|}{}                     & 512000   &  12.64      &  60.06     &  0.79   &  851.35      & 5.85E-04    \\  
        \multicolumn{1}{|l|}{}                     & 1000000  &  27.56      &  150.51    &  1.88   &  3829.64      & 1.24E-04    \\  
        \multicolumn{1}{|l|}{}                     & 1728000  &  50.79       &  267.6    &  3.61   &  11575.3      & 1.86E-03     \\ \hline \hline
        \multicolumn{1}{|l|}{\multirow{4}{*}{\rotatebox[origin=c]{90}{\scriptsize $\mathcal{H}_{*}$}}} & 64000  &  1.59   &   1.41   &  0.07  &  13.18   &   3.85E-05   \\  
        \multicolumn{1}{|l|}{}                     & 512000   &  20.07      &  45.67     &  0.88   &  851.35      & 4.89E-04    \\  
        \multicolumn{1}{|l|}{}                     & 1000000  &  50.21      &  109.16    &  4.18   &  3829.64      & 2.97E-04    \\ 
        \multicolumn{1}{|l|}{}                     & 1728000  &  86.94      &  192.64    &  7.78   &  11575.3      & 1.99E-03    \\ \hline \hline
        \multicolumn{1}{|l|}{\multirow{4}{*}{\rotatebox[origin=c]{90}{\scriptsize $\mathcal{H}_{\sqrt{d}}$}}}   & 64000  &  1.83   &  2.11   &   0.07  &  13.18    &  3.44E-05  \\  
        \multicolumn{1}{|l|}{}                     & 512000   &  23.77       &  47.23     &  1.06    &  851.35          &  2.98E-04   \\ 
        \multicolumn{1}{|l|}{}                     & 1000000  &  59.75       &  103.54    &  4.92    &   3829.64        &  2.43E-04   \\ 
        \multicolumn{1}{|l|}{}                     & 1728000  &  103.07      &  197.46    &   7.91   &   11575.3        &  5.69E-04   \\ \hline 
        \end{tabular}
}
    \caption{Comparison between the proposed algorithms and other \textbf{algebraic} algorithms as mentioned in \cref{tab:nex_notation} for the single layer Laplacian in $3$D with tolerance $\epsilon = 10^{-4}$. Rel. error refers to $2$-norm relative error in MVP, i.e., ${RE}_{MVP}$ (refer to \Cref{tab:app_notations}).}
    \label{table:smvp3d_lap_1}
\end{table} 
\begin{table}[H]
      \centering
      \resizebox{12cm}{!}{%
      \setlength\extrarowheight{0.9pt}
        \begin{tabular}{l|l|l|l|l|l|l|}
        \cline{2-7}
                                                   & \scriptsize $N$ & \scriptsize Mem. (GB) & \scriptsize $t_{init}$ (s) & \scriptsize $t_{MVP}$ (s) &\scriptsize $t_{direct}$ (s) & \scriptsize Rel. error $({RE}_{MVP})$ \\ \hline \cline{2-7}
        \multicolumn{1}{|l|}{\multirow{4}{*}{\rotatebox[origin=c]{90}{\scriptsize $\mathcal{H}^2_{*}$(b$+$t)}}}   & 64000  & 2.07   &  10.71   &  0.07  &  13.18   &  1.69E-06  \\  
        \multicolumn{1}{|l|}{}                     &  512000  &  22.98      &  237.27    &  0.67   &  851.35      &  4.41E-06   \\  
        \multicolumn{1}{|l|}{}                     &  1000000 &  43.38      &  528.45    &  1.34   &  3829.64    &  5.34E-06   \\  
        \multicolumn{1}{|l|}{}                     &  1728000 &  76.35      &   951.12   &  2.75   &  11575.3    &  1.30E-05   \\ \hline \hline
        \multicolumn{1}{|l|}{\multirow{4}{*}{\rotatebox[origin=c]{90}{\scriptsize $\bkt{\mathcal{H}^2 + \mathcal{H}}_{*}$}}} &  64000  & 2.27    & 8.73  & 0.10   &  13.18    &  5.91E-07   \\  
        \multicolumn{1}{|l|}{}                     & 512000   & 26.20       &  211.44    & 1.26    &  851.35     &  1.04E-06   \\  
        \multicolumn{1}{|l|}{}                     & 1000000  & 51.14       &  477.23    & 2.65    &  3829.64    &  1.64E-06   \\   
        \multicolumn{1}{|l|}{}                     & 1728000  & 91.67       &  863.76    & 5.21    &  11575.3    &  1.76E-05   \\ \hline \hline
        \multicolumn{1}{|l|}{\multirow{4}{*}{\rotatebox[origin=c]{90}{\scriptsize $\mathcal{H}^2_{\sqrt{d}}$(b)}}}   & 64000  & 2.47    & 12.27    &  0.11   &  13.18   &  2.02E-06  \\   
        \multicolumn{1}{|l|}{}                     & 512000   &  28.38      &  265.71    & 1.41    &  851.35     &  8.44E-06   \\  
        \multicolumn{1}{|l|}{}                     & 1000000  &  54.39      &  641.85    & 2.77    &  3829.64    &  1.41E-05   \\  
        \multicolumn{1}{|l|}{}                     & 1728000  &  96.25      &  1121.62   & 5.89    &  11575.3    &  3.91E-05   \\   \hline \hline
        \multicolumn{1}{|l|}{\multirow{4}{*}{\rotatebox[origin=c]{90}{\scriptsize $\mathcal{H}_{*}$}}} & 64000  &  2.78  & 5.78   & 0.24  & 13.18    &  2.95E-06  \\    
        \multicolumn{1}{|l|}{}                     & 512000   &  36.69       &  78.41     &  2.25    &  851.35     &  4.25E-06   \\  
        \multicolumn{1}{|l|}{}                     & 1000000  &  98.93       &  255.34    &  5.36    &  3829.64    &  8.33E-06   \\  
        \multicolumn{1}{|l|}{}                     & 1728000  &  172.26      &  589.91    &  11.35   &  11575.3    &  1.56E-05   \\   \hline \hline
        \multicolumn{1}{|l|}{\multirow{4}{*}{\rotatebox[origin=c]{90}{\scriptsize $\mathcal{H}_{\sqrt{d}}$}}}   & 64000  & 3.07    & 5.87   &  0.25   & 13.18   &  8.54E-07 \\    
        \multicolumn{1}{|l|}{}                     & 512000   & 41.91        &  89.41          &  2.87   &  851.35    & 1.27E-06    \\   
        \multicolumn{1}{|l|}{}                     & 1000000  & 113.34       &  257.31         &  5.88   &  3829.64    & 3.51E-06    \\  
        \multicolumn{1}{|l|}{}                     & 1728000  &  -           &   -             &  -   &   -   &  -   \\   \hline 
        \end{tabular}
}
    \caption{Comparison between the proposed algorithms and other \textbf{algebraic} algorithms as mentioned in \cref{tab:nex_notation} for the single layer Laplacian in $3$D with tolerance $\epsilon = 10^{-6}$. Rel. error refers to $2$-norm relative error in MVP, i.e., ${RE}_{MVP}$ (refer to \Cref{tab:app_notations}).}
    \label{table:smvp3d_lap_2}
\end{table}

\begin{table}[H]
      \centering
      \resizebox{12cm}{!}{%
        \begin{tabular}{l|l|l|l|l|l|l|}
        \cline{2-7}
                                                   & \scriptsize $N$ & \scriptsize Mem. (GB) & \scriptsize $t_{init}$ (s) & \scriptsize $t_{MVP}$ (s) &\scriptsize $t_{direct}$ (s) & \scriptsize Rel. error $({RE}_{MVP})$ \\ \hline \cline{2-7}
        \multicolumn{1}{|l|}{\multirow{4}{*}{\rotatebox[origin=c]{90}{\scriptsize $\mathcal{H}^2_{*}$(b$+$t)}}}   & 64000  &  2.94  & 21.16   &  0.10  &   13.18  & 7.96E-09  \\  
        \multicolumn{1}{|l|}{}                     &  512000  &  35.17      &  540.01     & 1.38    & 851.35      & 1.85E-08    \\  
        \multicolumn{1}{|l|}{}                     &  1000000 &  77.47      &  1260.17    & 2.84    & 3829.64     & 5.45E-08    \\  
        \multicolumn{1}{|l|}{}                     &  1728000 &  143.20      & 2502.12    & 5.11    & 11575.3     & 1.68E-07    \\ \hline \hline
        \multicolumn{1}{|l|}{\multirow{4}{*}{\rotatebox[origin=c]{90}{\scriptsize $\bkt{\mathcal{H}^2 + \mathcal{H}}_{*}$}}} &  64000  &  3.19  &  17.34   & 0.12    & 13.18   &  4.95E-08  \\  
        \multicolumn{1}{|l|}{}                     & 512000   &   39.52     &  498.17      &  1.43   & 851.35      &  5.12E-08   \\  
        \multicolumn{1}{|l|}{}                     & 1000000  &   88.22     &  1104.75     &  3.45   & 3829.64     &  1.03E-07   \\   
        \multicolumn{1}{|l|}{}                     & 1728000  &   165.07     &  2111.28    &  6.58   & 11575.3     &  1.33E-07    \\ \hline \hline
        \multicolumn{1}{|l|}{\multirow{4}{*}{\rotatebox[origin=c]{90}{\scriptsize $\mathcal{H}^2_{\sqrt{d}}$(b)}}}   & 64000  &  3.44   &  30.07    &  0.12   &  13.18  &  4.79E-08  \\   
        \multicolumn{1}{|l|}{}                     & 512000   &  42.99      &  645.32      &  1.68   &  851.35    &  8.53E-08   \\  
        \multicolumn{1}{|l|}{}                     & 1000000  &  96.54      &  1672.41     &  4.79   &  3829.64    &  1.27E-07   \\  
        \multicolumn{1}{|l|}{}                     & 1728000  &    -        &   -   &   -  &    -  &  -   \\   \hline \hline
        \multicolumn{1}{|l|}{\multirow{4}{*}{\rotatebox[origin=c]{90}{\scriptsize $\mathcal{H}_{*}$}}} & 64000  &  4.32    & 10.12    &  0.23   & 13.18    & 1.63E-08   \\    
        \multicolumn{1}{|l|}{}                     & 512000   &   58.89      &  195.06    &  3.42   &  851.35    &  2.60E-08   \\  
        \multicolumn{1}{|l|}{}                     & 1000000  &   153.92     &  565.13    &  8.78   &  3829.64    &  2.50E-07   \\  
        \multicolumn{1}{|l|}{}                     & 1728000  &   -          &   -   &  -   &   -   &  -   \\   \hline \hline
        \multicolumn{1}{|l|}{\multirow{4}{*}{\rotatebox[origin=c]{90}{\scriptsize $\mathcal{H}_{\sqrt{d}}$}}}   & 64000  &  4.71   &  9.15   &  0.25   & 13.18    &  8.47E-09  \\    
        \multicolumn{1}{|l|}{}                     & 512000   &   66.16     &  172.22    &  3.85      &  851.35    &  1.20E-08   \\   
        \multicolumn{1}{|l|}{}                     & 1000000  &   180.43     &  511.13    & 10.55     &  3829.64    & 4.15E-08    \\  
        \multicolumn{1}{|l|}{}                     & 1728000  &   -     &   -   &  -   &   -   &  -   \\   \hline 
        \end{tabular}
}
    \caption{Comparison between the proposed algorithms and other \textbf{algebraic} algorithms as mentioned in \cref{tab:nex_notation} for the single layer Laplacian in $3$D with tolerance $\epsilon = 10^{-8}$. Rel. error refers to $2$-norm relative error in MVP, i.e., ${RE}_{MVP}$ (refer to \Cref{tab:app_notations}).}
    \label{table:smvp3d_lap_3}
\end{table}

From the \Cref{table:smvp3d_lap_1}, \Cref{table:smvp3d_lap_2}, and \Cref{table:smvp3d_lap_3}, the following conclusions can be drawn.

\begin{enumerate}
    \item The scaling of the memory and MVP time of the $\mathcal{H}^2_{*}$(b$+$t) algorithm is similar to that of $\mathcal{H}^2_{\sqrt{d}}$(b), and notably, $\mathcal{H}^2_{*}$(b$+$t) shows less storage and improved MVP time. \Cref{table:smvp3d_lap_3} indicates that we get output for $N=1728000$ in the proposed algorithms ($\mathcal{H}^2_{*}$(b$+$t), $\bkt{\mathcal{H}^2 + \mathcal{H}}_{*}$). However, outputs corresponding to $N = 1728000$ are unavailable for the other algorithms as they ran out of memory.
    \item In $3$D, the initialization time of $\mathcal{H}^2_{*}$(b$+$t) is better than that of $\mathcal{H}^2_{\sqrt{d}}$(b) for a reasonably large value of $N$ (up to $N \leq 1728000$, the maximum value of $N$ achievable in our system), even though the initialization of $\mathcal{H}^2_{\sqrt{d}}$(b) scales linearly. The reason is as we go for the higher dimensional problems, the interaction list size as well as the far-field rank value, grow exponentially with the underlying dimension (the far-field rank scales as $\mathcal{O}(1)$ with $N$, but the value of the constant is large), and the $\mathcal{H}^2_{\sqrt{d}}$(b) has a larger interaction list compared to the $\mathcal{H}^2_{*}$(b$+$t) (refer to \Cref{tab:complexity_com}).
    \item $\mathcal{H}^2_{*}$(b$+$t) and $\mathcal{H}^2_{\sqrt{d}}$(b) are the nested versions of $\mathcal{H}_{*}$ and $\mathcal{H}_{\sqrt{d}}$ algorithms, respectively. We can see that the nested versions require much less storage and improved MVP time than the non-nested ones.
    \item It is noteworthy that in $3$D, the \emph{semi-nested} $\bkt{\mathcal{H}^2 + \mathcal{H}}_{*}$ is slightly better than the $\mathcal{H}^2_{\sqrt{d}}$(b) in terms of both memory and time.
\end{enumerate}

We arrange the above algorithms in ascending order with respect to memory and MVP time.\\
\textbf{Memory:} $\mathcal{H}^2_{*}$(b$+$t) < $\bkt{\mathcal{H}^2 + \mathcal{H}}_{*}$ < $\mathcal{H}^2_{\sqrt{d}}$(b) < $\mathcal{H}_{*}$ < $\mathcal{H}_{\sqrt{d}}$ \\ 
\textbf{MVP time:} $\mathcal{H}^2_{*}$(b$+$t) < $\bkt{\mathcal{H}^2 + \mathcal{H}}_{*}$ < $\mathcal{H}^2_{\sqrt{d}}$(b) < $\mathcal{H}_{*}$ < $\mathcal{H}_{\sqrt{d}}$

Therefore, though we do not have a strict theoretical error bound, the above numerical experiment demonstrates that in $3$D, the proposed $\mathcal{H}^2_{*}$(b$+$t) and $\bkt{\mathcal{H}^2 + \mathcal{H}}_{*}$ algorithms show less memory usage and improved MVP time than the standard $\mathcal{H}^2$ matrix algorithm with comparable accuracy.

\textbf{We also perform the fast MVPs for the Matérn and Helmholtz kernels with the performance analysis of both the proposed and various nested algorithms as mentioned in \Cref{tab:nex_notation}. Please refer to \Cref{other_kernels} for more details.}

\subsubsection{Fast MVP accelerated GMRES for integral equation in \texorpdfstring{$3$}{3}D}
\begin{figure}[H]
    \centering
    \subfloat[]{\includegraphics[height=3.8cm, width=5.5cm]{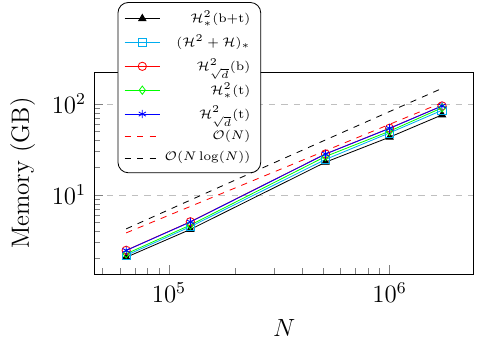}\label{green3d_mem}}%
    \subfloat[]{\includegraphics[height=3.8cm, width=5.5cm]{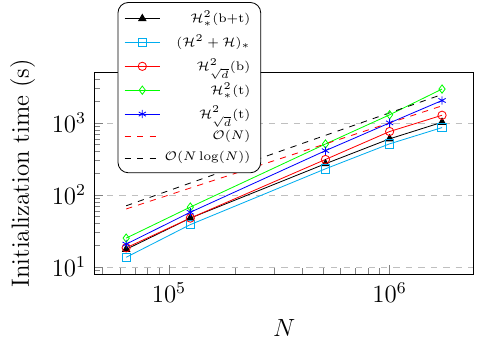}\label{green3d_assem}}%
    \subfloat[]{\includegraphics[height=3.8cm, width=5.5cm]{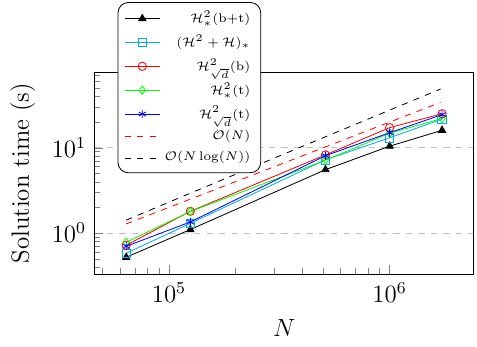}\label{green3d_sol}}%
    \caption{Plots of Memory, Initialization time and Solution time of the fast MVP algorithms.}
    \label{fig:green3d}
\end{figure}

We consider the Fredholm integral equation of the second kind over $C = [-1,1]^3 \subset \mathbb{R}^3$, which is given by
\begin{equation} \label{ieq3d_1}
    \sigma(\pmb{x})+\int_{C} F(\pmb{x},\pmb{y})\sigma(\pmb{y})d\pmb{y} = f(\pmb{x}) \qquad \pmb{x}, \pmb{y} \in C
\end{equation}
with $F(\pmb{x},\pmb{y}) = \dfrac{1}{4\pi \magn{\pmb{x}-\pmb{y}}_2}$ (Green's function for Laplace equation in $3$D). We follow a piece-wise constant collocation method with collocation points on a uniform grid in $C = [-1,1]^3$ to discretize \cref{ieq3d_1} and set the RHS $\pmb{f}$ as described in \Cref{tab:app_notations}. Therefore, we obtain the following linear system 
\begin{equation} \label{ieq3d_2}
    K \pmb{\sigma} = \pmb{f}
\end{equation}
The \Cref{ieq3d_2} is solved using \textbf{fast} GMRES, i.e., the MVP part of GMRES is accelerated using the $\mathcal{H}^2_{*}$(b$+$t), $\bkt{\mathcal{H}^2 + \mathcal{H}}_{*}$, $\mathcal{H}^2_{\sqrt{d}}$(b), $\mathcal{H}^2_{*}$(t) and $\mathcal{H}^2_{\sqrt{d}}$(t) algorithms. We set $\epsilon = 10^{-6}$ and $\epsilon_{GMRES} = 10^{-10}$ and report the memory, initialization time, and solution time in \Cref{table:green3d_int_eq}. The relative error in solution $(RE_{sol})$ is of order $10^{-10}$ in all cases. We also plot the memory (\Cref{green3d_mem}), initialization time (\Cref{green3d_assem}) and solution time (\Cref{green3d_sol}) in \Cref{fig:green3d}.
\begin{table}[H]
      \centering
      \resizebox{\textwidth}{!}{%
      \setlength\extrarowheight{0.9pt}
\begin{tabular}{|l|lllll|lllll|lllll|l|}
\hline
\multirow{2}{*}{N} & \multicolumn{5}{l|}{\qquad \qquad \qquad Memory (GB)}                                                                                           & \multicolumn{5}{l|}{\qquad \qquad \quad Initialization time (s)}                                                                                & \multicolumn{5}{l|}{\qquad \qquad \quad Solution time (s)}                                                                                      & \multirow{2}{*}{$\#iter$} \\ \cline{2-16}
                   & \multicolumn{1}{l|}{\scriptsize $\mathcal{H}^2_{*}$(b$+$t)} & \multicolumn{1}{l|}{\scriptsize $\bkt{\mathcal{H}^2 + \mathcal{H}}_{*}$} & \multicolumn{1}{l|}{\scriptsize $\mathcal{H}^2_{\sqrt{d}}$(b)} & \multicolumn{1}{l|}{\scriptsize $\mathcal{H}^2_{*}$(t)} & \scriptsize $\mathcal{H}^2_{\sqrt{d}}$(t) & \multicolumn{1}{l|}{\scriptsize $\mathcal{H}^2_{*}$(b$+$t)} & \multicolumn{1}{l|}{\scriptsize $\bkt{\mathcal{H}^2 + \mathcal{H}}_{*}$} & \multicolumn{1}{l|}{\scriptsize $\mathcal{H}^2_{\sqrt{d}}$(b)} & \multicolumn{1}{l|}{\scriptsize $\mathcal{H}^2_{*}$(t)} & \scriptsize $\mathcal{H}^2_{\sqrt{d}}$(t) & \multicolumn{1}{l|}{\scriptsize $\mathcal{H}^2_{*}$(b$+$t)} & \multicolumn{1}{l|}{\scriptsize $\bkt{\mathcal{H}^2 + \mathcal{H}}_{*}$} & \multicolumn{1}{l|}{\scriptsize $\mathcal{H}^2_{\sqrt{d}}$(b)} & \multicolumn{1}{l|}{\scriptsize $\mathcal{H}^2_{*}$(t)} & \scriptsize $\mathcal{H}^2_{\sqrt{d}}$(t) &                       \\ \hline 
                   64000  & \multicolumn{1}{l|}{2.07}     & \multicolumn{1}{l|}{2.17}       & \multicolumn{1}{l|}{2.47}      & \multicolumn{1}{l|}{2.26}       &   2.46    & \multicolumn{1}{l|}{17.6}     & \multicolumn{1}{l|}{13.66}       & \multicolumn{1}{l|}{18.5}      & \multicolumn{1}{l|}{25.22}        &   20.81    & \multicolumn{1}{l|}{0.52}     & \multicolumn{1}{l|}{0.55}       & \multicolumn{1}{l|}{0.72}      & \multicolumn{1}{l|}{0.78}        &   0.7    &   5                    \\ \hline
                  125000 & \multicolumn{1}{l|}{4.19}     & \multicolumn{1}{l|}{4.49}       & \multicolumn{1}{l|}{5.09}      & \multicolumn{1}{l|}{4.68}        &   5.06    & \multicolumn{1}{l|}{47.81}     & \multicolumn{1}{l|}{38.84}       & \multicolumn{1}{l|}{52.27}      & \multicolumn{1}{l|}{67.97}        &  57.6     & \multicolumn{1}{l|}{1.10}     & \multicolumn{1}{l|}{1.16}       & \multicolumn{1}{l|}{1.80}      & \multicolumn{1}{l|}{1.78}        &  1.35     &    4                   \\ \hline
                  512000 & \multicolumn{1}{l|}{22.99}     & \multicolumn{1}{l|}{24.64}       & \multicolumn{1}{l|}{28.33}      & \multicolumn{1}{l|}{26.54}    &   28.25    & \multicolumn{1}{l|}{271.53}     & \multicolumn{1}{l|}{230.77}       & \multicolumn{1}{l|}{311.19}      & \multicolumn{1}{l|}{510.19}        &  412.39    & \multicolumn{1}{l|}{5.54}     & \multicolumn{1}{l|}{7.23}       & \multicolumn{1}{l|}{8.25}      & \multicolumn{1}{l|}{7.15}        &  8.07     &   5                    \\ \hline
                  1000000 & \multicolumn{1}{l|}{43.34}     & \multicolumn{1}{l|}{48.30}       & \multicolumn{1}{l|}{54.38}      & \multicolumn{1}{l|}{50.40}   &   54.11    & \multicolumn{1}{l|}{601.7}     & \multicolumn{1}{l|}{512.91}       & \multicolumn{1}{l|}{760.36}      & \multicolumn{1}{l|}{1296.93}        &  1001.32   & \multicolumn{1}{l|}{10.44}     & \multicolumn{1}{l|}{13.09}       & \multicolumn{1}{l|}{17.22}      & \multicolumn{1}{l|}{14.81}        &   15.02    &  5                     \\ \hline
                  1728000 & \multicolumn{1}{l|}{76.12}     & \multicolumn{1}{l|}{85.52}       & \multicolumn{1}{l|}{95.56}      & \multicolumn{1}{l|}{91.11}    &  95.49     & \multicolumn{1}{l|}{1021.61}     & \multicolumn{1}{l|}{861.95}       & \multicolumn{1}{l|}{1278.25}      & \multicolumn{1}{l|}{2952.32}        & 2049.52 & \multicolumn{1}{l|}{16.03}     & \multicolumn{1}{l|}{21.68}       & \multicolumn{1}{l|}{25.06}      & \multicolumn{1}{l|}{22.08}        &  24.47     & 4                      \\ \hline
\end{tabular}
}
    \caption{Performance of the proposed and the nested algorithms as mentioned in \Cref{tab:nex_notation}. We set $\epsilon = 10^{-6}$ and $\epsilon_{GMRES} = 10^{-10}$.}
    \label{table:green3d_int_eq}
\end{table}

We arrange the above algorithms in ascending order with respect to memory and solution time.\\
\textbf{Memory:} $\mathcal{H}^2_{*}$(b$+$t) < $\bkt{\mathcal{H}^2 + \mathcal{H}}_{*}$ < $\mathcal{H}^2_{*}$(t) < $\mathcal{H}^2_{\sqrt{d}}$(b) $\leq$ $\mathcal{H}^2_{\sqrt{d}}$(t) \\ 
\textbf{Solution time:} $\mathcal{H}^2_{*}$(b$+$t) < $\bkt{\mathcal{H}^2 + \mathcal{H}}_{*}$ < $\mathcal{H}^2_{*}$(t) < $\mathcal{H}^2_{\sqrt{d}}$(t) < $\mathcal{H}^2_{\sqrt{d}}$(b)

\subsubsection{Fast MVP accelerated GMRES for RBF interpolation in \texorpdfstring{$3$}{3}D}
 \RK{Let the location of the particles $\{\pmb{x}_i \}_{i=1}^N$ be the $N^{1/3} \times N^{1/3} \times N^{1/3}$ Chebyshev grid on the domain $[-1,1]^3$. We choose the Chebyshev distribution of particles to study the performance of the nested algorithms over non-uniformly distributed particles. However, we use the uniform $2^d$ tree (oct tree in $3$D) as described in \Cref{tree_construction}. In this experiment, we consider the following non-translation invariant (NTI) scaled multiquadric radial basis function \cite{bozzini2015interpolation} 
\begin{equation}
  F \bkt{\pmb{x}, \pmb{y}} =  \bkt{\magn{\pmb{x} - \pmb{y}}_2^2 + \bkt{c(\pmb{x}) - c(\pmb{y})}^2 + 1}^{3/2}
\end{equation}
where $c(\pmb{x}) = \exp \bkt{-(\pmb{x}(1) + \pmb{x}(3))}$ and $\pmb{x} = \bkt{\pmb{x}(1), \pmb{x}(2), \pmb{x}(3)}$.

Consider the following linear system generated by the above NTI kernel function $F$
\begin{equation} \label{eq:nti_rbf_fmm_1}
    \lambda_i \sqrt{N} + \dsum_{j=1}^N F \bkt{\pmb{x}_i, \pmb{x}_j} \lambda_j = b_i, \qquad i = 1,2,\dots, N.
\end{equation}

By setting $\pmb{b}$ as described in the $11^{th}$ cell of \Cref{tab:app_notations}, the \Cref{eq:nti_rbf_fmm_1} can be written in the form 
\begin{equation} \label{eq:nti_rbf_system_1}
    K \pmb{\lambda} = \pmb{b}
\end{equation}
The \Cref{eq:nti_rbf_system_1} is solved using \textbf{fast} GMRES. We set $\epsilon = 10^{-6}$ and $\epsilon_{GMRES} = 10^{-10}$ and tabulate the memory, initialization time and solution time in \Cref{table:nti_rbf3d_alg} for all the nested algorithms. The relative error in solution $(RE_{sol})$ is of order $10^{-6}$ in all cases. We also plot the memory (\Cref{nti_alg_mem}), initialization time (\Cref{nti_alg_assem}) and solution time (\Cref{nti_alg_sol}) in \Cref{fig:nti_alg}.

\begin{table}[H]
      \centering
      \resizebox{\textwidth}{!}{%
      \setlength\extrarowheight{0.9pt}
\begin{tabular}{|l|lllll|lllll|lllll|l|}
\hline
\multirow{2}{*}{N} & \multicolumn{5}{l|}{\qquad \qquad \qquad Memory (GB)}                                                                                           & \multicolumn{5}{l|}{\qquad \qquad \quad Initialization time (s)}                                                                                & \multicolumn{5}{l|}{\qquad \qquad \quad Solution time (s)}                                                                                      & \multirow{2}{*}{$\#iter$} \\ \cline{2-16}
                   & \multicolumn{1}{l|}{\scriptsize $\mathcal{H}^2_{*}$(b$+$t)} & \multicolumn{1}{l|}{\scriptsize $\bkt{\mathcal{H}^2 + \mathcal{H}}_{*}$} & \multicolumn{1}{l|}{\scriptsize $\mathcal{H}^2_{\sqrt{d}}$(b)} & \multicolumn{1}{l|}{\scriptsize $\mathcal{H}^2_{*}$(t)} & \scriptsize $\mathcal{H}^2_{\sqrt{d}}$(t) & \multicolumn{1}{l|}{\scriptsize $\mathcal{H}^2_{*}$(b$+$t)} & \multicolumn{1}{l|}{\scriptsize $\bkt{\mathcal{H}^2 + \mathcal{H}}_{*}$} & \multicolumn{1}{l|}{\scriptsize $\mathcal{H}^2_{\sqrt{d}}$(b)} & \multicolumn{1}{l|}{\scriptsize $\mathcal{H}^2_{*}$(t)} & \scriptsize $\mathcal{H}^2_{\sqrt{d}}$(t) & \multicolumn{1}{l|}{\scriptsize $\mathcal{H}^2_{*}$(b$+$t)} & \multicolumn{1}{l|}{\scriptsize $\bkt{\mathcal{H}^2 + \mathcal{H}}_{*}$} & \multicolumn{1}{l|}{\scriptsize $\mathcal{H}^2_{\sqrt{d}}$(b)} & \multicolumn{1}{l|}{\scriptsize $\mathcal{H}^2_{*}$(t)} & \scriptsize $\mathcal{H}^2_{\sqrt{d}}$(t) &                       \\ \hline 
            64000      & \multicolumn{1}{l|}{0.69}     & \multicolumn{1}{l|}{0.69}       & \multicolumn{1}{l|}{0.79}      & \multicolumn{1}{l|}{0.76}        &    0.79   & \multicolumn{1}{l|}{26.25}     & \multicolumn{1}{l|}{25.12}       & \multicolumn{1}{l|}{30.32}      & \multicolumn{1}{l|}{27.86}        &  39.97     & \multicolumn{1}{l|}{1.03}     & \multicolumn{1}{l|}{1.11}       & \multicolumn{1}{l|}{1.57}      & \multicolumn{1}{l|}{1.14}        &   1.56      &  34            \\ \hline
            125000    & \multicolumn{1}{l|}{0.92}     & \multicolumn{1}{l|}{0.92}       & \multicolumn{1}{l|}{1.08}      & \multicolumn{1}{l|}{1.07}        &    1.12   & \multicolumn{1}{l|}{52.21}     & \multicolumn{1}{l|}{49.79}       & \multicolumn{1}{l|}{62.04}      & \multicolumn{1}{l|}{68.54}        & 134.92      & \multicolumn{1}{l|}{2.96}     & \multicolumn{1}{l|}{3.18}       & \multicolumn{1}{l|}{4.31}      & \multicolumn{1}{l|}{3.26}        &   3.88    &  37                     \\ \hline
           512000    & \multicolumn{1}{l|}{7.53}     & \multicolumn{1}{l|}{7.65}       & \multicolumn{1}{l|}{8.73}      & \multicolumn{1}{l|}{8.39}        &  9.02   & \multicolumn{1}{l|}{243.02}     & \multicolumn{1}{l|}{205.45}       & \multicolumn{1}{l|}{277.25}      & \multicolumn{1}{l|}{417.53}     & 562.48      & \multicolumn{1}{l|}{12.11}     & \multicolumn{1}{l|}{17.06}       & \multicolumn{1}{l|}{20.98}      & \multicolumn{1}{l|}{14.49}        &  19.02     & 42                      \\ \hline
         1000000    & \multicolumn{1}{l|}{7.79}     & \multicolumn{1}{l|}{8.69}       & \multicolumn{1}{l|}{9.81}      & \multicolumn{1}{l|}{9.59}        &   11.01  & \multicolumn{1}{l|}{326.58}     & \multicolumn{1}{l|}{288.06}       & \multicolumn{1}{l|}{382.72}      & \multicolumn{1}{l|}{828.57}        & 1070.25      & \multicolumn{1}{l|}{32.46}     & \multicolumn{1}{l|}{40.50}       & \multicolumn{1}{l|}{46.58}      & \multicolumn{1}{l|}{43.48}        &  45.32     &  46                     \\ \hline
        1728000    & \multicolumn{1}{l|}{17.51}     & \multicolumn{1}{l|}{19.56}       & \multicolumn{1}{l|}{21.15}      & \multicolumn{1}{l|}{20.88}        &  23.66    & \multicolumn{1}{l|}{552.26}     & \multicolumn{1}{l|}{492.06}       & \multicolumn{1}{l|}{656.43}      & \multicolumn{1}{l|}{1754.62}        & 1879.72      & \multicolumn{1}{l|}{57.48}     & \multicolumn{1}{l|}{68.81}       & \multicolumn{1}{l|}{76.97}      & \multicolumn{1}{l|}{72.96}        & 74.46      &  52                       \\ \hline
\end{tabular}
}
    \caption{Performance of the proposed and the nested algorithms as mentioned in \Cref{tab:nex_notation}. We set $\epsilon = 10^{-6}$ and $\epsilon_{GMRES} = 10^{-10}$.}
    \label{table:nti_rbf3d_alg}
\end{table}

 \begin{figure}[H]
    \centering
    \subfloat[]{\includegraphics[height=3.8cm, width=5.5cm]{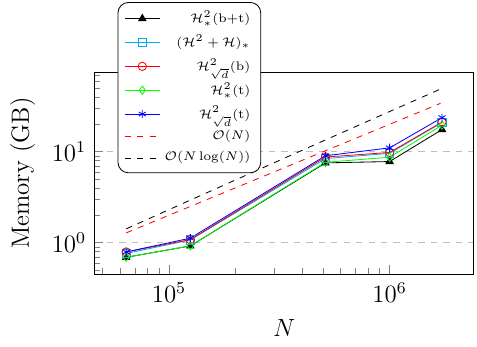}\label{nti_alg_mem}}%
    \subfloat[]{\includegraphics[height=3.8cm, width=5.5cm]{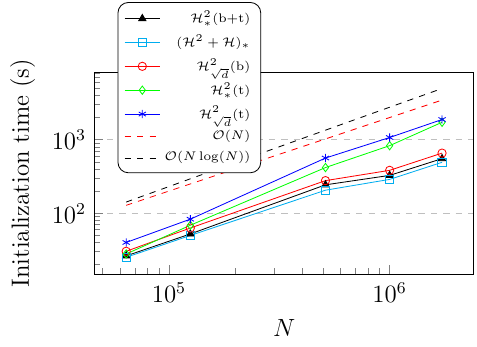}\label{nti_alg_assem}}%
    \subfloat[]{\includegraphics[height=3.8cm, width=5.5cm]{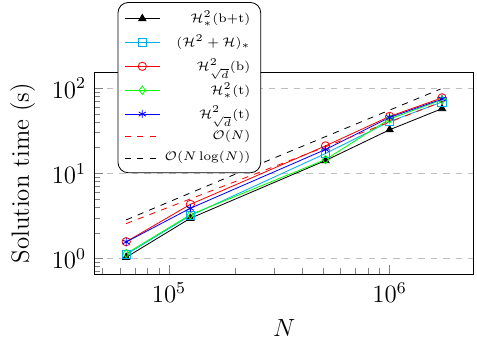}\label{nti_alg_sol}}%
    \caption{Plots of Memory, Initialization time and Solution time of the fast MVP algorithms.}
    \label{fig:nti_alg}
\end{figure}

The kinks in the plots are due to the additional levels in the oct tree structure. We arrange the above algorithms in ascending order with respect to memory and solution time.\\
\textbf{Memory:} $\mathcal{H}^2_{*}$(b$+$t) < $\bkt{\mathcal{H}^2 + \mathcal{H}}_{*}$ < $\mathcal{H}^2_{*}$(t)  < $\mathcal{H}^2_{\sqrt{d}}$(b) < $\mathcal{H}^2_{\sqrt{d}}$(t) \\ 
\textbf{Solution time:} $\mathcal{H}^2_{*}$(b$+$t) < $\bkt{\mathcal{H}^2 + \mathcal{H}}_{*}$ < $\mathcal{H}^2_{*}$(t) < $\mathcal{H}^2_{\sqrt{d}}$(t) < $\mathcal{H}^2_{\sqrt{d}}$(b)
}

\emph{\textbf{Summary of the experiments in $3$D.}} Here is a summary of the results obtained from all the experiments performed in $3$D.
\begin{enumerate}
    \item In all the experiments, the proposed $\mathcal{H}^2_{*}$(b$+$t) outperforms all the other algebraic algorithms (\Cref{tab:nex_notation}) in terms of memory and MVP/Solution time. Also, the $\mathcal{H}^2_{*}$(b$+$t) has the least initialization time among the fully nested algorithms. Hence, $\mathcal{H}^2_{*}$(b$+$t) could be an attractive alternative to the standard $\mathcal{H}^2$ matrix algorithm.
    \item The initialization time of the $\mathcal{H}^2_{*}$(t) algorithm (\Cref{nca_t2b_entire_hodlrdd}) is the highest and the proposed $\mathcal{H}^2_{*}$(b$+$t) performs way better than the $\mathcal{H}^2_{*}$(t) algorithm. 
    \item Note that in $3$D, the $\bkt{\mathcal{H}^2 + \mathcal{H}}_{*}$, $\mathcal{H}^2_{\sqrt{d}}$(b) and $\mathcal{H}^2_{\sqrt{d}}$(t) exhibit similar scaling in memory. The numerical results demonstrate that the $\bkt{\mathcal{H}^2 + \mathcal{H}}_{*}$ beats the $\mathcal{H}^2_{\sqrt{d}}$(b) and $\mathcal{H}^2_{\sqrt{d}}$(t) in terms of memory, initialization time and MVP/Solution time. 
\end{enumerate}

Therefore, in both $2$D and $3$D, the proposed $\mathcal{H}^2_{*}$(b$+$t) algorithm outperforms the NCA-based standard $\mathcal{H}^2$ matrix algorithms \cite{zhao2019fast} in the context of memory and MVP/Solution time. Also, in $3$D, the proposed $\bkt{\mathcal{H}^2 + \mathcal{H}}_{*}$ algorithm performs slightly better than the standard $\mathcal{H}^2$ matrix algorithms.
\section{Conclusion} \label{conclusion}
\rk{We have presented an efficient \emph{\textbf{nested}} hierarchical matrix algorithm, $\mathcal{H}^2_{*}$(b$+$t), and a \emph{\textbf{semi-nested}} hierarchical matrix algorithm, $\bkt{\mathcal{H}^2 + \mathcal{H}}_{*}$. In our proposed hierarchical matrix algorithms, the admissible clusters are the far-field and the vertex-sharing clusters}. Due to the use of nested form of the bases, the computational cost is significantly reduced. We evaluate the performance of the proposed algorithms through extensive numerical experiments conducted in both $2$D and $3$D, which include time taken for various kernel MVPs and solving linear systems using the fast GMRES. We also compare the proposed algorithms with different related algebraic fast MVP algorithms and present various benchmarks. Notably, all these algorithms are developed purely algebraically, making them kernel-independent. To the best of our knowledge, this is the first work to study the performance analysis of a wide range of algebraic fast MVP algorithms in $2$D and $3$D. Numerical results show that the proposed algorithms are competitive with the standard $\mathcal{H}^2$ algorithm with respect to the memory and MVP time. Therefore, they could be an attractive alternative to the standard $\mathcal{H}^2$ matrix algorithm. Finally, we would like to release the implementation of the proposed algorithms made available at \texttt{\url{https://github.com/riteshkhan/H2weak/}}.
\section*{Acknowledgments}
We acknowledge the use of the computing resources at HPCE, IIT Madras.
\bibliographystyle{siamplain}
\bibliography{refs}
\appendix
\section{Purely algebraic \texorpdfstring{$\mathcal{H}^2_{\sqrt{d}}$}{H2}(b) algorithm \cite{zhao2019fast, gujjula2022new}} \label{b2t_h2matrix} 
In this section, we discuss a purely algebraic standard $\mathcal{H}^2$ MVP algorithm. Since the construction of the standard $\mathcal{H}^2$ matrix representation is based on the B$2$T NCA, we denote this algorithm as $\mathcal{H}^2_{\sqrt{d}}$(b). We discuss $\mathcal{H}^2_{\sqrt{d}}$(b) MVP algorithm in two steps:
\begin{enumerate}
    \item Initialization of $\mathcal{H}^2_{\sqrt{d}}$(b) representation. (\Cref{b2t_h2matrix_init})
    \item Calculation of the potential (MVP). (\Cref{b2t_h2matrix_mvp})
\end{enumerate}

\subsection{Initialization of \texorpdfstring{$\mathcal{H}^2_{\sqrt{d}}$(b)}{H2} representation} \label{b2t_h2matrix_init}
In the B$2$T pivot selection, the pivots of a cluster at a parent level are obtained from the pivots at its child level. Therefore, one needs to traverse the $2^d$ uniform tree from bottom to top direction (starting at the leaf level) to find pivots of all clusters of the tree. The detailed procedure to obtain the four sets of pivots corresponding to a cluster $X$ is given below. 

\begin{itemize}
        \item If $X$ is a leaf cluster (childless), then construct the following four sets
    \end{itemize}
    \begin{align}
        \Tilde{t}^{X,i} := t^X \qquad \qquad \text{ and } \qquad \qquad \Tilde{s}^{X,i} := \bigcup_{Y \in \mathcal{IL}_{\sqrt{d}}(X)} s^Y
    \end{align}
     \begin{align}
        \Tilde{t}^{X,o} := \bigcup_{Y \in \mathcal{IL}_{\sqrt{d}}(X)} t^Y \qquad \qquad \text{ and } \qquad \qquad \Tilde{s}^{X,o} := s^X 
    \end{align}
    \begin{itemize}
        \item If $X$ is a non-leaf cluster, then construct the following four sets
    \end{itemize}
    \begin{align}
        \Tilde{t}^{X,i} := \bigcup_{X_c \in child(X)} t^{X_c, i} \qquad \qquad \text{ and } \qquad \qquad \Tilde{s}^{X,i} := \bigcup_{Y \in \mathcal{IL}_{\sqrt{d}}(X)} \bigcup_{Y_c \in child(Y)} s^{Y_c,o}
    \end{align}
     \begin{align}
        \Tilde{t}^{X,o} := \bigcup_{Y \in \mathcal{IL}_{\sqrt{d}}(X)} \bigcup_{Y_c \in child(Y)} t^{Y_c,i} \qquad \qquad \text{ and } \qquad \qquad \Tilde{s}^{X,o} := \bigcup_{X_c \in child(X)} s^{X_c,o} 
    \end{align}
    
To obtain the pivots $t^{X,i}$, $s^{X,i}$, $t^{X,o}$ and $s^{X,o}$, one needs to perform ACA. We perform ACA on the matrix $K_{\Tilde{t}^{X,i}, \Tilde{s}^{X,i}}$ with user-given tolerance $\epsilon$. The sets $t^{X,i}$ and $s^{X,i}$ are the row and column pivots chosen by the ACA. The search spaces of the pivots for a particular cluster are illustrated in \Cref{fig:B2T_pivot}. Similarly, perform ACA on the matrix $K_{\Tilde{t}^{X,o}, \Tilde{s}^{X,o}}$ to obtain the other two sets of pivots $t^{X,o}$ and $s^{X,o}$. 

\begin{figure}[H]
\begin{center}
\captionsetup[subfloat]{labelformat=empty}
    \subfloat[]{
    \begin{tikzpicture}[scale=0.5]
        \draw[draw=black,fill=red] (-8,1) rectangle (-7,2);
        \pgfmathsetseed{58678246}
        \foreach \x in {1,...,200}
        {
		\pgfmathsetmacro{\x}{rand/2 -7.5}
		\pgfmathsetmacro{\y}{rand/2 +1.5}
		\fill[black]    (\x,\y) circle (0.015);
        };
        \pgfmathsetseed{58678246}
        \foreach \x in {1,...,80}
        {
		\pgfmathsetmacro{\x}{rand/2 -7.5}
		\pgfmathsetmacro{\y}{rand/2 +1.5}
		\fill[blue]    (\x,\y) circle (0.035);
        };

        \pgfmathsetseed{58678246}
        \foreach \x in {1,...,70}
        {
		\pgfmathsetmacro{\x}{rand/5 +0.25}
		\pgfmathsetmacro{\y}{rand/5 +0.25}
		\fill[black]    (\x,\y) circle (0.015);
        };
        \pgfmathsetseed{58678246}
        \foreach \x in {1,...,8}
        {
		\pgfmathsetmacro{\x}{rand/5 +0.25}
		\pgfmathsetmacro{\y}{rand/5 +0.25}
		\fill[blue]    (\x,\y) circle (0.035);
        };
    
        \fill [white] (0,0) rectangle (1,1);
        \fill [white] (0,1) rectangle (1,2);
        \fill [white] (0,2) rectangle (1,3);
        \fill [white] (0,3) rectangle (1,4);
        \fill [white] (1,3) rectangle (2,4);
        \fill [white] (2,3) rectangle (3,4);
        \fill [white] (3,3) rectangle (4,4);
        \node[anchor=north] at (1.5,2.8) {};
        \node[anchor=north] at (1.5,1.8) {};
        \node[anchor=north] at (.5,1.8) {};
        \draw[black] (0, 0) grid (4, 4);
        \draw[line width=0.5mm,  fill=green] (1, 0) rectangle (2, 1);
        \draw[line width=0.5mm,  fill=green] (1, 1) rectangle (2, 2);
        \draw[line width=0.5mm,  fill=green] (2, 2) rectangle (3, 3);
        \draw[line width=0.5mm,  fill=green] (3, 2) rectangle (4, 3);
        \draw[line width=0.5mm,  fill=green] (2, 0) rectangle (3, 1);
        \draw[line width=0.5mm,  fill=green] (3, 0) rectangle (4, 1);
        \draw[line width=0.5mm,  fill=green] (3, 1) rectangle (4, 2);
        \draw[line width=0.5mm,  fill=green] (1, 2) rectangle (2, 3);
        \fill [white] (1.5,1) rectangle (3,2.5);
        \fill [red] (2,1.5) rectangle (2.5,2);
        \node[anchor=north] at (2.27,2.15) {\tiny $X$};
        \draw[step=0.5cm, black] (1, 0) grid (4, 3);
        \draw[line width=0.5mm] (1, 0) rectangle (2, 1);
        \draw[line width=0.5mm] (1, 1) rectangle (2, 2);
        \draw[line width=0.5mm] (2, 2) rectangle (3, 3);
        \draw[line width=0.5mm] (3, 2) rectangle (4, 3);
        \draw[line width=0.5mm] (2, 0) rectangle (3, 1);
        \draw[line width=0.5mm] (3, 0) rectangle (4, 1);
        \draw[line width=0.5mm] (3, 1) rectangle (4, 2);
        \draw[line width=0.5mm] (1, 2) rectangle (2, 3);

        \pgfmathsetseed{58678246}
        \foreach \x in {1,...,400}
        {
		\pgfmathsetmacro{\x}{1.5*rand+2.5}
		\pgfmathsetmacro{\y}{rand/2+0.5}
		\fill[black]    (\x,\y) circle (0.015);
        };

        \pgfmathsetseed{52213446}
        \foreach \x in {1,...,200}
        {
		\pgfmathsetmacro{\x}{1.5*rand+2.5}
		\pgfmathsetmacro{\y}{rand/4+2.75}
		\fill[black]    (\x,\y) circle (0.015);
        };

        \pgfmathsetseed{52213446}
        \foreach \x in {1,...,200}
        {
		\pgfmathsetmacro{\x}{rand/2.5+3.5}
		\pgfmathsetmacro{\y}{rand+2}
		\fill[black]    (\x,\y) circle (0.015);
        };

        \pgfmathsetseed{52213446}
        \foreach \x in {1,...,200}
        {
		\pgfmathsetmacro{\x}{rand/4+1.25}
		\pgfmathsetmacro{\y}{rand+1.75}
		\fill[black]    (\x,\y) circle (0.015);
        };

        \draw[->, line width=0.35mm] (2,1.7) -- (-7,1.7);
        \node at (-3,2.2) {Zoomed in view};

        \pgfmathsetseed{58678246}
        \foreach \x in {1,...,80}
        {
		\pgfmathsetmacro{\x}{1.5*rand+2.5}
		\pgfmathsetmacro{\y}{rand/2+0.5}
		\fill[blue]    (\x,\y) circle (0.035);
        };

        \pgfmathsetseed{58678244}
        \foreach \x in {1,...,80}
        {
		\pgfmathsetmacro{\x}{1.5*rand+2.5}
		\pgfmathsetmacro{\y}{rand/4+2.75}
		\fill[blue]    (\x,\y) circle (0.035);
        };

        \pgfmathsetseed{58678246}
        \foreach \x in {1,...,80}
        {
		\pgfmathsetmacro{\x}{rand/2.5+3.5}
		\pgfmathsetmacro{\y}{rand+2}
		\fill[blue]    (\x,\y) circle (0.035);
        };

        \pgfmathsetseed{58678248}
        \foreach \x in {1,...,40}
        {
		\pgfmathsetmacro{\x}{rand/4+1.25}
		\pgfmathsetmacro{\y}{rand+1.75}
		\fill[blue]    (\x,\y) circle (0.035);
        };

    \end{tikzpicture}
    \label{far_pivots_b2t}
    } \qquad \qquad
    \subfloat{
        \begin{tikzpicture}
            [
            box/.style={rectangle,draw=black, minimum size=0.25cm},scale=0.2
            ]
            \node[box,fill=red,,font=\tiny,label=right:Cluster considered,  anchor=west] at (-4,8){};
            \node[box,fill=green,,font=\tiny,label=right:Far-field interaction of the cluster,  anchor=west] at (-4,6){};
            \node[circle,fill=black,inner sep=0pt,minimum size=2pt,label=right:{Pivot}] (a) at (-4,4) {};
        \end{tikzpicture}
    }    
    \caption{Illustration of search spaces of the pivots in B$2$T NCA for $\mathcal{H}^2_{\sqrt{d}}$(b) construction.}
    \label{fig:B2T_pivot}
\end{center}
\end{figure}

As we discussed before, these four sets are the main components for constructing the operators. We construct different operators as described in \Cref{operator_construction}. Therefore, we get the following sets of operators:
 \begin{enumerate}
     \item P$2$M $\bkt{V_{X}}$ / M$2$M $\bkt{\Tilde{V}_{X X_c}}$, \quad $X_c \in child \bkt{X}$.
     \item M$2$L $\bkt{T_{X,Y}}, \quad Y \in \mathcal{IL}_{\sqrt{d}} \bkt{X}$.
     \item L$2$L $\bkt{\Tilde{U}_{X_c X}}$, \quad $X \in parent \bkt{X_c}$ / L$2$P $\bkt{U_{X}}$.
 \end{enumerate}

\subsection{Calculation of the potential (MVP)} \label{b2t_h2matrix_mvp}
The procedure to compute the potential is as follows:

\begin{enumerate}
    \item \textbf{Upward traversal:} 
        \begin{itemize}
            \item Particles to multipole \emph{(P$2$M) at leaf level $\kappa$} : For all leaf clusters $X$, calculate \\ $v^{\bkt{\kappa}}_{X} = V_{X}^{*} q_X^{\bkt{\kappa}}$
            \item Multipole to multipole \emph{(M$2$M) at non-leaf level} : For all non-leaf $X$ clusters, calculate \\ $v^{\bkt{l}}_{X} = \dsum_{X_c \in child \bkt{X}}  \Tilde{V}_{X X_c}^{*} v_{X_c}^{\bkt{l+1}} , \quad \kappa-1 \geq l \geq 2$
        \end{itemize}
    \item  \textbf{Transverse traversal:}
        \begin{itemize}
            \item Multipole to local \emph{(M$2$L) at all levels and for all clusters} : For all cluster $X$, calculate \\
             $u^{\bkt{l}}_{X} = \dsum_{Y \in \mathcal{IL}_{\sqrt{d}} \bkt{X}}  T_{X,Y} v_{Y}^{\bkt{l}} , \quad 2 \leq l \leq \kappa$ 
        \end{itemize}
    \item  \textbf{Downward traversal:}
        \begin{itemize}
            \item Local to local \emph{(L$2$L) at non-leaf level} : For all non-leaf clusters $X$, calculate \\
             $u^{\bkt{l+1}}_{X_c} := u^{\bkt{l+1}}_{X_c} +    \Tilde{U}_{X_c X} u_{X}^{\bkt{l}} , \quad 2 \leq l \leq \kappa-1$ and $X \in parent \bkt{X_c} $
            \item Local to particles \emph{(L$2$P) at leaf level $\kappa$} : For all leaf clusters $X$, calculate \\ $\phi ^{\bkt{\kappa}}_{X} = U_{X} u^{\bkt{\kappa}}_{X}$
        \end{itemize}
\end{enumerate}

\textbf{Near-field potential and the total potential at leaf level.}
For each leaf cluster $X$, we add the near-field (neighbor$+$self) potential, which is a direct computation. Hence, the final computed potential of a leaf cluster is given by
\begin{align}
    \phi ^{\bkt{\kappa}}_{X} := \phi ^{\bkt{\kappa}}_{X} +   \underbrace{\dsum_{X' \in \mathcal{N}_{\sqrt{d}} \bkt{X}}  K_{t^X,s^{X'}} q_{X'}^{\bkt{\kappa}}}_{\text{ Near-field potential (direct)}}
\end{align}
We know that $X_i$ is the $i^{th}$ leaf cluster and $\phi ^{\bkt{\kappa}}_{X_i}$ represents the potential corresponding to it, $1 \leq i \leq 2^{d \kappa}$. 

Therefore, the computed potential is given by $\pmb{\Tilde{\phi}} =\bigg [\phi ^{\bkt{\kappa}}_{X_1};\phi ^{\bkt{\kappa}}_{X_2}; \cdots ;\phi ^{\bkt{\kappa}}_{X_{2^{d \kappa}}} \bigg ]$ (MATLAB notation).

We refer the reader to Algorithm $2$ of  \cite{zhao2019fast} and \cite{gujjula2022new} for more details.

\section{Purely algebraic \texorpdfstring{$\mathcal{H}^2_{\sqrt{d}}$}{H2}(t) algorithm \cite{zhao2019fast}} \label{t2b_h2matrix} 
In this section, we discuss a purely algebraic standard $\mathcal{H}^2$ MVP algorithm. Since the construction of the standard $\mathcal{H}^2$ matrix representation is based on the T$2$B NCA, we denote this algorithm as $\mathcal{H}^2_{\sqrt{d}}$(t). 
We discuss the $\mathcal{H}^2_{\sqrt{d}}$(t) MVP algorithm in the following steps:
\begin{enumerate}
    \item Initialization of $\mathcal{H}^2_{\sqrt{d}}$(t) representation. (\Cref{t2b_h2matrix_init})
    \item Calculation of the potential (MVP). (\Cref{t2b_h2matrix_mvp})
\end{enumerate}

\subsection{Initialization of \texorpdfstring{$\mathcal{H}^2_{\sqrt{d}}$(t)}{H2} representation} \label{t2b_h2matrix_init}
In the T$2$B pivot selection, the pivots of a cluster at a child level are obtained from its own index set and the pivots at its parent level. Therefore, one needs to traverse the $2^d$ uniform tree from top to bottom direction to find pivots of all clusters. The detailed procedure to obtain the four sets of pivots corresponding to a cluster $X$ is given below. 
\begin{itemize}
        \item If $X$ has no parent, i.e., parent($X$) = NULL (parentless), then construct the following four sets
\end{itemize}
    \begin{align}
        \Tilde{t}^{X,i} := t^X \qquad \qquad \text{ and } \qquad \qquad \Tilde{s}^{X,i} := \bigcup_{Y \in \mathcal{IL}_{\sqrt{d}}(X)} s^Y
    \end{align}
     \begin{align}
        \Tilde{t}^{X,o} := \bigcup_{Y \in \mathcal{IL}_{\sqrt{d}}(X)} t^Y \qquad \qquad \text{ and } \qquad \qquad \Tilde{s}^{X,o} := s^X 
    \end{align}
    \begin{itemize}
        \item If $X$ has parent, i.e., parent($X$) $\neq$ NULL, then construct the following four sets
    \end{itemize}
    \begin{align}
        \Tilde{t}^{X,i} := t^X \qquad \qquad \text{ and } \qquad \qquad \Tilde{s}^{X,i} := \bigcup_{Y \in \mathcal{IL}_{\sqrt{d}}(X)} s^{Y} \bigcup s^{parent(X),i}
    \end{align}
     \begin{align}
        \Tilde{t}^{X,o} := \bigcup_{Y \in \mathcal{IL}_{\sqrt{d}}(X)} t^Y \bigcup t^{parent(X),o} \qquad \qquad \text{ and } \qquad \qquad \Tilde{s}^{X,o} := s^X 
    \end{align}
    
To obtain the pivots $t^{X,i}$, $s^{X,i}$, $t^{X,o}$ and $s^{X,o}$, one needs to perform ACA \cite{aca}. We perform ACA on the matrix $K_{\Tilde{t}^{X,i}, \Tilde{s}^{X,i}}$ with user-given tolerance $\epsilon$. The sets $t^{X,i}$ and $s^{X,i}$ are the row and column pivots chosen by the ACA. The search spaces of the pivots for a particular cluster are illustrated in \Cref{fig:T2B_pivot}. Similarly, perform ACA on the matrix $K_{\Tilde{t}^{X,o}, \Tilde{s}^{X,o}}$ to obtain the other two sets of pivots $t^{X,o}$ and $s^{X,o}$.

\begin{figure}[H]
\begin{center}
\captionsetup[subfloat]{labelformat=empty}
    \subfloat[]{
    \begin{tikzpicture}[scale=0.5]
        \draw[draw=black,fill=red] (-8,1) rectangle (-7,2);
        \pgfmathsetseed{58678246}
        \foreach \x in {1,...,200}
        {
		\pgfmathsetmacro{\x}{rand/2 -7.5}
		\pgfmathsetmacro{\y}{rand/2 +1.5}
		\fill[black]    (\x,\y) circle (0.015);
        };
        \pgfmathsetseed{58678246}
        \foreach \x in {1,...,80}
        {
		\pgfmathsetmacro{\x}{rand/2 -7.5}
		\pgfmathsetmacro{\y}{rand/2 +1.5}
		\fill[blue]    (\x,\y) circle (0.035);
        };
        \fill [green!50] (0,0) rectangle (1,1);
        \fill [green!50] (0,1) rectangle (1,2);
        \fill [green!50] (0,2) rectangle (1,3);
        \fill [green!50] (0,3) rectangle (1,4);
        \fill [green!50] (1,3) rectangle (2,4);
        \fill [green!50] (2,3) rectangle (3,4);
        \fill [green!50] (3,3) rectangle (4,4);
        \node[anchor=north] at (1.5,2.8) {};
        \node[anchor=north] at (1.5,1.8) {};
        \node[anchor=north] at (.5,1.8) {};
        \draw[black] (0, 0) grid (4, 4);
        \draw[line width=0.5mm,  fill=green] (1, 0) rectangle (2, 1);
        \draw[line width=0.5mm,  fill=green] (1, 1) rectangle (2, 2);
        \draw[line width=0.5mm,  fill=green] (2, 2) rectangle (3, 3);
        \draw[line width=0.5mm,  fill=green] (3, 2) rectangle (4, 3);
        \draw[line width=0.5mm,  fill=green] (2, 0) rectangle (3, 1);
        \draw[line width=0.5mm,  fill=green] (3, 0) rectangle (4, 1);
        \draw[line width=0.5mm,  fill=green] (3, 1) rectangle (4, 2);
        \draw[line width=0.5mm,  fill=green] (1, 2) rectangle (2, 3);
        \fill [white] (1.5,1) rectangle (3,2.5);
        \fill [red] (2,1.5) rectangle (2.5,2);
        \node[anchor=north] at (2.27,2.15) {\tiny $X$};
        \draw[step=0.5cm, black] (1, 0) grid (4, 3);
        \draw[line width=0.5mm] (1, 0) rectangle (2, 1);
        \draw[line width=0.5mm] (1, 1) rectangle (2, 2);
        \draw[line width=0.5mm] (2, 2) rectangle (3, 3);
        \draw[line width=0.5mm] (3, 2) rectangle (4, 3);
        \draw[line width=0.5mm] (2, 0) rectangle (3, 1);
        \draw[line width=0.5mm] (3, 0) rectangle (4, 1);
        \draw[line width=0.5mm] (3, 1) rectangle (4, 2);
        \draw[line width=0.5mm] (1, 2) rectangle (2, 3);

        \pgfmathsetseed{58678246}
        \foreach \x in {1,...,400}
        {
		\pgfmathsetmacro{\x}{rand/2+0.5}
		\pgfmathsetmacro{\y}{2*rand+2}
		\fill[black]    (\x,\y) circle (0.015);
        };

        \pgfmathsetseed{52213446}
        \foreach \x in {1,...,400}
        {
		\pgfmathsetmacro{\x}{2*rand+2}
		\pgfmathsetmacro{\y}{rand/2+3.5}
		\fill[black]    (\x,\y) circle (0.015);
        };

        \pgfmathsetseed{58678226}
        \foreach \x in {1,...,200}
        {
		\pgfmathsetmacro{\x}{1.5*rand+2.5}
		\pgfmathsetmacro{\y}{rand/2+0.5}
		\fill[black]    (\x,\y) circle (0.015);
        };

        \pgfmathsetseed{52213436}
        \foreach \x in {1,...,200}
        {
		\pgfmathsetmacro{\x}{1.5*rand+2.5}
		\pgfmathsetmacro{\y}{rand/4+2.75}
		\fill[black]    (\x,\y) circle (0.015);
        };

        \pgfmathsetseed{52213448}
        \foreach \x in {1,...,200}
        {
		\pgfmathsetmacro{\x}{rand/2.5+3.5}
		\pgfmathsetmacro{\y}{rand+2}
		\fill[black]    (\x,\y) circle (0.015);
        };

        \pgfmathsetseed{52213456}
        \foreach \x in {1,...,200}
        {
		\pgfmathsetmacro{\x}{rand/4+1.25}
		\pgfmathsetmacro{\y}{rand+1.75}
		\fill[black]    (\x,\y) circle (0.015);
        };

        \draw[->, line width=0.35mm] (2,1.7) -- (-7,1.7);
        \node at (-3,2.2) {Zoomed in view};

        \pgfmathsetseed{58678248}
        \foreach \x in {1,...,80}
        {
		\pgfmathsetmacro{\x}{rand/2+0.5}
		\pgfmathsetmacro{\y}{2*rand+2}
		\fill[blue]    (\x,\y) circle (0.035);
        };

        \pgfmathsetseed{52213446}
        \foreach \x in {1,...,80}
        {
		\pgfmathsetmacro{\x}{2*rand+2}
		\pgfmathsetmacro{\y}{rand/2+3.5}
		\fill[blue]    (\x,\y) circle (0.035);
        };

        \pgfmathsetseed{58678246}
        \foreach \x in {1,...,80}
        {
		\pgfmathsetmacro{\x}{1.5*rand+2.5}
		\pgfmathsetmacro{\y}{rand/2+0.5}
		\fill[blue]    (\x,\y) circle (0.035);
        };

        \pgfmathsetseed{52213446}
        \foreach \x in {1,...,80}
        {
		\pgfmathsetmacro{\x}{1.5*rand+2.5}
		\pgfmathsetmacro{\y}{rand/4+2.75}
		\fill[blue]    (\x,\y) circle (0.035);
        };

        \pgfmathsetseed{52213446}
        \foreach \x in {1,...,80}
        {
		\pgfmathsetmacro{\x}{rand/2.5+3.5}
		\pgfmathsetmacro{\y}{rand+2}
		\fill[blue]    (\x,\y) circle (0.035);
        };

        \pgfmathsetseed{52213446}
        \foreach \x in {1,...,80}
        {
		\pgfmathsetmacro{\x}{rand/4+1.25}
		\pgfmathsetmacro{\y}{rand+1.75}
		\fill[blue]    (\x,\y) circle (0.035);
        };

    \end{tikzpicture}
    \label{far_pivots_t2b}
    }\qquad \qquad
    \subfloat{
        \begin{tikzpicture}
            [
            box/.style={rectangle,draw=black, minimum size=0.25cm},scale=0.2
            ]
            \node[box,fill=red,,font=\tiny,label=right:Cluster considered,  anchor=west] at (-4,8){};
            \node[box,fill=green!50,,font=\tiny,label=right:Far-field interaction at parent level,  anchor=west] at (-4,6){};
            \node[box,fill=green,,font=\tiny,label=right:Far-field interaction of the cluster,  anchor=west] at (-4,4){};
            \node[circle,fill=black,inner sep=0pt,minimum size=2pt,label=right:{Pivot}] (a) at (-4,2) {};
        \end{tikzpicture}
    }   
    \caption{Illustration of search spaces of the pivots in T$2$B NCA for $\mathcal{H}^2_{\sqrt{d}}$(t) construction.}
    \label{fig:T2B_pivot}
\end{center}
\end{figure}

As discussed before, these four sets are the main components for constructing the operators. We construct different operators as described in \Cref{operator_construction}. Therefore, we get the following sets of operators:
 \begin{enumerate}
     \item P$2$M $\bkt{V_{X}}$ / M$2$M $\bkt{\Tilde{V}_{X X_c}}$, \quad $X_c \in child \bkt{X}$.
     \item M$2$L $\bkt{T_{X,Y}}, \quad Y \in \mathcal{IL}_{\sqrt{d}} \bkt{X}$.
     \item L$2$L $\bkt{\Tilde{U}_{X_c X}}$, \quad $X \in parent \bkt{X_c}$ / L$2$P $\bkt{U_{X}}$.
 \end{enumerate}

\subsection{Calculation of the potential (MVP)} \label{t2b_h2matrix_mvp}
We perform the potential (MVP) calculation by following upward, transverse and downward tree traversal, which is similar to \Cref{b2t_h2matrix_mvp}.

We refer the reader to Algorithm $1$ of \cite{zhao2019fast} for more details.

\section{\texorpdfstring{$\mathcal{H}^2_{*}$}{H2*}(t): T\texorpdfstring{$2$}{2}B NCA on our \emph{weak admissibility} condition} \label{t2b_nHODt} 
 This section discusses a fast MVP algorithm when applying the T$2$B NCA on our \emph{weak admissibility} condition in higher dimensions (\Cref{weak_admis}), i.e., the admissible clusters are the far-field and vertex-sharing clusters. We denote this algorithm as $\mathcal{H}^2_{*}$(t). We discuss the $\mathcal{H}^2_{*}$(t) algorithm in the following steps: 
\begin{enumerate}
    \item Initialization of $\mathcal{H}^2_{*}$(t) representation. (\Cref{t2b_nHODt_init})
    \item Calculation of the potential (MVP). (\Cref{t2b_nHODt_mvp})
\end{enumerate}

\subsection{Initialization of \texorpdfstring{$\mathcal{H}^2_{*}$}{H2*}(t) representation} \label{t2b_nHODt_init}
In the T$2$B pivot selection, the pivots of a cluster at a child level are obtained from its own index set and the pivots at its parent level. Therefore, one needs to traverse the $2^d$ uniform tree from top to bottom direction to find pivots of all clusters. The detailed procedure to obtain the pivots corresponding to a cluster $X$ is given below. 
\begin{itemize}
        \item If $X$ has no parent, i.e., parent($X$) = NULL (parentless), then construct the following four sets
\end{itemize}
    \begin{align}
        \Tilde{t}^{X,i} := t^X \qquad \qquad \text{ and } \qquad \qquad \Tilde{s}^{X,i} := \bigcup_{Y \in \mathcal{IL}_{*}(X)} s^Y
    \end{align}
     \begin{align}
        \Tilde{t}^{X,o} := \bigcup_{Y \in \mathcal{IL}_{*}(X)} t^Y \qquad \qquad \text{ and } \qquad \qquad \Tilde{s}^{X,o} := s^X 
    \end{align}
    \begin{itemize}
        \item If $X$ has parent, i.e., parent($X$) $\neq$ NULL, then construct the following four sets
    \end{itemize}
    \begin{align}
        \Tilde{t}^{X,i} := t^X \qquad \qquad \text{ and } \qquad \qquad \Tilde{s}^{X,i} := \bigcup_{Y \in \mathcal{IL}_{*}(X)} s^{Y} \bigcup s^{parent(X),i}
    \end{align}
     \begin{align}
        \Tilde{t}^{X,o} := \bigcup_{Y \in \mathcal{IL}_{*}(X)} t^Y \bigcup t^{parent(X),o} \qquad \qquad \text{ and } \qquad \qquad \Tilde{s}^{X,o} := s^X 
    \end{align}
    
To obtain the pivots $t^{X,i}$, $s^{X,i}$, $t^{X,o}$ and $s^{X,o}$, one needs to perform ACA \cite{aca}. We perform ACA on the matrix $K_{\Tilde{t}^{X,i}, \Tilde{s}^{X,i}}$ with user-given tolerance $\epsilon$. The sets $t^{X,i}$ and $s^{X,i}$ are the row and column pivots chosen by the ACA. The search spaces of the pivots for a particular cluster are illustrated in \Cref{fig:T2B_pivot_entire_interaction}. Similarly, perform ACA on the matrix $K_{\Tilde{t}^{X,o}, \Tilde{s}^{X,o}}$ to obtain the other two sets of pivots $t^{X,o}$ and $s^{X,o}$.

These four sets are the main components for constructing the operators. We construct different operators as described in \Cref{operator_construction}. Therefore, we get the following sets of operators:
 \begin{enumerate}
     \item P$2$M $\bkt{V_{X}}$ / M$2$M $\bkt{\Tilde{V}_{X X_c}}$, \quad $X_c \in child \bkt{X}$.
     \item M$2$L $\bkt{T_{X,Y}}, \quad Y \in \mathcal{IL}_{*} \bkt{X}$.
     \item L$2$L $\bkt{\Tilde{U}_{X_c X}}$, \quad $X \in parent \bkt{X_c}$ / L$2$P $\bkt{U_{X}}$.
 \end{enumerate}
 
\subsection{Calculation of the potential (MVP)} \label{t2b_nHODt_mvp}
The procedure to compute the potential is as follows:

\begin{enumerate}
    \item \textbf{Upward traversal:} 
        \begin{itemize}
            \item Particles to multipole \emph{(P$2$M) at leaf level $\kappa$} : For all leaf clusters $X$, calculate \\ $v^{\bkt{\kappa}}_{X} = V_{X}^{*} q_X^{\bkt{\kappa}}$
            \item Multipole to multipole \emph{(M$2$M) at non-leaf level} : For all non-leaf $X$ clusters, calculate \\ $v^{\bkt{l}}_{X} = \dsum_{X_c \in child \bkt{X}}  \Tilde{V}_{X X_c}^{*} v_{X_c}^{\bkt{l+1}} , \quad \kappa-1 \geq l \geq 1$
        \end{itemize}
    \item  \textbf{Transverse traversal:}
        \begin{itemize}
            \item Multipole to local \emph{(M$2$L) at all levels and for all clusters} : For all cluster $X$, calculate \\
             $u^{\bkt{l}}_{X} = \dsum_{Y \in \mathcal{IL}_{*} \bkt{X}}  T_{X,Y} v_{Y}^{\bkt{l}} , \quad 1 \leq l \leq \kappa$ 
        \end{itemize}
    \item  \textbf{Downward traversal:}
        \begin{itemize}
            \item Local to local \emph{(L$2$L) at non-leaf level} : For all non-leaf clusters $X$, calculate \\
             $u^{\bkt{l+1}}_{X_c} := u^{\bkt{l+1}}_{X_c} +    \Tilde{U}_{X_c X} u_{X}^{\bkt{l}} , \quad 1 \leq l \leq \kappa-1$ and $X \in parent \bkt{X_c} $
            \item Local to particles \emph{(L$2$P) at leaf level $\kappa$} : For all leaf clusters $X$, calculate \\ $\phi ^{\bkt{\kappa}}_{X} = U_{X} u^{\bkt{\kappa}}_{X}$
        \end{itemize}
\end{enumerate}

\textbf{Near-field potential and the total potential at leaf level.}
For each leaf cluster $X$, we add the near-field (neighbor$+$self) potential, which is a direct computation. Hence, the final computed potential of a leaf cluster is given by
\begin{align}
    \phi ^{\bkt{\kappa}}_{X} := \phi ^{\bkt{\kappa}}_{X} +   \underbrace{\dsum_{X' \in \mathcal{N}_{*} \bkt{X}}  K_{t^X,s^{X'}} q_{X'}^{\bkt{\kappa}}}_{\text{ Near-field potential (direct)}}
\end{align}
We know that $X_i$ is the $i^{th}$ leaf cluster and $\phi ^{\bkt{\kappa}}_{X_i}$ represents the potential corresponding to it, $1 \leq i \leq 2^{d \kappa}$. 

\begin{figure}[H]
    \centering
    \captionsetup[subfloat]{labelformat=empty}
    \subfloat{
        \begin{tikzpicture}
        \node[anchor=west] at (0,-1.25){$\Tilde{\phi} = $};
        \filldraw  (1,0) to[out=260,in=100]  (1,-2.5) to [out=98,in=262] cycle;
        \end{tikzpicture}
    }
    \subfloat[]{
        \includegraphics[scale=0.25]{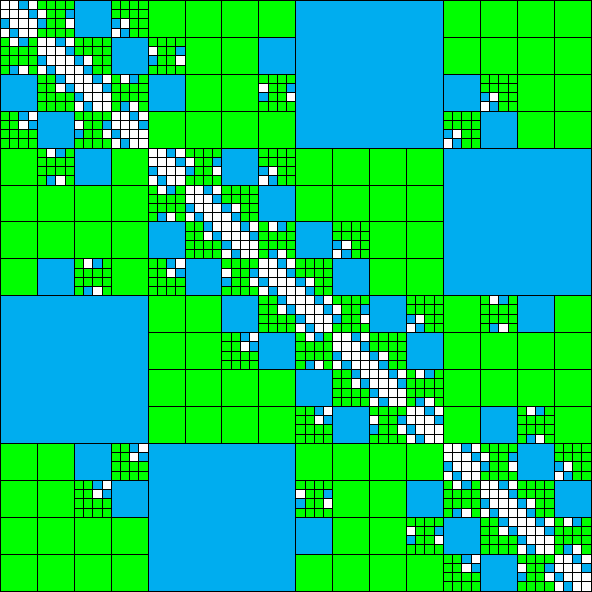}
        }\quad
    \subfloat[]{
        \begin{tikzpicture}
        \node at (-3,0) {};
        \node at (-3,1.1) {$+$};
        \end{tikzpicture}
    }\quad
    \subfloat[]{
        \includegraphics[scale=0.25]{HODLR_images/HODLR2d_nbd_3.pdf}
        }
    \subfloat{
        \begin{tikzpicture}
        \filldraw (2.75,0) to[out=-80,in=80]  (2.75,-2.5) to [out=82,in=-82] cycle;
        \end{tikzpicture}
    }\quad
    \subfloat{
        \begin{tikzpicture}
        \node[anchor=west] at (-1,-2.75){$\times$};
        \draw[draw=black, fill=blue!40] (0, -4) rectangle (0.12, -1.5);
        \node[anchor=west] at (-0.15,-2.75){$\pmb{q}$};
        \end{tikzpicture}
    }\quad
    \caption{In the $\mathcal{H}^2_{*}$(t) algorithm, the operators (P$2$M/M$2$M, M$2$L and L$2$L/L$2$P) are constructed using T$2$B NCA. After that, we calculate the potential using Upward, Transverse and Downward tree traversal. The near-field potential is added at the leaf level.}
    \label{fig:nHOD_t_construction}
\end{figure}

Therefore, the computed potential is given by $\pmb{\Tilde{\phi}} =\bigg [\phi ^{\bkt{\kappa}}_{X_1};\phi ^{\bkt{\kappa}}_{X_2}; \cdots ;\phi ^{\bkt{\kappa}}_{X_{2^{d \kappa}}} \bigg ]$ (MATLAB notation). The visual representation of the $\mathcal{H}^2_{*}$(t) algorithm is given in \Cref{fig:nHOD_t_construction}.

\section{Additional experiments in \texorpdfstring{$3$}{3}D}
In this section, we perform additional experiments using different kernel matrices and compare the performance of the proposed algorithms with the standard $\mathcal{H}^2$ matrix algorithm \cite{zhao2019fast}.
\subsection{Fast MVP with different kernel functions} \label{other_kernels}
We consider $N$ uniformly distributed particles with location at $\{\pmb{x}_i \}_{i=1}^N$ in the cube $[-1,1]^3$. The kernel matrix $(K)$ is generated by the Matérn covariance kernel and Helmholtz kernel. The fast MVP $(\Tilde{K} \pmb{q})$ is performed using different nested hierarchical representations (refer to \Cref{tab:nex_notation}).  

\subsubsection{Matérn covariance kernel with \texorpdfstring{$\nu=1/2$}{nu equal to 1/2}}
We choose the kernel function as the Matérn covariance kernel. The $(i,j)^{th}$ entry of the kernel matrix $K \in \mathbb{R}^{N \times N}$ is given by
\begin{equation}
K(i,j) =  \exp{(-r)} = \exp{\bkt{-\magn{\pmb{x}_i - \pmb{x}_j}_2}} 
\end{equation}

\begin{figure}[H]
    \centering
    \subfloat[]{\includegraphics[height=4.5cm, width=6cm]{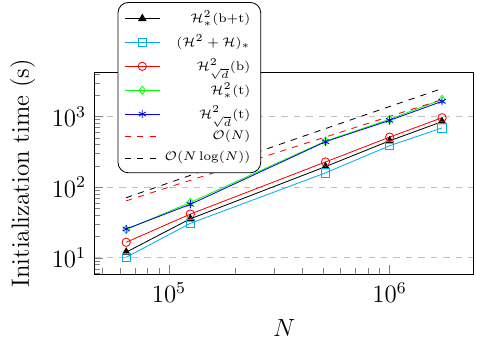}\label{matern3d_assem}} \qquad
    \subfloat[]{\includegraphics[height=4.5cm, width=6cm]{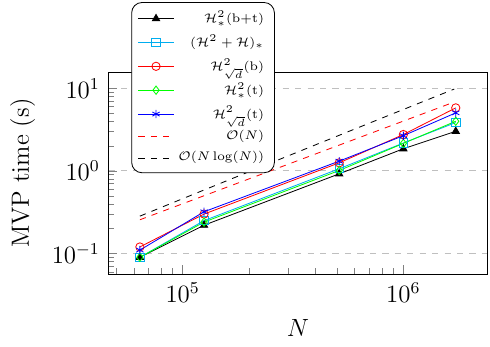}\label{matern3d_mvp}}
    \caption{Plots of Initialization time and MVP time of different fast MVP algorithms.}
    \label{fig:matern3d_1}
\end{figure}

\begin{figure}[H]
    \centering
    \subfloat[]{\includegraphics[height=4.5cm, width=6cm]{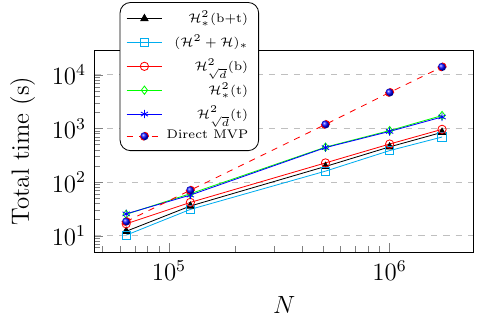}\label{matern3d_total}} \qquad
    \subfloat[]{\includegraphics[height=4.5cm, width=6cm]{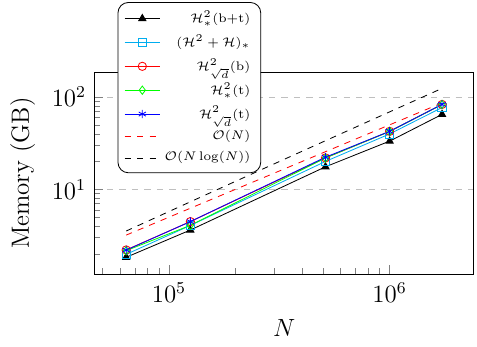}\label{matern3d_mem}}
    \caption{Plots of total time $(t_{init} + t_{MVP})$ and memory of the algorithms. We also plot the time for Direct MVP.}
    \label{fig:matern3d_2}
\end{figure}

We randomly select $5$ different column vectors $(\pmb{q})$ and perform the fast MVPs $(\Tilde{K} \pmb{q})$ using $\mathcal{H}^2_{*}$(b$+$t), $\bkt{\mathcal{H}^2 + \mathcal{H}}_{*}$, $\mathcal{H}^2_{\sqrt{d}}$(b), $\mathcal{H}^2_{*}$(t) and $\mathcal{H}^2_{\sqrt{d}}$(t) algorithms. The average initialization time, MVP time, memory and relative error in MVP $({RE}_{MVP})$ are tabulated in \Cref{table:matern3d_a} and \Cref{table:matern3d_b}. We also plot the initialization time (\Cref{matern3d_assem}), MVP time (\Cref{matern3d_mvp}), total time (\Cref{matern3d_total}) and memory (\Cref{matern3d_mem}). Please refer to \Cref{fig:matern3d_1} and \Cref{fig:matern3d_2}.

\begin{table}[H]
      \centering
      \resizebox{\textwidth}{!}{%
      \setlength\extrarowheight{0.9pt}
\begin{tabular}{|l|lllll|lllll|}
\hline
\multirow{2}{*}{$N$} & \multicolumn{5}{l|}{\qquad \qquad Initialization time (s)}                                                                                           & \multicolumn{5}{l|}{\qquad \qquad  \qquad \quad MVP time (s)}                                                                                \\ \cline{2-11} 
                   & \multicolumn{1}{l|}{\scriptsize $\mathcal{H}^2_{*}$(b$+$t)} & \multicolumn{1}{l|}{\scriptsize $\bkt{\mathcal{H}^2 + \mathcal{H}}_{*}$} & \multicolumn{1}{l|}{\scriptsize $\mathcal{H}^2_{\sqrt{d}}$(b)} & \multicolumn{1}{l|}{\scriptsize $\mathcal{H}^2_{*}$(t)} &\scriptsize $\mathcal{H}^2_{\sqrt{d}}$(t) & \multicolumn{1}{l|}{\scriptsize $\mathcal{H}^2_{*}$(b$+$t)} & \multicolumn{1}{l|}{\scriptsize $\bkt{\mathcal{H}^2 + \mathcal{H}}_{*}$} & \multicolumn{1}{l|}{\scriptsize $\mathcal{H}^2_{\sqrt{d}}$(b)} & \multicolumn{1}{l|}{\scriptsize $\mathcal{H}^2_{*}$(t)} & \scriptsize $\mathcal{H}^2_{\sqrt{d}}$(t) \\ \hline
64000   &\multicolumn{1}{|l|}{12.04}     & \multicolumn{1}{l|}{10.18}       & \multicolumn{1}{l|}{16.59}      & \multicolumn{1}{l|}{25.15}        &  25.59     & \multicolumn{1}{l|}{0.09}     & \multicolumn{1}{l|}{0.09}       & \multicolumn{1}{l|}{0.12}      & \multicolumn{1}{l|}{0.09}        &   0.11    \\ \hline
125000  &\multicolumn{1}{|l|}{35.65}     & \multicolumn{1}{l|}{30.77}       & \multicolumn{1}{l|}{41.58}         & \multicolumn{1}{l|}{60.64}        & 57.17      & \multicolumn{1}{l|}{0.22}     & \multicolumn{1}{l|}{0.25}       & \multicolumn{1}{l|}{0.3}      & \multicolumn{1}{l|}{0.24}        &  0.32      \\ \hline
512000  &\multicolumn{1}{|l|}{195.09}     & \multicolumn{1}{l|}{159.19}       & \multicolumn{1}{l|}{226.98}      & \multicolumn{1}{l|}{448.17}    &  439.29    & \multicolumn{1}{l|}{0.92}     & \multicolumn{1}{l|}{1.06}       & \multicolumn{1}{l|}{1.25}      & \multicolumn{1}{l|}{1.01}        &  1.31     \\ \hline
1000000 &\multicolumn{1}{|l|}{452.44}     & \multicolumn{1}{l|}{385.12}       & \multicolumn{1}{l|}{510.36}      & \multicolumn{1}{l|}{911.33}    &  879.86    & \multicolumn{1}{l|}{1.85}     & \multicolumn{1}{l|}{2.17}       & \multicolumn{1}{l|}{2.75}      & \multicolumn{1}{l|}{2.15}        &  2.66     \\ \hline
1728000 &\multicolumn{1}{|l|}{844.16}     & \multicolumn{1}{l|}{688.95}       & \multicolumn{1}{l|}{958.17}      & \multicolumn{1}{l|}{1738.65}   &  1635.63   & \multicolumn{1}{l|}{3.01}     & \multicolumn{1}{l|}{3.84}       & \multicolumn{1}{l|}{5.78}      & \multicolumn{1}{l|}{3.98}        &  5.03     \\ \hline
\end{tabular}
}
    \caption{Performance of the proposed and the nested algorithms as mentioned in \Cref{tab:nex_notation} in terms of the initialization time and MVP time.}
    \label{table:matern3d_a}
\end{table}

\begin{table}[H]
      \centering
      \resizebox{\textwidth}{!}{%
      \setlength\extrarowheight{0.9pt}
\begin{tabular}{|l|lllll|lllll|}
\hline
\multirow{2}{*}{$N$} & \multicolumn{5}{l|}{\qquad \qquad \qquad Memory (GB)}                                                                                           & \multicolumn{5}{l|}{\qquad \qquad  Relative error in MVP ($2$-norm)}                                                                                \\ \cline{2-11} 
                   & \multicolumn{1}{l|}{\scriptsize $\mathcal{H}^2_{*}$(b$+$t)} & \multicolumn{1}{l|}{\scriptsize $\bkt{\mathcal{H}^2 + \mathcal{H}}_{*}$} & \multicolumn{1}{l|}{\scriptsize $\mathcal{H}^2_{\sqrt{d}}$(b)} & \multicolumn{1}{l|}{\scriptsize $\mathcal{H}^2_{*}$(t)} &\scriptsize $\mathcal{H}^2_{\sqrt{d}}$(t) & \multicolumn{1}{l|}{\scriptsize $\mathcal{H}^2_{*}$(b$+$t)} & \multicolumn{1}{l|}{\scriptsize $\bkt{\mathcal{H}^2 + \mathcal{H}}_{*}$} & \multicolumn{1}{l|}{\scriptsize $\mathcal{H}^2_{\sqrt{d}}$(b)} & \multicolumn{1}{l|}{\scriptsize $\mathcal{H}^2_{*}$(t)} & \scriptsize $\mathcal{H}^2_{\sqrt{d}}$(t) \\ \hline
64000    &\multicolumn{1}{|l|}{1.84}     & \multicolumn{1}{l|}{1.98}       & \multicolumn{1}{l|}{2.21}      & \multicolumn{1}{l|}{2.17}        &  2.2      & \multicolumn{1}{l|}{8.13E-07}     & \multicolumn{1}{l|}{5.78E-07}       & \multicolumn{1}{l|}{4.93E-07}      & \multicolumn{1}{l|}{1.10E-06}        &   5.64E-07    \\ \hline
125000   &\multicolumn{1}{|l|}{3.63}     & \multicolumn{1}{l|}{4.07}       & \multicolumn{1}{l|}{4.48}      & \multicolumn{1}{l|}{4.08}        &  4.48     & \multicolumn{1}{l|}{1.53E-06}     & \multicolumn{1}{l|}{7.32E-07}       & \multicolumn{1}{l|}{7.89E-07}      & \multicolumn{1}{l|}{2.20E-06}        &  8.62E-07     \\ \hline
512000   &\multicolumn{1}{|l|}{17.56}     & \multicolumn{1}{l|}{20.08}       & \multicolumn{1}{l|}{21.92}      & \multicolumn{1}{l|}{21.89}    &  22.31    & \multicolumn{1}{l|}{2.32E-06}     & \multicolumn{1}{l|}{1.10E-06}       & \multicolumn{1}{l|}{1.01E-06}      & \multicolumn{1}{l|}{3.39E-06}        &  1.13E-06     \\ \hline
1000000  &\multicolumn{1}{|l|}{33.3}     & \multicolumn{1}{l|}{39.83}       & \multicolumn{1}{l|}{42.49}      & \multicolumn{1}{l|}{42.44}     &  42.51    & \multicolumn{1}{l|}{3.35E-06}     & \multicolumn{1}{l|}{1.35E-06}       & \multicolumn{1}{l|}{1.34E-06}      & \multicolumn{1}{l|}{5.47E-06}        &  1.19E-06     \\ \hline
1728000  &\multicolumn{1}{|l|}{64.87}     & \multicolumn{1}{l|}{77.03}       & \multicolumn{1}{l|}{83.12}      & \multicolumn{1}{l|}{83.05}    &  83.34    & \multicolumn{1}{l|}{3.86E-06}     & \multicolumn{1}{l|}{1.28E-06}       & \multicolumn{1}{l|}{1.31E-06}      & \multicolumn{1}{l|}{6.23E-06}        &  1.73E-06     \\ \hline
\end{tabular}
}
    \caption{Performance of the proposed and the nested algorithms as mentioned in \Cref{tab:nex_notation} in terms of the memory and relative error. We set the tolerance $\epsilon = 10^{-6}$.}
    \label{table:matern3d_b}
\end{table}

Let us arrange the above algorithms in ascending order with respect to memory and MVP time.\\
\textbf{Memory:} $\mathcal{H}^2_{*}$(b$+$t) < $\bkt{\mathcal{H}^2 + \mathcal{H}}_{*}$ < $\mathcal{H}^2_{*}$(t) < $\mathcal{H}^2_{\sqrt{d}}$(b) $\leq$ $\mathcal{H}^2_{\sqrt{d}}$(t) \\ 
\textbf{MVP time:} $\mathcal{H}^2_{*}$(b$+$t) < $\bkt{\mathcal{H}^2 + \mathcal{H}}_{*}$ < $\mathcal{H}^2_{*}$(t) < $\mathcal{H}^2_{\sqrt{d}}$(t) < $\mathcal{H}^2_{\sqrt{d}}$(b)

\subsubsection{\texorpdfstring{$3$}{3}D Helmholtz kernel with \texorpdfstring{$k=1$}{k equal to 1}}
We choose the kernel function as the Helmholtz kernel with wave number $k=1$. The $(\alpha, \beta)^{th}$ entry of the kernel matrix $K \in \mathbb{C}^{N \times N}$ is given by
\begin{equation}
K(\alpha,\beta) = \begin{cases}
            \dfrac{\exp{(ir)}}{r} = \dfrac{\exp{\bkt{i\magn{\pmb{x}_{\alpha} - \pmb{x}_{\beta}}_2}}}{\magn{\pmb{x}_{\alpha} - \pmb{x}_{\beta}}_2} & \text{ if } \alpha \neq \beta\\
             0 & \text{ otherwise}
            \end{cases}    
\end{equation}

\begin{figure}[H]
    \centering
    \subfloat[]{\includegraphics[height=4.5cm, width=6cm]{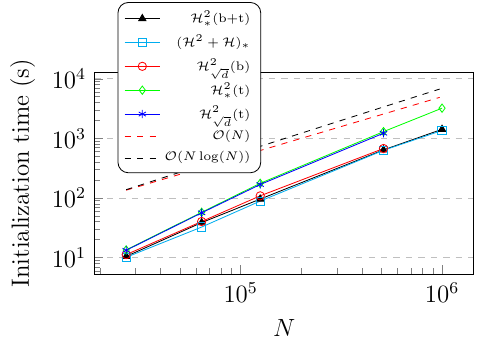}\label{helm_assem}} \qquad
    \subfloat[]{\includegraphics[height=4.5cm, width=6cm]{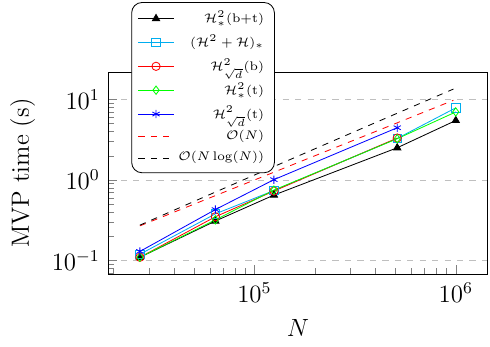}\label{helm_mvp}}
    \caption{Plots of Initialization time and MVP time of different fast MVP algorithms.}
    \label{fig:helm3d_1}
\end{figure}

\begin{figure}[H]
    \centering
    \subfloat[]{\includegraphics[height=4.5cm, width=6cm]{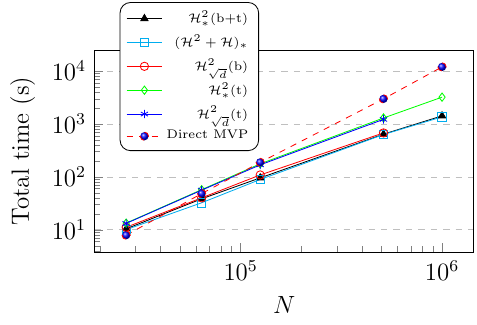}\label{helm_total}} \qquad
    \subfloat[]{\includegraphics[height=4.5cm, width=6cm]{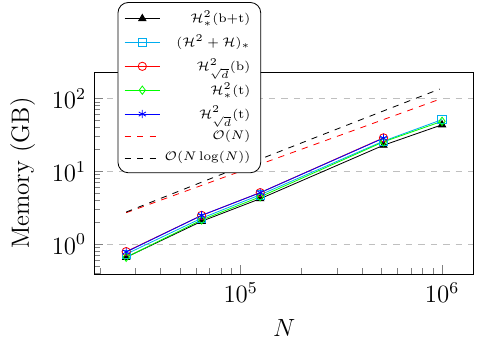}\label{helm_mem}}
    \caption{Plots of total time $(t_{init} + t_{MVP})$ and memory of the algorithms. We also plot the time for Direct MVP.}
    \label{fig:helm3d_2}
\end{figure}

We randomly select $5$ different column vectors $(\pmb{q})$ and perform the fast MVPs $(\Tilde{K} \pmb{q})$ using $\mathcal{H}^2_{*}$(b$+$t), $\bkt{\mathcal{H}^2 + \mathcal{H}}_{*}$, $\mathcal{H}^2_{\sqrt{d}}$(b), $\mathcal{H}^2_{*}$(t) and $\mathcal{H}^2_{\sqrt{d}}$(t) algorithms. The average initialization time, MVP time, memory and relative error in MVP $({RE}_{MVP})$ are tabulated in \Cref{table:helm3d_a} and \Cref{table:helm3d_b}. We also plot the initialization time (\Cref{helm_assem}), MVP time (\Cref{helm_mvp}), total time (\Cref{helm_total}) and memory (\Cref{helm_mem}). Please refer to \Cref{fig:helm3d_1} and \Cref{fig:helm3d_2}.

\begin{table}[H]
      \centering
      \resizebox{\textwidth}{!}{%
      \setlength\extrarowheight{0.9pt}
\begin{tabular}{|l|lllll|lllll|}
\hline
\multirow{2}{*}{$N$} & \multicolumn{5}{l|}{\qquad \qquad Initialization time (s)}                                                                                           & \multicolumn{5}{l|}{\qquad \qquad  \qquad \quad MVP time (s)}                                                                                \\ \cline{2-11} 
                   & \multicolumn{1}{l|}{\scriptsize $\mathcal{H}^2_{*}$(b$+$t)} & \multicolumn{1}{l|}{\scriptsize $\bkt{\mathcal{H}^2 + \mathcal{H}}_{*}$} & \multicolumn{1}{l|}{\scriptsize $\mathcal{H}^2_{\sqrt{d}}$(b)} & \multicolumn{1}{l|}{\scriptsize $\mathcal{H}^2_{*}$(t)} &\scriptsize $\mathcal{H}^2_{\sqrt{d}}$(t) & \multicolumn{1}{l|}{\scriptsize $\mathcal{H}^2_{*}$(b$+$t)} & \multicolumn{1}{l|}{\scriptsize $\bkt{\mathcal{H}^2 + \mathcal{H}}_{*}$} & \multicolumn{1}{l|}{\scriptsize $\mathcal{H}^2_{\sqrt{d}}$(b)} & \multicolumn{1}{l|}{\scriptsize $\mathcal{H}^2_{*}$(t)} & \scriptsize $\mathcal{H}^2_{\sqrt{d}}$(t) \\ \hline
27000   &\multicolumn{1}{|l|}{10.47}     & \multicolumn{1}{l|}{10.12}       & \multicolumn{1}{l|}{11.22}         & \multicolumn{1}{l|}{13.44}        &  13.26     & \multicolumn{1}{l|}{0.11}     & \multicolumn{1}{l|}{0.12}       & \multicolumn{1}{l|}{0.11}      & \multicolumn{1}{l|}{0.11}        &  0.13    \\ \hline
64000  &\multicolumn{1}{|l|}{38.65}     & \multicolumn{1}{l|}{32.55}       & \multicolumn{1}{l|}{40.33}         & \multicolumn{1}{l|}{57.54}        &  56.29     & \multicolumn{1}{l|}{0.31}     & \multicolumn{1}{l|}{0.38}       & \multicolumn{1}{l|}{0.35}      & \multicolumn{1}{l|}{0.32}        &  0.43     \\ \hline
125000  &\multicolumn{1}{|l|}{95.82}     & \multicolumn{1}{l|}{88.65}       & \multicolumn{1}{l|}{108.59}        & \multicolumn{1}{l|}{176.07}       & 168.23     & \multicolumn{1}{l|}{0.65}     & \multicolumn{1}{l|}{0.74}       & \multicolumn{1}{l|}{0.73}      & \multicolumn{1}{l|}{0.75}        &  1.01     \\ \hline
512000 &\multicolumn{1}{|l|}{640.44}     & \multicolumn{1}{l|}{630.25}       & \multicolumn{1}{l|}{665.35}      & \multicolumn{1}{l|}{1295.68}      & 1219.55    & \multicolumn{1}{l|}{2.51}     & \multicolumn{1}{l|}{3.28}       & \multicolumn{1}{l|}{3.31}      & \multicolumn{1}{l|}{3.25}        &  4.45      \\ \hline
1000000 &\multicolumn{1}{|l|}{1418.08}     & \multicolumn{1}{l|}{1355.88}       & \multicolumn{1}{l|}{-}         & \multicolumn{1}{l|}{3208.47}      &   -        & \multicolumn{1}{l|}{5.47}     & \multicolumn{1}{l|}{7.85}       & \multicolumn{1}{l|}{-}      & \multicolumn{1}{l|}{6.98}           &   -    \\ \hline
\end{tabular}
}
    \caption{Performance of the proposed and the nested algorithms as mentioned in \Cref{tab:nex_notation} in terms of the initialization time and MVP time.}
    \label{table:helm3d_a}
\end{table}

\begin{table}[H]
      \centering
      \resizebox{\textwidth}{!}{%
      \setlength\extrarowheight{0.9pt}
\begin{tabular}{|l|lllll|lllll|}
\hline
\multirow{2}{*}{$N$} & \multicolumn{5}{l|}{\qquad \qquad \qquad Memory (GB)}                                                                                           & \multicolumn{5}{l|}{\qquad \qquad  Relative error in MVP ($2$-norm)}                                                                                \\ \cline{2-11} 
                   & \multicolumn{1}{l|}{\scriptsize $\mathcal{H}^2_{*}$(b$+$t)} & \multicolumn{1}{l|}{\scriptsize $\bkt{\mathcal{H}^2 + \mathcal{H}}_{*}$} & \multicolumn{1}{l|}{\scriptsize $\mathcal{H}^2_{\sqrt{d}}$(b)} & \multicolumn{1}{l|}{\scriptsize $\mathcal{H}^2_{*}$(t)} &\scriptsize $\mathcal{H}^2_{\sqrt{d}}$(t) & \multicolumn{1}{l|}{\scriptsize $\mathcal{H}^2_{*}$(b$+$t)} & \multicolumn{1}{l|}{\scriptsize $\bkt{\mathcal{H}^2 + \mathcal{H}}_{*}$} & \multicolumn{1}{l|}{\scriptsize $\mathcal{H}^2_{\sqrt{d}}$(b)} & \multicolumn{1}{l|}{\scriptsize $\mathcal{H}^2_{*}$(t)} & \scriptsize $\mathcal{H}^2_{\sqrt{d}}$(t) \\ \hline
27000   &\multicolumn{1}{|l|}{0.66}     & \multicolumn{1}{l|}{0.73}       & \multicolumn{1}{l|}{0.79}      & \multicolumn{1}{l|}{0.66}        &  0.78     & \multicolumn{1}{l|}{5.60E-07}     & \multicolumn{1}{l|}{3.26E-07}       & \multicolumn{1}{l|}{3.77E-07}      & \multicolumn{1}{l|}{8.42E-07}        &  3.41E-07     \\ \hline
64000  &\multicolumn{1}{|l|}{2.06}     & \multicolumn{1}{l|}{2.27}       & \multicolumn{1}{l|}{2.48}      & \multicolumn{1}{l|}{2.15}        &  2.48     & \multicolumn{1}{l|}{1.25E-06}     & \multicolumn{1}{l|}{6.64E-07}       & \multicolumn{1}{l|}{5.84E-07}      & \multicolumn{1}{l|}{1.74E-06}        &  7.03E-07     \\ \hline
125000  &\multicolumn{1}{|l|}{4.2}     & \multicolumn{1}{l|}{4.75}        & \multicolumn{1}{l|}{5.09}      & \multicolumn{1}{l|}{4.46}        &  5.08     & \multicolumn{1}{l|}{3.09E-06}     & \multicolumn{1}{l|}{8.26E-07}       & \multicolumn{1}{l|}{1.02E-06}      & \multicolumn{1}{l|}{3.32E-06}        &  1.11E-06     \\ \hline
512000 &\multicolumn{1}{|l|}{22.75}     & \multicolumn{1}{l|}{25.94}     & \multicolumn{1}{l|}{28.74}      & \multicolumn{1}{l|}{25.35}      &  28.48    & \multicolumn{1}{l|}{4.07E-06}     & \multicolumn{1}{l|}{9.99E-07}       & \multicolumn{1}{l|}{1.45E-06}      & \multicolumn{1}{l|}{6.11E-06}        &  1.34E-06     \\ \hline
1000000 &\multicolumn{1}{|l|}{43.41}     & \multicolumn{1}{l|}{51.12}     & \multicolumn{1}{l|}{-}          & \multicolumn{1}{l|}{48.45}      &   -       & \multicolumn{1}{l|}{7.19E-06}     & \multicolumn{1}{l|}{1.09E-06}       & \multicolumn{1}{l|}{-}           & \multicolumn{1}{l|}{1.74E-05}        &    -   \\ \hline
\end{tabular}
}
    \caption{Performance of the proposed and the nested algorithms as mentioned in \Cref{tab:nex_notation} in terms of the memory and relative error. We set the tolerance $\epsilon = 10^{-6}$.}
    \label{table:helm3d_b}
\end{table}

We arrange the above algorithms in ascending order with respect to memory and MVP time.\\
\textbf{Memory:} $\mathcal{H}^2_{*}$(b$+$t) < $\mathcal{H}^2_{*}$(t) < $\bkt{\mathcal{H}^2 + \mathcal{H}}_{*}$ < $\mathcal{H}^2_{\sqrt{d}}$(b) $\leq$ $\mathcal{H}^2_{\sqrt{d}}$(t) \\ 
\textbf{MVP time:} $\mathcal{H}^2_{*}$(b$+$t) < $\mathcal{H}^2_{*}$(t) < $\bkt{\mathcal{H}^2 + \mathcal{H}}_{*}$  < $\mathcal{H}^2_{\sqrt{d}}$(b) < $\mathcal{H}^2_{\sqrt{d}}$(t)

\subsubsection{Fast MVP accelerated GMRES for RBF interpolation}
Let the location of the particles $\{\pmb{x}_i \}_{i=1}^N$ be the $N^{1/3} \times N^{1/3} \times N^{1/3}$ Chebyshev grid on the domain $[-1,1]^3$. We consider the Chebyshev distribution of particles to study the performance of the nested algorithms over slightly non-uniformly distributed particles. However, we use the uniform $2^d$ tree (oct tree in $3$D) as described in \Cref{tree_construction}. 

\begin{figure}[H]
    \centering
    \subfloat[]{\includegraphics[height=3.8cm, width=5.5cm]{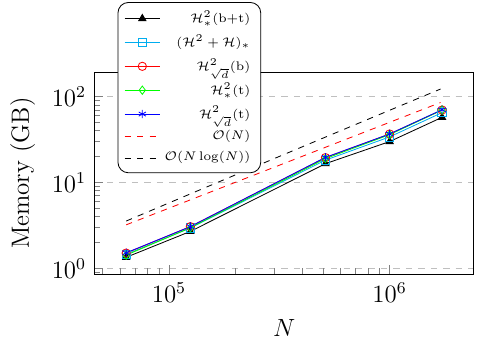}\label{rbf3d_lna_mem}}%
    \subfloat[]{\includegraphics[height=3.8cm, width=5.5cm]{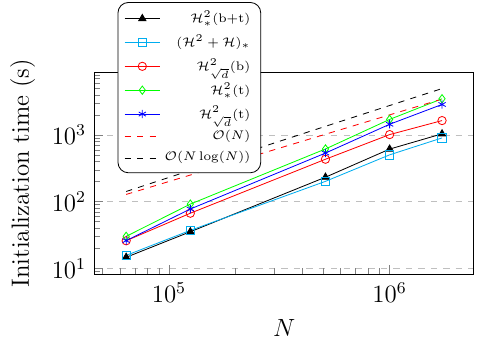}\label{rbf3d_lna_assem}}%
    \subfloat[]{\includegraphics[height=3.8cm, width=5.5cm]{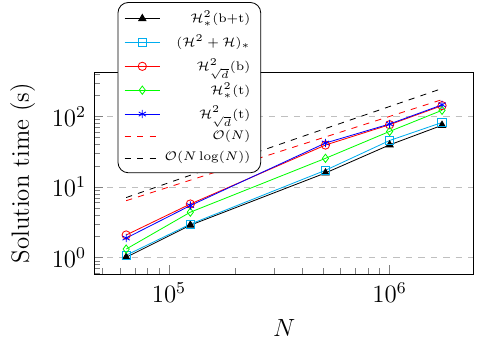}\label{rbf3d_lna_sol}}%
    \caption{Plots of Memory, Initialization time and Solution time of the fast MVP algorithms.}
    \label{fig:rbf3d_lna}
\end{figure}

Let us consider the following radial basis function
\begin{equation}
  G(r) = \begin{cases}
            \dfrac{\log{r}}{\log{a}} & \text{ if } r \geq a\\
            \dfrac{r \log{r} - 1}{a \log{a} -1} & \text{ if } r < a
            \end{cases}  
\end{equation}
Consider the linear system (\Cref{rbf3d_eq2}) generated by the above kernel function $G$
\begin{equation} \label{rbf3d_eq2}
    \alpha \lambda_i + \dsum_{j=1, j \neq i}^N G \bkt{\magn{\pmb{x}_i - \pmb{x}_j}_2} \lambda_j = b_i, \qquad i = 1,2,\dots, N.
\end{equation}
We set $a=0.0001$ and $\alpha=\sqrt{N}$. By setting $\pmb{b}$ as described in the $11^{th}$ cell of \Cref{tab:app_notations}, the \Cref{rbf3d_eq2} can be written in the form 
\begin{equation} \label{rbf3d_eq3}
    K \pmb{\lambda} = \pmb{b}
\end{equation}
The \Cref{rbf3d_eq3} is solved using \textbf{fast} GMRES. We set $\epsilon = 10^{-6}$ and $\epsilon_{GMRES} = 10^{-10}$ and tabulate the memory, initialization time and solution time in \Cref{table:rbf3d} for all the nested algorithms. The relative error in solution $(RE_{sol})$ is of order $10^{-6}$ in all cases. We also plot the memory (\Cref{rbf3d_lna_mem}), initialization time (\Cref{rbf3d_lna_assem}) and solution time (\Cref{rbf3d_lna_sol}) in \Cref{fig:rbf3d_lna}.

\begin{table}[H]
      \centering
      \resizebox{\textwidth}{!}{%
      \setlength\extrarowheight{0.9pt}
\begin{tabular}{|l|lllll|lllll|lllll|l|}
\hline
\multirow{2}{*}{N} & \multicolumn{5}{l|}{\qquad \qquad \qquad Memory (GB)}                                                                                           & \multicolumn{5}{l|}{\qquad \qquad \quad Initialization time (s)}                                                                                & \multicolumn{5}{l|}{\qquad \qquad \quad Solution time (s)}                                                                                      & \multirow{2}{*}{$\#iter$} \\ \cline{2-16}
                   & \multicolumn{1}{l|}{\scriptsize $\mathcal{H}^2_{*}$(b$+$t)} & \multicolumn{1}{l|}{\scriptsize $\bkt{\mathcal{H}^2 + \mathcal{H}}_{*}$} & \multicolumn{1}{l|}{\scriptsize $\mathcal{H}^2_{\sqrt{d}}$(b)} & \multicolumn{1}{l|}{\scriptsize $\mathcal{H}^2_{*}$(t)} & \scriptsize $\mathcal{H}^2_{\sqrt{d}}$(t) & \multicolumn{1}{l|}{\scriptsize $\mathcal{H}^2_{*}$(b$+$t)} & \multicolumn{1}{l|}{\scriptsize $\bkt{\mathcal{H}^2 + \mathcal{H}}_{*}$} & \multicolumn{1}{l|}{\scriptsize $\mathcal{H}^2_{\sqrt{d}}$(b)} & \multicolumn{1}{l|}{\scriptsize $\mathcal{H}^2_{*}$(t)} & \scriptsize $\mathcal{H}^2_{\sqrt{d}}$(t) & \multicolumn{1}{l|}{\scriptsize $\mathcal{H}^2_{*}$(b$+$t)} & \multicolumn{1}{l|}{\scriptsize $\bkt{\mathcal{H}^2 + \mathcal{H}}_{*}$} & \multicolumn{1}{l|}{\scriptsize $\mathcal{H}^2_{\sqrt{d}}$(b)} & \multicolumn{1}{l|}{\scriptsize $\mathcal{H}^2_{*}$(t)} & \scriptsize $\mathcal{H}^2_{\sqrt{d}}$(t) &                       \\ \hline 
            64000      & \multicolumn{1}{l|}{1.34}     & \multicolumn{1}{l|}{1.45}       & \multicolumn{1}{l|}{1.5}      & \multicolumn{1}{l|}{1.42}        &    1.51   & \multicolumn{1}{l|}{14.58}     & \multicolumn{1}{l|}{15.57}       & \multicolumn{1}{l|}{25.81}      & \multicolumn{1}{l|}{30.06}        &  25.93     & \multicolumn{1}{l|}{1.01}     & \multicolumn{1}{l|}{1.08}       & \multicolumn{1}{l|}{2.09}      & \multicolumn{1}{l|}{1.32}        &   1.89    &  19            \\ \hline
            125000    & \multicolumn{1}{l|}{2.68}     & \multicolumn{1}{l|}{2.94}       & \multicolumn{1}{l|}{3.05}      & \multicolumn{1}{l|}{2.96}        &    3.05   & \multicolumn{1}{l|}{35.24}     & \multicolumn{1}{l|}{36.98}       & \multicolumn{1}{l|}{66.85}      & \multicolumn{1}{l|}{91.65}        & 77.34      & \multicolumn{1}{l|}{2.88}     & \multicolumn{1}{l|}{2.98}       & \multicolumn{1}{l|}{5.74}      & \multicolumn{1}{l|}{4.4}        &   5.46    &  20                     \\ \hline
           512000    & \multicolumn{1}{l|}{16.5}     & \multicolumn{1}{l|}{18.15}       & \multicolumn{1}{l|}{19.4}      & \multicolumn{1}{l|}{18.91}        &  19.54   & \multicolumn{1}{l|}{232.01}     & \multicolumn{1}{l|}{201.11}       & \multicolumn{1}{l|}{431.82}      & \multicolumn{1}{l|}{615.79}     & 534.27      & \multicolumn{1}{l|}{15.8}     & \multicolumn{1}{l|}{17.25}       & \multicolumn{1}{l|}{39.32}      & \multicolumn{1}{l|}{25.48}        &  41.97     & 25                      \\ \hline
         1000000    & \multicolumn{1}{l|}{29.85}     & \multicolumn{1}{l|}{33.85}       & \multicolumn{1}{l|}{36.55}      & \multicolumn{1}{l|}{35.95}        &   36.60  & \multicolumn{1}{l|}{613.21}     & \multicolumn{1}{l|}{502.63}       & \multicolumn{1}{l|}{1011.81}      & \multicolumn{1}{l|}{1709.84}        & 1445.45      & \multicolumn{1}{l|}{39.51}     & \multicolumn{1}{l|}{45.58}       & \multicolumn{1}{l|}{76.31}      & \multicolumn{1}{l|}{61.29}        &  78.07     &  27                     \\ \hline
        1728000    & \multicolumn{1}{l|}{57.15}     & \multicolumn{1}{l|}{64.87}       & \multicolumn{1}{l|}{69.62}      & \multicolumn{1}{l|}{69.23}        &  70.13    & \multicolumn{1}{l|}{1038.16}     & \multicolumn{1}{l|}{901.72}       & \multicolumn{1}{l|}{1641.97}      & \multicolumn{1}{l|}{3490.25}        & 2867.08      & \multicolumn{1}{l|}{74.86}     & \multicolumn{1}{l|}{82.02}       & \multicolumn{1}{l|}{142.17}      & \multicolumn{1}{l|}{122.08}        & 145.93      &  30                       \\ \hline
\end{tabular}
}
    \caption{Performance of the proposed and the nested algorithms as mentioned in \Cref{tab:nex_notation}. We set $\epsilon = 10^{-6}$ and $\epsilon_{GMRES} = 10^{-10}$.}
    \label{table:rbf3d}
\end{table}

We arrange the above algorithms in ascending order with respect to memory and solution time.\\
\textbf{Memory:} $\mathcal{H}^2_{*}$(b$+$t) < $\bkt{\mathcal{H}^2 + \mathcal{H}}_{*}$ < $\mathcal{H}^2_{*}$(t) < $\mathcal{H}^2_{\sqrt{d}}$(b) $\leq$ $\mathcal{H}^2_{\sqrt{d}}$(t) \\ 
\textbf{Solution time:} $\mathcal{H}^2_{*}$(b$+$t) < $\bkt{\mathcal{H}^2 + \mathcal{H}}_{*}$ < $\mathcal{H}^2_{*}$(t) < $\mathcal{H}^2_{\sqrt{d}}$(b) < $\mathcal{H}^2_{\sqrt{d}}$(t)

In all the experiments, the proposed $\mathcal{H}^2_{*}$(b$+$t) algorithm outperforms the NCA-based standard $\mathcal{H}^2$ matrix algorithms \cite{zhao2019fast} in the context of memory and MVP/Solution time. Also, $\mathcal{H}^2_{*}$(b$+$t) has the least initialization time among the fully nested algorithms. Therefore, the proposed $\mathcal{H}^2_{*}$(b$+$t) could be an attractive alternative to the standard $\mathcal{H}^2$ matrix algorithms. Additionally, in $3$D, the proposed $\bkt{\mathcal{H}^2 + \mathcal{H}}_{*}$ algorithm performs slightly better than the standard $\mathcal{H}^2$ matrix algorithms.
\end{document}

%% file: include.tex
\usepackage[T1]{fontenc}
\usepackage{xr}
\usepackage[utf8]{inputenc}
\usepackage{amsfonts}
\usepackage{amsmath,amssymb}
\usepackage[caption=false]{subfig} 
\usepackage{tikz,pgfplots}
\usepackage{multirow}
\usepackage{makecell}
\usepackage{tikz-3dplot}
\usetikzlibrary{math}
\usepackage{bm}
\usepackage{rotating}
\usepackage{algorithm}
\usepackage{algpseudocode}
\usepackage[top=1.06in,bottom=1.06in,left=1in,right=1in]{geometry}
\usepackage{graphicx,epstopdf} 
\usepackage{enumerate}
\usepackage{enumitem}
\usepackage{graphicx}

\usepackage{hyperref}
\usepackage{relsize}
\usepackage{adjustbox}
\usepackage{cleveref}
\usepackage{tcolorbox}
\usepackage{dsfont}
\usepackage{accents}

\usepackage{multirow}
\usepackage{booktabs}
\usepackage{graphicx}

\usepackage{pdflscape}

\usepackage{perpage}
\MakePerPage{footnote}

\usepackage{colortbl}
\usetikzlibrary{fadings}
\ifpdf
  \DeclareGraphicsExtensions{.eps,.pdf,.png,.jpg}
\else
  \DeclareGraphicsExtensions{.eps}
\fi

\newcommand{\rk}[1]{{#1}}
\newcommand{\RK}[1]{{#1}}

\headers{New algebraic fast algorithms for $N$-body problems in two and three dimensions}{Ritesh Khan, Sivaram Ambikasaran}

\title{New algebraic fast algorithms for $N$-body problems in two and three dimensions}

\author{Ritesh Khan\thanks{Department of Mathematics, IIT Madras, Chennai, India. \email{\lowercase{khanritesh28@gmail.com}}}
\and Sivaram Ambikasaran\thanks{
Department of Data Science \& Artificial Intelligence, Wadhwani School of Data Science \& Artificial Intelligence,
Department of Mathematics and  IIT Madras, Chennai, India. \email{\lowercase{sivaambi@alumni.stanford.edu}}}}

\newsiamremark{remark}{Remark}

\hypersetup{
    colorlinks=true,
    urlcolor=blue,
    citecolor = magenta
    }
\usetikzlibrary{arrows} 
\usetikzlibrary{decorations.markings}

\pgfplotsset{compat=1.15} 

\newcommand{\dsum}{\displaystyle\sum}

\newcommand{\bkt}[1]{\left(#1\right)}

\newcommand{\magn}[1]{\left\lVert#1\right\rVert}

\newcommand{\ceil}[1]{\left\lceil#1\right\rceil}

\newcommand{\cmt}[1]{\iffalse {#1} \fi}

\makeatletter
\newcommand*{\addFileDependency}[1]{
  \typeout{(#1)}
  \@addtofilelist{#1}
  \IfFileExists{#1}{}{\typeout{No file #1.}}
}
\makeatother